\documentclass[opre,nonblindrev]{informs3_hide}
\usepackage{endnotes}
\let\enotesize=\normalsize
\def\notesname{Endnotes}%
\def\makeenmark{$^{\theenmark}$}
\def\enoteformat{\rightskip0pt\leftskip0pt\parindent=1.75em
  \leavevmode\llap{\theenmark.\enskip}}

\usepackage{natbib}
 \bibpunct[, ]{(}{)}{,}{a}{}{,}%
 \def\bibfont{\small}%
\TheoremsNumberedBySection
\ECRepeatTheorems

\EquationsNumberedThrough    % Default: (1), (2), ...
\usepackage{amssymb,amsmath,amsxtra,amstext,xfrac,mathtools}
\usepackage{graphicx,ccaption,booktabs,enumerate,paralist,mdwlist}
\usepackage{dsfont,verbatim,latexsym}
\usepackage[titletoc,toc,page]{appendix}
\usepackage{mathrsfs}
\usepackage{bbm}
\usepackage{bm}
\usepackage{multirow}
\usepackage[table]{xcolor}
\usepackage{longtable,multirow,threeparttable}
\usepackage{float}
\usepackage{tikz,qtree}
\usepackage{mwe} 
\usepackage{subfigure}
\usepackage{eurosym}
\usepackage[english]{babel}
\usepackage[utf8]{inputenc}
\usepackage{fancyhdr}
\usepackage{color}
\usepackage{hyperref}
\hypersetup{
    colorlinks=true,
    citecolor=green,
    filecolor=blue,
    linkcolor=blue,
    urlcolor=blue
}
\usepackage{multirow}
\usepackage{caption}

\usepackage{algorithm,algorithmic}

\definecolor{MyDarkBlue}{RGB}{158,0,0}

\newenvironment{proof1}{\paragraph{\textit{Proof.}}}{\hfill$\square$}

\def\EE{\mathbb{E}}
\def\PP{\mathbb{P}}
\def\RR{\mathbb{R}}
\def\NN{\mathbb{N}}
\def\eps{\epsilon}
\def\cO{\mathcal{O}}
\usepackage{accents}
\newcommand\munderbar[1]{%
  \underaccent{\bar}{#1}}
\renewcommand{\theARTICLETOP}{}

\RUNTITLE{Efficient Algorithms for Minimizing Compositions of Convex Functions and Random Functions}
\RUNAUTHOR{Chen et al.}
\TITLE{Efficient Algorithms for A Class of Stochastic Hidden Convex Optimization and Its Applications in Network Revenue Management
}

\ARTICLEAUTHORS{
\AUTHOR{Xin Chen$^1$,  Niao He$^2$, Yifan Hu$^3$, Zikun Ye$^4$}
\AFF{
$^1$ ISyE, Georgia Institute of Technology, Atlanta, US, \EMAIL{xin.chen@isye.gatech.edu}\\
$^2$ Department of Computer Science, ETH Z\"urich, Switzerland,
\EMAIL{niao.he@inf.ethz.ch}\\
$^3$ College of Management of Technology, EPFL, Switzerland, \EMAIL{yifan.hu@epfl.ch}\\
$^4$
Michael G. Foster School of Business, University of Washington, Seattle, US,
\EMAIL{zikunye@uw.edu}
\\
}
}

\begin{document}

%%%%%%%%%%%%%%%%%%%%%%%%%%%%%%%%%%%%%%%%%%%%%%%%%%%%%%%%%%%%%%%%%%%%%%%%%%%%%%%%%%%%

\OneAndAHalfSpacedXI 
\ABSTRACT{
\noindent  We study a class of stochastic nonconvex optimization in the form of $\min_{x\in\mathcal{X}} F(x):=\EE_\xi [f(\phi(x,\xi))]$, i.e., $F$ is a composition of a convex function $f$ and a random function $\phi$. Leveraging an (implicit) convex reformulation via a variable transformation $u=\EE[\phi(x,\xi)]$, we develop stochastic gradient-based algorithms and establish their sample and gradient complexities for achieving an $\eps$-global optimal solution. Interestingly, our proposed Mirror Stochastic Gradient (MSG) method operates only in the original $x$-space using gradient estimators of the original nonconvex objective $F$ and achieves $\tilde \cO(\eps^{-2})$ complexities, which matches the lower bounds for solving stochastic convex optimization problems. Under booking limits control, we formulate the air-cargo network revenue management (NRM) problem with random two-dimensional capacity, random consumption, and routing flexibility as a special case of the stochastic nonconvex optimization, where the random function $\phi(x,\xi)=x\wedge\xi$, i.e., the random demand $\xi$ truncates the booking limit decision $x$. Extensive numerical experiments demonstrate the superior performance of our proposed MSG algorithm for booking limit control with higher revenue and lower computation cost than state-of-the-art bid-price-based control policies, especially when the variance of random capacity is large.
}


\maketitle

\section{Introduction}\label{sec:intro}
A wide range of operations management problems are special cases of the following stochastic optimization model,
\begin{equation}\label{problem:math_original_extension}
\min_{x\in \mathcal{X}}~F(x) \coloneqq \mathbb{E}_{\xi\sim \PP(\xi)} [f(\phi(x,\xi))],
\end{equation}
where $\mathcal{X}  = [\munderbar{X}_1,\bar{X}_1]\times\ldots\times [\munderbar{X}_d,\bar{X}_d] \subseteq \RR^d$, $\xi\in\Xi \subseteq \RR^d$ is a random vector, $\phi(x,\xi)=(\phi_1(x_1,\xi_1),\dots,\phi_d(x_d,\xi_d))^\top$ is component-wise non-decreasing in $x$, and $f$ is convex.
Throughout the paper, we assume that the distribution $\PP(\xi)$ remains unknown, and we can only generate independent and identically distributed (i.i.d.) samples from $\PP(\xi)$.

The optimization problem \eqref{problem:math_original_extension} arises pervasively in supply chain management and revenue management. A notable example is $\phi(x,\xi)=x\wedge\xi$, where $\wedge$ denotes component-wise minimum. In inventory control problems with supply capacity uncertainty, the amount delivered by suppliers is the minimum of the replenishment order quantity $x$ and the realized random capacity $\xi$, i.e., $\phi(x,\xi)=x\wedge\xi$ \citep[see, e.g.,][]{ciarallo1994periodic,chen2015new, chen2018preservation,feng2018supply,chen2019stochastic,feng2019dynamic}, \textcolor{black}{ in network revenue management problems using booking limit control policies, the accepted reservation is the minimum of the booking limit decision $x$ and the random demand  $\xi$} \citep[see, e.g.,][]{brumelle1993airline,karaesmen2004overbooking,li2017dynamic,chen2018preservation}.

For these applications, an intrinsic challenge is that the random function $\phi(x,\xi)$ is generally nonlinear in $x$; thus, the objective function $F$ is nonconvex in $x$ even if $f$ is (strongly) convex. For example, when $\phi(x,\xi)=x\wedge \xi$ and $f(x) = \|x\|^2$, it is easy to verify that $F(x)$ is nonconvex. As a result, how to efficiently solve the nonconvex problem \eqref{problem:math_original_extension} to global optimality remains unclear. 
{\color{black} The aim of this paper is to design efficient algorithms that solve the optimization problem to global optimality. We focus on stochastic gradient-based methods as they can handle online data and are suitable for large-scale problems.
We measure the efficiency of our proposed algorithms by sample complexity and gradient complexity, i.e., the number of samples and the number of evaluations of $\nabla f$ needed to achieve an $\eps$-optimal solution, respectively. Despite the aforementioned advantages of gradient-based methods, they generally can only converge to approximate stationary points of nonconvex objectives. }

Interestingly, under some technical conditions on the random function $\phi(x,\xi)$, an equivalent convex reformulation of the nonconvex problem \eqref{problem:math_original_extension} exists~\citep{feng2018supply}:
\begin{equation}
\label{problem:transformation_short}
\min_{u\in\mathcal{U}}~G(u) \coloneqq \mathbb{E}_{\xi\sim \PP(\xi)} [f(\phi( g^{-1}(u), \xi))],
\end{equation}
where $g(x):=\mathbb{E}_{\xi\sim \PP(\xi)} [\phi(x, \xi)]$, $\mathcal{U}=g(\mathcal{X})$ is the image of $g$ on $\mathcal{X}$ and is convex, and $g^{-1}$ is the inverse of $g$ for simplicity.
With the convex reformulation, it is promising to solve the original nonconvex objective to global optimality. However, to the best of our knowledge, no algorithm has been developed to solve either the original nonconvex problem or the convex reformulation, which we address in this paper.

Note that existing gradient-based methods, like projected \textit{stochastic gradient descent} (SGD)~\citep{nemirovski2009robust}, are not directly applicable to solving the convex reformulation $\min_{u\in\mathcal{U}} G(u) = F(g^{-1}(u))$. Indeed, since $g(x) = \EE [\phi(x,\xi)]$ involves the unknown distribution $\PP(\xi)$, $g^{-1}(u)$ is unknown. As a result, it is hard to build unbiased gradient estimators for $G(u)$. For the same reason, the closed-form of $\mathcal{U}$ is unknown; hence, it is hard to perform projection onto $\mathcal{U}$.  

To address these issues, a natural idea is to utilize \emph{sample average approximation} (SAA)~\citep{kleywegt2002sample} and apply projected SGD on the empirical convex reformulation constructed via the empirical distribution. We denote such a method as SAA+SG (Algorithm \ref{alg:sgd_on_g} with detailed discussion in Appendix \ref{sec:sgd_u}). Although SAA+SG is intuitive and converges globally, it has several drawbacks. First, it requires access to offline data, limiting its applicability to the online setting. Second, the algorithm needs to estimate $g^{-1}(u)$ at each iteration, which requires solving an additional optimization subproblem and adds additional computational costs.

Instead, we consider algorithms that only operate in the original $x$ space, i.e., algorithms that do not require updating the transformed decision variable $u$. Specifically, we propose the \emph{Regularized Stochastic Gradient descent} (RSG) in Algorithm \ref{alg:RSG} and the \emph{Mirror Stochastic Gradient descent} (MSG) in Algorithm \ref{alg:MSG}. \textcolor{black}{ RSG performs regularized projected SGD updates on problem \eqref{problem:math_original_extension} and converges to a specific approximate stationary point. Under mild conditions, we show that the converging point corresponds to an approximate global optimal solution. MSG performs an update on the original $x$ space that mirrors a virtual SGD update on the convex reformation problem \eqref{eqn:trans-model-2-extension} and finds an approximate global optimal solution of the convex reformulation \eqref{problem:transformation_short} but only updates $u$ implicitly via updating $x$. To achieve such mirror behavior, MSG uses $\cO(\log(\eps^{-1}))$ number of samples at each iteration to build a preconditioning matrix (see \eqref{eq:MSG_desired_update}).}

\begin{table}[t]\centering 
\caption{Summary of complexities of proposed algorithms for achieving an $\eps$-optimal solution} 
\label{tab:algos_complx} %\begin{tabular}
    \begin{tabular}{cccc}
    \hline
    \hline\\[-1.8ex]
    {Algorithm}  & {Properties} &{Sample Complexity} & {Gradient Complexity}\\
    \hline\\[-1.8ex]

    {SAA+SG }   &  Requires a batch of offline data & \multirow{2}{*}{$\tilde\cO(d\eps^{-2})$ }   & \multirow{2}{*}{$\tilde\cO(d^2\eps^{-4})$}        \\
    {(Algorithm \ref{alg:sgd_on_g})}  & requires less assumption on $\PP(\xi)$ & & \\ 
    \\[-1.8ex]\hline\\[-1.8ex]
    
    RSG   &  simple to implement & \multirow{2}{*}{$\cO(\eps^{-4})$ }   & \multirow{2}{*}{$\cO(\eps^{-4})$}         \\
    (Algorithm \ref{alg:RSG})  & one sample per-iteration & &  \\ 
    \\[-1.8ex]\hline\\[-1.8ex]
    
    MSG   &  $\cO(\log(\eps^{-1}))$ samples per-iteration & \multirow{2}{*}{$\tilde \cO(\eps^{-2})$}    & \multirow{2}{*}{$\tilde \cO(\eps^{-2})$}        \\
    (Algorithm \ref{alg:MSG})  & optimal sample \& gradient complexity & &  \\
    \\[-1.8ex]\hline
    \hline\\[-1.8ex]
    \end{tabular}
    \vskip -0.1in
    \end{table}
    
We establish the global convergence of SAA+SG, RSG, and MSG. Table \ref{tab:algos_complx} summarizes the sample and gradient complexities. {For a more comprehensive summary of assumptions and complexity bounds, see Table \ref{tab:algos_assumption} in Appendix \ref{appendix:discussion_on_assupmtions}. } RSG achieves a $\cO(\eps^{-4})$ complexity bound while MSG  achieves better $\tilde \cO(\eps^{-2})$ sample and gradient complexities, where $\cO(\cdot)$ hides the constant terms and $\tilde \cO(\cdot)$ additionally hides the logarithmic dependence. In terms of lower bounds for solving the original problem \eqref{problem:math_original_extension}, utilizing the analysis of lower bounds on black-box stochastic gradient methods for stochastic convex optimization~\citep{agarwal2009information}, we obtain a $\cO(\eps^{-2})$ lower bounds for problem \eqref{problem:math_original_extension}. It implies that the performance of MSG matches the optimal possible black-box stochastic gradient method in terms of the dependence on accuracy if ignoring the logarithmic factor.

We apply the proposed methods to solve an air-cargo network revenue management (NRM) problem with booking limit control, as well as a passenger NRM problem as a special case of the air-cargo NRM problem. We first formulate the NRM problem as two-stage stochastic programming and show that it is a special case of the nonconvex stochastic problem \eqref{problem:math_original_extension}. Specifically, during the first stage, we first set up booking limits and then accept the amount of demand up to booking limits $x$ truncated by uncertain demand $\xi$, i.e., $\phi(x,\xi) = x\wedge \xi$. In the second stage, we make capacity allocation and routing decisions to serve the accepted demand after the reveal of the random capacity. A notable advantage of such modeling is that it only requires the aggregated level of demand over the whole reservation period rather than the arrival rate in each period, which is typically assumed to be accessible in dynamic models. {Note that the focus of the paper is to find an optimal booking limit control. Thus we reduce the usual online NRM problem to a stochastic optimization problem, and we are interested in analyzing the sample efficiency for solving the problem. We leave the design of the online booking limit control policy as an interesting future direction. }

We further conduct extensive numerical experiments to demonstrate the effectiveness and generalizability of the booking limit control in NRM problems, comparing to several different bid-price based control policies, including deterministic linear programming (DLP), dynamic programmings decomposition (DPD) \citep{erdelyi2010dynamic}, and the state-of-the-art virtual capacity and bid price control (VCBP) \citep{previgliano2021managing} in passenger NRM problem, and the state-of-art DPD \citep{barz2016air} specifically designed for air-cargo NRM problem. 

\subsection{Contributions}
%Our contributions are summarized as follows.

{\color{black}
\noindent\textbf{Algorithm design and global convergence with non-asymptotic guarantees.}  {We propose three algorithms, SAA+SG, RSG, and MSG, and establish the first non-asymptotic global convergence guarantees for the nonconvex stochastic optimization \eqref{problem:math_original_extension}. } In addition, RSG and MSG operate only in the original $x$-space, and MSG achieves the optimal complexity bounds. 
}

\noindent \textbf{NRM modeling and algorithm.} 
To the best of our knowledge, we are the first in the literature to propose booking limit control for the air-cargo network revenue problem that takes into account random show-ups, random capacity, random consumption, and routing flexibility at the same time. 
Our algorithms provide non-asymptotic global guarantees under some mild assumptions, while the VCBP algorithm only converges asymptotically to stationary points. 

\noindent \textbf{Numerical Results.}
Our numerical results demonstrate the superior performance of the proposed algorithms and provide strong justification for utilizing booking limit control in these NRM problems. In passenger NRM (a special case of air-cargo NRM), the booking limit control policy obtained by MSG significantly outperforms bid-price-based methods DLP, DPD, and VCBP with 43.6\%, 8.3\%, and 4.8\% revenue improvement, respectively, and achieves the lowest computation time.
In air-cargo NRM, our method outperforms the state-of-the-art DPD \citep{barz2016air} by 12.86\% under the fixed-route setting and 17.22\% under the routing flexibility setting, which indicates the advantage of the booking limit control policy in dealing with routing flexibility. In addition, the numerical results indicate that booking limit control gains more revenue improvement against bid-price-based control policies, especially when the random capacity has a large variance.

\subsection{Literature Review}\label{sec:literature}
We next review three streams of related literature.

\subsubsection{Stochastic Gradient-Based Algorithms.} (Projected) SGD and its numerous variants form one of the most important families of algorithms for solving classical stochastic optimization. For strongly convex and convex stochastic optimization \citep{nemirovski2009robust, bottou2018optimization},  the complexity to achieve an $\epsilon$-global optimality is $\mathcal{O}(\epsilon^{-1})$ and $\mathcal{O}(\epsilon^{-2})$ , respectively (For SGD, sample complexity equals gradient complexity). For nonconvex stochastic optimization, the gradient complexity to achieve an $\epsilon$-stationary point is $\mathcal{O}(\epsilon^{-4})$~\citep{ghadimi2013stochastic, ghadimi2016mini}. An extension of the stochastic gradient method is the stochastic primal-dual method~\citep{agrawal2014fast,li2022online}. They are usually designed to handle functional constraints that do not admit easy projection.

\subsubsection{Solving Nonconvex Optimization to Global Optimality.}
For nonconvex optimization, there are several conditions that allow design efficient algorithms with global optimality guarantees, 1) hidden convexity, i.e., the problem admits a convex reformulation~\citep{ben1996hidden}, 2) bisection methods for low-dimensional nonconvex problems, 3) Polyak-Łojasiewicz (PL) condition~\citep{karimi2016linear}, 4) structured nonconvex optimization~\citep{sun2023nonconvex}.

For hidden convex optimization problems, 
\citet{chen2019network,miao2021network}  considered a pricing-based NRM problem with a nonconvex objective and nonconvex constraints that admits a convex reformulation. Due to the nonconvex constraint, their algorithm performs updates on the space of the convex reformulation.  \citet{chen2019network} achieved $\cO(\eps^{-5})$ sample complexity in terms of the number of demands and \citet{miao2021network} improved the complexity to $\cO(\eps^{-2})$ via a sophisticated ellipsoid method with cutting planes. Differently, our problem has a convex box constraint, and thus, our algorithms operate in the original space. The proposed MSG method achieves $\tilde \cO(\eps^{-2})$ sample complexity. \citet{chen2022learning} considered a lost-sale inventory control problem with random supply in the limiting regime, which becomes a special case of problem \eqref{problem:math_original_extension} when the dimension $d=1$. To achieve a global solution, they used a bisection method rather than leveraging the hidden convexity. Note that the complexity of the bisection method scales exponentially in $d$, and thus, it is not suitable when $d$ is large. {There are other hidden convex optimization that arise from linear quadratic control~\citep{anderson2019system}, policy optimization in constrained Markov decision process~\citep{ying2023policy} and convex reinforcement learning~\citep{zhang2020variational}. After the initial version of this paper was released, several follow-up papers considered general stochastic hidden convex optimization~\citep{fatkhullin2023stochastic}, which focuses on the performance of classical SGD methods following a different analysis of our RSG method, and hidden monotone games~\citep{sakos2024exploiting}, which use a similar preconditioning idea as our proposed MSG method. Note that  \citep{sakos2024exploiting} can only achieve the same complexity bound in the unconstrained setting using the true inverse of the gradient of the transformation, which is much more costly than our MSG method. \citet{ghai2022non} explored the algorithmic equivalence between the SGD on the original nonconvex optimization and the mirror descent on the convex reformulation. However, their results have very restrictive requirements  on the transformation function that it has to be the gradient of a distance generating function of the Bregman divergence. \citet{miao2024demand} simplifies the sophisticated ellipsoid method in \citet{miao2021network} to a primal-dual based approach with $\cO(\eps^{-2})$ complexity bounds for the pricing-based NRM problem.

\citet{kunnumkal2008using} studied finding an optimal base-stock policy in an inventory system with lost sales. 
They demonstrated the asymptotic stationary convergence of a stochastic approximation method and established the relationship between an approximate stationary point and an approximate global solution. The analysis of the RSG algorithm follows a similar idea, yet we characterize the non-asymptotic sample complexity.

\citet{balseiro2023best} designed a primal-dual method to solve a nonconvex online resource allocation problem to global optimality. They address the nonconvexity via the Shapley-Folkman Theorem~\citep{starr1969quasi} that the primal-dual gap can be upper bounded by a constant that is independent of iterations. \citet{han2020optimal} designed an elimination method to achieve global solutions in nonconvex auctions. \citet{yuan2021marrying} combined SGD with bandit algorithm to search for optimal $(s,S)$ policy in the nonconvex inventory systems with fixed costs. The key to overcoming nonconvexity is that the objective is convex in $S$ while nonconvex in $s$. Thus one could discretize the $s$ space and adapt a bandit algorithm to find approximate optimality for each $S$.  

Recently, various papers studied the Polyak-Łojasiewicz (PL) condition and other error bound conditions~\citep{karimi2016linear}. These conditions ensure that first-order stationary points are also globally optimal, and thus one could utilize first-order methods to find global optimality despite nonconvexity. The sample complexity of SGD to achieve an $\eps$-global optimal solution is $\cO(\eps^{-1})$~\citep{hu2021bias}. Note that policy optimization with certain policy parameterization for reinforcement learning problems satisfies the PL condition~\citep{bhandari2019global,agarwal2020optimality}. However, one can easily verify that the PL condition does not hold for \eqref{problem:math_original_extension} when $\phi(x,\xi)=x\wedge\xi$. 

For more structured nonconvex optimization problems that admit efficient algorithms with global optimality guarantees, we refer interested readers to the website \citep{sun2023nonconvex}. Although our problem belongs to the hidden convex problem class, the transformation function in our problem is unknown, and thus the methodology developed therein is generally not applicable.

\subsubsection{Network Revenue Management.}

One popular approach to network revenue management problems is booking limit control. It sets a threshold for each reservation class and accepts all requests until the threshold is met. \citet{karaesmen2004overbooking} solved a two-stage stochastic model via SGD to obtain an overbooking limit in the setting where the demand is assumed to be infinity (i.e., no truncation), and they demonstrate asymptotic convergence of the algorithm. \citet{wang2016optimal} obtained integral booking limits from a one-period stochastic integer programming considering discrete random demands with the truncation and discrete random resource capacities. However, their focus is on the integral decision space and does not consider routing flexibility as we do in the paper. \citet{wang2021two} modeled the problem as two-stage stochastic programming under the special case when there is only one-dimensional deterministic capacity and one fixed route. In contrast, our models with the booking limit control for network revenue problems can incorporate the multi-dimensional random demand and capacity and allow flexibility in routing. Furthermore, the proposed algorithms are readily applicable and have non-asymptotic global convergence guarantees. 

In addition to booking limit control, another pervasively applied approach uses bid price control, which can be derived from deterministic linear programming (DLP) \citep{talluri1998analysis}. One can treat bid prices as prices for the resources, and the reservation is accepted if its revenue is higher than the sum of the bid prices of the required resources. More sophisticated time-dependent bid price control can be obtained from dynamic programming. \cite{erdelyi2010dynamic} propose a DPD approach to jointly make overbooking and capacity allocation decisions in passenger revenue management, but they do not consider the random capacity. Compared with passenger revenue management, air-cargo network revenue management problems receive significantly less attention in the literature because of their complication, which prevents the direct application of existing techniques developed in the passenger NRM problems. \cite{barz2016air} consider an air-cargo network setting with both random capacities and routing flexibility. They develop a DPD approach to obtain the bid price policy, which depends on the time and the expected consumption of total accepted requests. However, they only deal with the routing flexibility in a heuristic way, while our approach considers optimal routing decisions after the realization of the random demand and capacity.

\subsection*{Organizations}

The rest of the paper is organized as follows. In Section \ref{sec:alg}, we discuss the convex reformulation and the intuition behind RSG and MSG. In Section \ref{sec:global_convergence}, we demonstrate the sample and gradient complexities for RSG and MSG. In Section \ref{sec:application}, {We discuss a number of operations management applications to illustrate the broad applicability of our algorithms.}. We further formulate the NRM problem as a two-stage stochastic model, a special case of the studied stochastic nonconvex optimization.  In Section \ref{sec:numerics}, we present numerical experiments in various NRM settings.

\subsection*{Preliminaries}
\textcolor{black}{For an abuse of notation, let $\nabla$ denote derivative, (sub)gradient, {Clarke subdifferential}, and Jacobian. 
For $x,~u\in\RR^d$, and $\xi\in\Xi\subseteq \RR^d$, let $x=(x_1,\ldots,x_d)^\top$, $u=(u_1,\ldots,u_d)^\top$, $\xi=(\xi_1,\ldots,\xi_d)^\top$, where a subscript denotes the corresponding coordinate of a vector.  
We use \textcolor{black}{$\|\cdot\|$ to denote $l_2$ norm for vector and matrix. Note that the $l_2$ norm for a matrix is also known as the spectral norm, i.e., the largest singular value of a matrix. In addition, it holds that 
$\|\Lambda_1\Lambda_2\|\leq \|\Lambda_1\|\|\Lambda_2\|$ for any $\Lambda_1, \Lambda_2\in\RR^{d\times d}$.} Let $\Pi_\mathcal{X}(x):= \argmin_{y\in\mathcal{X}}\|y-x\|^2$ denote projection from $x$ onto set $\mathcal{X}$.
A function $f$ is $L_f$-Lipschitz continuous on $\mathcal{X}$ if it holds that $\|f(x)-f(y)\|\leq L_f\|x-y\|$ for any $x,y\in\mathcal{X}$. If the gradient of a function is Lipschitz continuous, we also call this function smooth. If a function $f$ satisfies $f(x)-f(y)-\nabla f(y)^\top (x-y)\geq \mu \|x-y\|^2$ for some constant $\mu$, we say $f$ is $\mu$-strongly convex if $\mu>0$, $f$ is convex if $\mu=0$, and $f$ is $\mu$-weakly convex if $\mu<0$. Note that any $S_f$-smooth function is also $\mu$-weakly convex by definition.
We use $[N] :=\{1,\dots,N\},~N\in\mathbb{N}^+$ to denote the set of subscript.} We use $\Lambda^{-\top}$ to denote the transpose of the inverse of a matrix $\Lambda$. We mainly focus on the complexity bounds in terms of the accuracy $\eps$: we use $\cO$ to hide constants that do not depend on the desired accuracy $\eps$ and use $\tilde \cO$ to further hide the $\log(\eps^{-1})$ term.

\section{Convex Reformulation and Algorithmic Design}
\label{sec:alg}

In this section, we first formally state the convex reformulation of the optimization problem \eqref{problem:math_original_extension} and the corresponding conditions. Then, we discuss the intuition behind the algorithmic design of our proposed gradient-based methods. 
Recall the transformed problem:
\begin{equation}\label{eqn:trans-model-2-extension}
\min_{u\in\mathcal{U}}~G(u) \coloneqq \mathbb{E}_{\xi\sim \PP(\xi)} [f(\phi( g^{-1}(u), \xi))],
\end{equation}
where $g(x):=\mathbb{E}_{\xi\sim \PP(\xi)} [\phi(x, \xi)]$, $\mathcal{U}:=\{u\mid  u_i\in\mathcal{U}_i:=[\mathbb{E}[\phi_i(\munderbar{X}_i, \xi_i)], \mathbb{E}[\phi_i(\bar{X}_i, \xi_i)]], \text{ for } i\in[d]\}$, and $g^{-1}(u):=(g^{-1}_1(u_1),\dots,g^{-1}_d(u_d))$ with $g^{-1}_i(u_i)=\inf_{x_i\in[ \munderbar{X}_i,\bar{X}_i]} \{x_i\mid g_i(x_i)\geq u_i\} \text{ for } i\in[d]$. {\color{black}Next, we list conditions for problem \eqref{eqn:trans-model-2-extension} to be an equivalent convex reformulation of problem \eqref{problem:math_original_extension}.
\begin{assumption}
\label{assumption:reformulation}
We assume
\begin{itemize}
    \item[(a).] Random vector $\xi\in\Xi\subseteq\RR^d$ is coordinate-wise independent.
    \item[(b).] Function $\phi_i(x_i,\xi_i)$ is non-decreasing in $x_i$ for any given $\xi_i\in\Xi_i$ and any $i\in[d]$.
    \item[(c).] Function $\{\phi_i(g^{-1}_i(u_i),\xi_i), u_i\in\mathcal{U}_i\}$ is \textit{stochastic linear in midpoint}\footnote{Definition 1 in \cite{feng2018supply}: A function $\{Y(x),x\in\mathcal{X}\}$ for some convex $\mathcal{X}$ is \textit{stochastically linear in midpoint} if, for any $x_1,x_2\in \mathcal{X}$, there exist $\hat{Y}(x_1)$ and $\hat{Y}(x_2)$ defined on a common probability space such that (i) $\hat{Y}(x_i)=^{d} Y(x_i),~i=1,2$ and (ii) $(\hat{Y}(x_1)+\hat{Y}(x_2))/2\leq_{cv} Y((x_1+x_2)/2)$ where $=^{d}$ denotes equal in distribution and $\leq_{cv}$ denotes concave order.} for any $i\in[d]$.
\end{itemize}
\end{assumption} 
\cite{feng2018supply} showed that stochastic linearity in midpoint property holds for various functions $\phi(x,\xi)$ used in supply chain management applications with dimension $d=1$. Below we list four examples of function $\phi$ commonly used in operations management, including $\phi(x,\xi)=x\wedge \xi$ in our NRM applications.

\begin{itemize}
\item[(i)] $\phi_i(x_i,\xi_i)=x_i\xi_i$, $\quad\quad\quad\quad\quad$ (ii) $\phi_i(x_i,\xi_i)=x_i\xi_i/(x_i+\alpha\xi_i^\kappa)$ for $\kappa\leq 1$, $\alpha>0$, $i\in[d]$,
\item[(iii)] $\phi_i(x_i,\xi_i)=x_i\wedge\xi_i$,
$\quad\quad\quad~~$ (iv) $\phi_i(x_i,\xi_i)=(x_i/(x_i+\xi_i))k$ for some $k\geq0$ and $i\in[d]$.
\end{itemize}
Here $\xi\in\RR_+^d$ is a non-negative random vector, and $x\in \RR_+^d$; thus, $\phi_i$ is  non-decreasing in $x_i$. {Example (i) appears in inventory problems with random yield. Example (ii) appears in a supply function in procurement from multiple suppliers \citep{dada2007newsvendor}. Example (iii) is the mostly studied random function. Example (iv) appears in the random supply from one producer with total production quantity $k$ to multiple firms \citep{tang2014pay}. Firm $i$ orders quantity $x_i$, and other firms order $\xi_i$ in total, which is unobserved to firm $i$. So the proportional delivery quantity to firm $i$ is $(x_i/(x_i+\xi_i))k$. We elaborate on more applications satisfying these
$\phi$ functions in Section \ref{sec:examples}.}

\begin{proposition} [\citealt{feng2018supply}] \label{thm:equivalent-trans-extension} 
Under Assumption \ref{assumption:reformulation}(a)(b), problem \eqref{eqn:trans-model-2-extension} has the same objective value as problem \eqref{problem:math_original_extension} via the variable change, i.e., $F(x)=G(g(x)),~\forall x\in \mathcal{X}$ and $G(u)=F(g^{-1}(u)),~\forall u\in \mathcal{U}$. Additionally, if Assumption \ref{assumption:reformulation}(c) holds and $f$ is convex (component-wise convex) in $x\in\mathcal{X}$, then $G$ is convex (component-wise convex) in $u\in\mathcal{U}$.
\end{proposition}
The proposition shows that for convex $f$, the reformulated problem $\min_{u\in\mathcal{U}} G(u)$ is a convex optimization problem under certain conditions. \citet{feng2018supply} demonstrated the proof of the proposition when dimension $d=1$. Since the random vector $\xi$ is component-wise independent, the proof of the one-dimensional case can be extended to the high-dimensional setting, following Theorem 7.A.8 and Theorem 7.A.24 in \cite{shaked2007stochastic} for convex $f$ and component-wise convex $f$, respectively. In addition, we make the following technical assumptions. }

\begin{assumption}
\label{assumption:general} We assume
\begin{enumerate}
    \item [(a)] Domain $\mathcal{X}$ has a finite radius $D_\mathcal{X}$, i.e., $\forall x\in\mathcal{X}, \|x\|\leq D_\mathcal{X}$.
    \item[(b)] Function $f$ is convex, $L_f$-Lipschitz continuous, and continuously differentiable.
    \item[(c)] Random function $\phi(x,\xi)$ is $L_\phi$-Lipschitz continuous in $x$ for any $\xi\in \Xi$.
    \item[(d)] {\color{black} For any $\xi\in\Xi$, random function $\phi(x,\xi)$ is differentiable in $x$ almost surely.}
\end{enumerate}
\end{assumption} 
Assumption \ref{assumption:general}(a), that domain $\mathcal{X}$ is bounded, is widely seen in supply chain management and revenue management. 
Assumption \ref{assumption:general}(b) about convexity is necessary for the convex reformulation \eqref{eqn:trans-model-2-extension}. For our NRM applications in Section \ref{sec:application}, as we will show in Lemma \ref{lem:compo-convexity}, the function $f$ is convex if all accepted demands show up and is component-wise convex otherwise. The assumption that $f$ is $L_f$-Lipschitz continuous and continuously differentiable is standard. One can easily verify that all four widely-used $\phi$ functions mentioned above satisfy Assumption \ref{assumption:general}(c)(d) under mild conditions, e.g., when $x>0$ and $\xi>0$ is a nonnegative random vector. Below we list a key assumption on the transformation function $g$. Various combinations of $\phi$ and $\PP(\xi)$ can guarantee it. 

\begin{assumption}
\label{assumption:general_2} For the transformation function $g:\mathcal{X}\rightarrow\mathcal{U}$, we assume 
\begin{enumerate}
    \item[(a)] {Function $g(x)$ is continuously differentiable for any $x\in\mathcal{X}$.}
    \item[(b)] Matrix $\nabla g(x)-\mu_g I$ is positive semi-definite for any $x\in\mathcal{X}$ and some constant $\mu_g>0$.
    \item[(c)]  Jacobian matrix  $\nabla g(x)$ is $S_g$-Lipschitz continuous in $x\in\mathcal{X}$, i.e., $\|\nabla g(x) - \nabla g(y)\|\leq S_g\|x-y\|$ for any $x,y\in\mathcal{X}$.
\end{enumerate}
\end{assumption}

We show in Appendix \ref{appendix:verify_all_phi} the general conditions to ensure Assumption \ref{assumption:general_2}. Further, Table \ref{tab:phi_assumption} in Appendix \ref{appendix:verify_all_phi} summarizes the conditions needed for all four $\phi$ functions and $\PP(\xi)$ to ensure  Assumption \ref{assumption:general_2}. For $\phi_i(x_i,\xi_i)=x_i\xi_i$, $\phi_i(x_i,\xi_i)=x_i\xi_i/(x_i+\alpha\xi_i^\kappa)$, and $\phi_i(x_i,\xi_i)=(x_i/(x_i+\xi_i))k$, they satisfy Assumption \ref{assumption:general_2} for a compact domain $\mathcal{X}\subset\RR^d_+$ and any distribution $\PP(\xi)$ that admits a nonnegative bounded support, which is common in operation management literature. For $\phi(x,\xi)=x\wedge\xi$, we characterize the conditions on $\PP(\xi)$ to ensure Assumption \ref{assumption:general_2} in Lemma \ref{lemma:verify_x_wedge_xi}. In Section \ref{sec:convergence_NRM},  we further characterize the performance of the proposed algorithm for $\phi(x,\xi)=x\wedge\xi$ when needed conditions on $\PP(\xi)$ do not hold, with the analysis given in Appendix \ref{appendix:a2.3_fails}. {Given Assumption \ref{assumption:general_2}(a)(b), the inverse function theorem ensures that $g^{-1}(x)$ is differentiable.}

Next, we provide closed forms of the gradients of $F$ and $G$. The proof is in Appendix \ref{proof_of_closed_form_gradient}.
\begin{lemma}
\label{lemma:gradient_advanced}
Under Assumptions \ref{assumption:general}(b)(c)(d) and \ref{assumption:general_2}{(a)(b)}, for any $x\in\mathcal{X}$ and any $u\in \mathcal{U}$, we have
\begin{align}
\label{eq:gradient_F}
&    \nabla F(x) 
= 
    \EE_{\xi} [\nabla \phi(x,\xi)^\top\nabla f(\phi(x, \xi))],\\ 
&    \nabla G(u) 
= 
    [\nabla g(g^{-1}(u))]^{-\top}\EE_{\xi} [\nabla \phi(g^{-1}(u),\xi)^\top\nabla f(\phi(g^{-1}(u), \xi)) ].
\end{align}
\end{lemma}

\subsection{Algorithmic Design of Global Converging Algorithms }
\label{sec:RSG}
In this subsection, we discuss the motivation for the global converging algorithm design. 

\subsubsection{Intuition of SAA+SG (Algorithm \ref{alg:sgd_on_g})}
Since Proposition \ref{thm:equivalent-trans-extension} provides an equivalent finite-dimensional convex reformulation \eqref{eqn:trans-model-2-extension} of the original nonconvex problem \eqref{problem:math_original_extension}, intuitively, one may design gradient-based methods on $G$ to solve \eqref{eqn:trans-model-2-extension}. 
A straightforward way is to perform projected stochastic gradient descent (PSGD)~\citep{nemirovski2009robust} on the convex reformulation, i.e., 
$
u^{t+1} = \Pi_{\mathcal{U}}(u^t - \gamma  v(u^t)),
$
where $v(u^t) := [\nabla g(g^{-1}(u^t))]^{-\top} \nabla \phi(g^{-1}(u^t),\xi^t)^\top\nabla f(\phi(g^{-1}(u^t), \xi^t))$  is an unbiased gradient estimator of $G(u^t)$ and $\xi^t$ is drawn independently from $\PP(\xi)$.  

However, since $\PP(\xi)$ is unknown, the closed-forms of  $g(x)=\EE [\phi(x,\xi)]$, $\mathcal{U}$, and $g^{-1}$ remain unknown. It leads to two challenges: 1) it is hard to construct unbiased stochastic gradients of $G$ since we do not know $x^t$; 2) it is hard to perform projections onto $\mathcal{U}$. Thus the classical PSGD is not implementable on $G$. 

We can utilize SAA to estimate the unknown $\PP(\xi)$ with the empirical distribution. {Based on this idea, we design a SAA+SG algorithm. To the best of our knowledge, it is the first of its kind in the literature.} Due to the page limit, we defer the detailed algorithmic construction, global convergence 
complexities, and related discussion in Appendix \ref{sec:sgd_u}. {Note that SAA+SG only requires Assumptions \ref{assumption:reformulation} and \ref{assumption:general} to achieve a global convergence. Thus it requires fewer assumptions. However, SAA+SG requires access to a batch of samples in the beginning {to build an estimated transformation function $\hat h$} and requires computing the inverse of the estimated transformation function at each iteration, which can be costly.} 

In the following subsections, we propose two algorithms, regularized stochastic gradient method (RSG) and mirror stochastic gradient method (MSG), to solve problem \eqref{problem:math_original_extension}. A key property of RSG and MSG is that both algorithms operate only in the original space on $x$, thus avoiding the indirect estimation of $x$ from $u$. {Figure \ref{fig:demo} compares the difference in the updating procedure in RSG and MSG compared to SAA+SG. It remains to figure out what is the gradien estimator $v(x)$ for RSG and MSG.}
%Figure \ref{fig:RSG} and Figure \ref{fig:SAASGD} in Appendix \ref{appendix:algorithmic_difference} illustrate the difference in the updating procedure in RSG and MSG compared to SAA+SG.

\begin{figure}[t]
    \centering
    \includegraphics[width=0.45\linewidth]{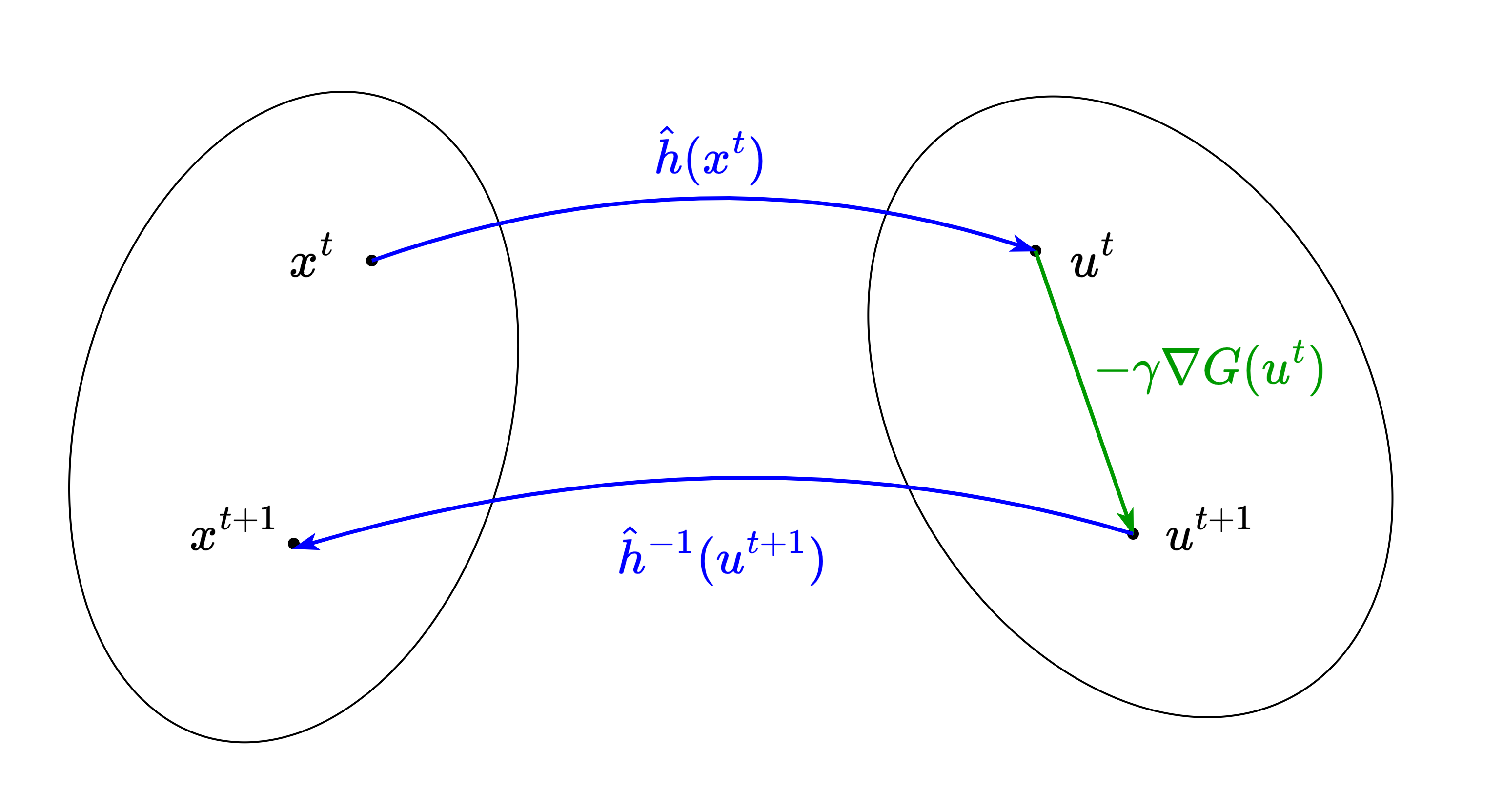}~
    \includegraphics[width=0.45\linewidth]{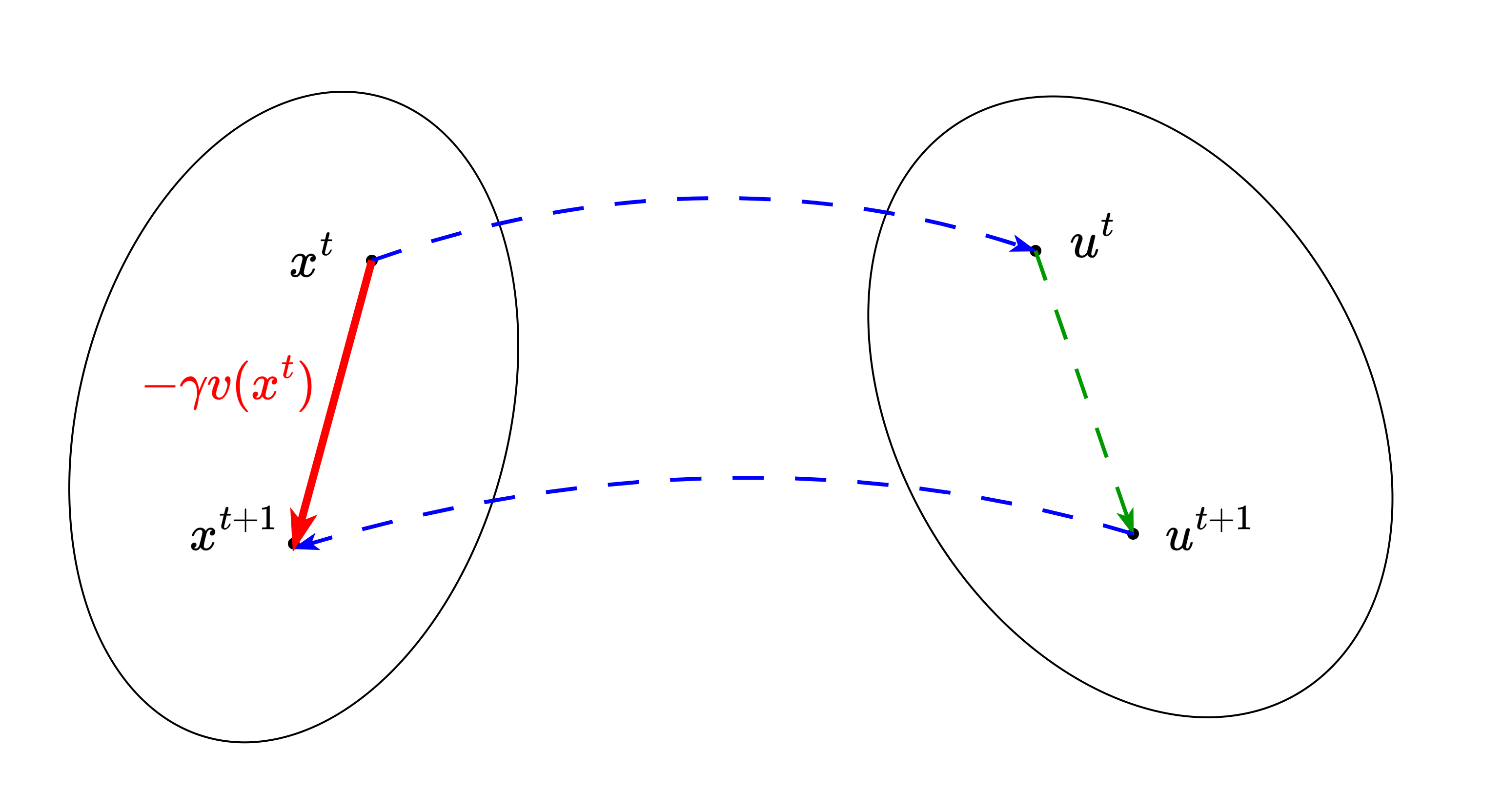}
    \caption{{Left: updating procedure of SAA+SG. Right: updating procedure of RSG and MSG.}}
    \label{fig:demo}
    %\vskip -0.1in
\end{figure}

\subsubsection{Intuition of RSG (Algorithm \ref{alg:RSG})}
By Lemma \ref{lemma:gradient_advanced}, it holds for $x = g^{-1}(u)$  that
$
\nabla G(u) =  [\nabla g(x)]^{-\top} \nabla F(x).
$
Let $x^*\in \argmin_{x\in\mathcal{X}} F(x)$ and $u^*: = g(x^*)$. By Proposition \ref{thm:equivalent-trans-extension}, we have $u^*\in\argmin_{u\in\mathcal{U}} G(u)$. Utilizing the convexity of $G$, we have
\begin{equation}
\label{intuition:RSG}
\begin{aligned}
    F(x)-F(x^*) 
=&
    G(u) -G(u^*)
\leq 
    \nabla G(u)^\top ({u}-u^*) 
\leq 
    \|\nabla G({u})\|~\|{u}-u^*\|\\
= & 
    \|[\nabla g(x)]^{-\top}\nabla F({x})\|~\|{u}-u^*\|
\leq 
    \| [\nabla g(x)]^{-1}\|~\|\nabla F({x})\|~ \|{u}-u^*\|,
\end{aligned}
\end{equation}
where the first inequality uses convexity of $G$, the second inequality uses the Cauchy-Schwarz inequality, the second equality holds by the relationship between $\nabla F(x)$ and $\nabla G(u)$, and the third inequality uses the property of spectral norm.
It implies that if $x\in\mathcal{X}$ is a stationary point of $F$ such that $\nabla F(x)=0$ and $\|[\nabla g(x)]^{-1}\|$ is finite-valued,  $x$ is also a global optimal solution.

To find an approximate stationary point of the original problem, we propose to solve problem \eqref{problem:math_original_extension} via regularized stochastic gradient method (RSG), 
$$
v_\lambda(x)=\nabla \phi(x,\xi)^\top\nabla f(\phi(x, \xi)) + \lambda x,
$$
which is an unbiased gradient estimator of a regularized objective $
F^\lambda(x) = F(x)+\frac{\lambda}{2}\|x\|^2,
$ with a regularization parameter $\lambda\geq 0$.  
Intuitively, any approximate stationary points of $F^\lambda$ are approximate stationary points of $F$ for small $\lambda$. Next, we use an example with $\phi(x,\xi) = x\wedge \xi$ to illustrate why we add this regularization.
\begin{algorithm}[t]%\footnotesize
	\caption{Regularized Stochastic Gradient (RSG)}
	\label{alg:RSG}
	\begin{algorithmic}[1]
		\REQUIRE  Number of iterations $ T $, stepsizes $\{\gamma_t\}_{t=1}^{T}$, initialization $x^1$, regularization $\lambda$.
		\FOR{$t=1$ to $T$ \do} 
		\STATE Draw a sample $\xi^t$ from $\mathbb{P}(\xi)$ and construct a gradient estimator
		\begin{center}
		$	
		v_\lambda(x^t)=\nabla \phi(x^t,\xi^t)^\top \nabla f(\phi(x^t,\xi^t)) + \lambda x^t.
		$   
		\end{center}
		\STATE Update $ x^{t+1}=\Pi_{\mathcal{X}}(x^t-\gamma_t v_\lambda(x^t))$.
		\ENDFOR
		\ENSURE $\hat x^T$ is selected uniformly from $ \{x^t\}_{t=1}^T$.
	\end{algorithmic}
\end{algorithm}

\begin{example}[Example of RSG on $\phi(x,\xi)=x\wedge\xi$]
\label{example:RSG_on_truncation}
When $\phi(x,\xi)=x\wedge\xi$, an unbiased gradient estimator of $F$ is $v(x) := \nabla \phi(x,\xi)^\top\nabla f(\phi(x, \xi))=\mathbb{I}(\xi\geq x)^\top\nabla f(x\wedge \xi)$, where $\mathbb{I}(\xi\geq x)$ denotes a diagonal matrix with the $i$-th diagonal entry being the indicator function $\mathbf{1}(\xi_i\geq x_i)$. If $x_i\geq \mathrm{ess}\sup \xi_i$ for some $i\in[d]$, the $i$-th coordinate of the gradient estimator $v(x)$ is  $0$ for any realization of $\xi$.  As a result,  projected SGD may not perform any update on the $i$-th coordinate and get stuck. 
RSG addresses this issue by adding a regularization. For $x^t$ such that $x_i^t\geq \mathrm{ess}\sup \xi_i$ for some $i\in[d]$,  RSG would perform the update on the $i$-th coordinate such that
\begin{equation}
\label{eq:RSG_phase_one}
x_i^{t+1} = (1-\gamma_t \lambda) x_i^t.
\end{equation}
\textcolor{black}{Denote $ \mathcal{X}^*_\mathrm{local}:=\{x\mid x_i\geq  \mathrm{ess}\sup \xi_i, \text{ for some } i\in[d], x\in\mathcal{X}\}$.  
\eqref{eq:RSG_phase_one} implies that RSG first shrinks the decision variable component-wisely to find a $x^t$ such that $x^t_i\leq \mathrm{ess}\sup \xi_i$ for any $i\in[d]$, and hence avoid convergence to any points in $\mathcal{\mathcal{X}^*_\mathrm{local}}$. According to \eqref{intuition:RSG}, $\mathcal{X}^*_\mathrm{local}$ is exactly the set of local solutions  that we intend to avoid. Hence, regularization ensures convergence to an approximate stationary point of $F^\lambda$ such that $\|[\nabla g(x)]^{-1}\|$ is finite.  For a small $\lambda>0$, such a $x$ is also an approximate stationary point of $F$, and thus an approximate global optimal solution of $F$ by \eqref{intuition:RSG}. In summary, RSG uses regularization to avoid vanishing gradient and ensure global convergence, while regularization in statistical learning literature is usually for avoiding overfit~\citep{vapnik1999nature}. 
}
\end{example}

\subsubsection{Intuition of MSG (Algorithm \ref{alg:MSG})}
As for MSG, the key step is to design a gradient estimator of  $F$ such that each update on the original space $x$ mirrors the gradient update of the convex objective $G$ on the reformulated space $u$ {as illustrated by Figure \ref{fig:demo}}. Denote the stepsize as $\gamma$.

Next, we illustrate how to build the gradient estimator $v_F$. For ease of demonstration, consider the simplified setting when $\mathcal{X}=\mathcal{U}=\RR^d$. The gradient descent update on $G$ for a point $u= g(x)$ is $u^{\prime} = u - \gamma \nabla G(u).$ Denote $x^\prime := g^{-1}(u^\prime)$.
If MSG mirrors the exact gradient descent update on $G$ using a gradient estimator $\tilde v_F(x)$, one should have
$
x^\prime = x -\gamma \tilde v_F(x).
$
Therefore, it holds that
\begin{equation}
\label{eq:MSG_desired_update}
\begin{aligned}
    \tilde v_F(x) 
= & 
    -\frac{x^\prime-x}{\gamma} = -\frac{g^{-1}(u^\prime)-g^{-1}(u)}{\gamma}\\
\approx & 
    -[\nabla g^{-1}(u)]^{\top} (u^\prime-u) = [\nabla g(x)]^{-\top} \nabla G(u) = 
    [\nabla g(x)]^{-\top}[\nabla g(x)]^{-\top} \nabla F(x),
\end{aligned}
\end{equation}
where we assume $\nabla g^{-1}(u)$ exists and use the first-order approximation in the second line. Notice that $u^\prime - u=-\gamma \nabla G(u)$. As long as the stepsize $\gamma$ is small, the approximation error is controlled by $\cO(\gamma^2)$. It motivates us to design a stochastic estimator $v_F(x)$ for $[\nabla g(x)]^{-\top}[\nabla g(x)]^{-\top} \nabla F(x)$. One can also interpret $[\nabla g(x)]^{-\top}[\nabla g(x)]^{-\top}$ as a pre-conditioning matrix.

It remains to build efficient estimators of $[\nabla g(x)]^{-1}$ with small bias and small variance at a low sampling cost. To achieve that, we utilize the well-known equality for infinite series of matrices. Let $\Lambda\in\RR^{d\times d}$ be a symmetric random matrix, $\EE \Lambda$ be its expectation and $I$ denote the identity matrix. Suppose that $\EE \Lambda$ is invertible and  $0\prec \EE \Lambda\prec I $. It holds that 
$$
[\EE \Lambda]^{-1} = \sum_{i=0}^\infty (I-\EE \Lambda)^i = \sum_{k=0}^\infty \prod_{i=1}^k (I-\EE \Lambda^i) = \sum_{k=0}^\infty \EE\prod_{i=1}^k  (I-\Lambda^i){\approx\sum_{k=0}^K \EE\prod_{i=1}^k  (I-\Lambda^i),}
$$
where $\prod_{i=1}^k (I-\EE \Lambda^i)=I$ if $k=0$ and $\{\Lambda^i\}_{i=1}^k$ are i.i.d. samples.
Utilizing a randomization scheme over $k\in \NN$, one can construct an estimator of $[\EE \Lambda]^{-1}$. Particularly, for an integer $K>0$, to estimate $[\nabla g(x)]^{-1}$,  we construct the following estimator:
generate $k$ uniformly from $\{0,\ldots,K-1\}$, generate i.i.d. sample $\{\xi^{i}\}_{i=1}^k$ from $\PP(\xi)$, and form the following estimator 
\begin{equation}
\label{eq:matrix_inverse_estimator}
[\nabla \hat g(x)]^{-1} = 
\begin{cases}
\frac{K}{cL_\phi} \Pi_{i=1}^{k} \Big(I - \frac{\nabla \phi(x,\xi^{i})}{cL_\phi}\Big) & \text{ for $k\geq 1$,}\\
 \frac{K}{cL_\phi} I & \text{ for $k=0$},
\end{cases} 
\end{equation}
where $c>1$ is to ensure that ${\nabla \phi(x,\xi^{i})}/{cL_\phi}\prec I$.
Although $\nabla g$ is a diagonal matrix in our problem, such estimators are used for more general matrix inverse estimation, e.g., estimating inverse Hessian matrix in bilevel optimization~\citep{hong2020two,hu2024contextual}.
\begin{lemma}
\label{lm:properties_of_inverse_expectation_estimation}
Under Assumption \ref{assumption:general}(c)(d) and Assumption \ref{assumption:general_2}(b),  the bias and the second moment of  the estimator \eqref{eq:matrix_inverse_estimator} with a constant $c>1$ satisfy:
$$
\|\EE [\nabla \hat g(x)]^{-1}-[\nabla g(x)]^{-1}\|\leq \frac{1}{\mu_g}\Big(1- \frac{\mu_g}{cL_\phi}\Big)^K,\quad
\EE \|[\nabla \hat g(x)]^{-1}\|^2 \leq \frac{K^2}{c^2L_\phi^2}.
$$
Moreover, the number of samples to construct the estimator in expectation is $(K-1)/2$.
\end{lemma}
Note that one could also use other distributions rather than the uniform distribution over $\{0,\ldots, K-1\}$. 
We defer related discussions and the proof to Appendix \ref{proof:inverse_estimation}.

\begin{algorithm}[t]%\footnotesize
	\caption{Mirror Stochastic Gradient (MSG)}
	\label{alg:MSG}
	\begin{algorithmic}[1]
		\REQUIRE  Number of iterations $ T $, stepsizes $\{\gamma_t\}_{t=1}^{T}$, initialization $x^1$, regularization parameter $\lambda$.
		\FOR{$t=1$ to $T$ \do} 
		\STATE Draw two independent samples $k_1$, $k_2$ uniformly from $\{0,\ldots,K-1\}$, draw i.i.d. samples $\{\tilde \xi^{ti}\}_{i=1}^{k_1}$, $\{\tilde \xi^{tj}\}_{j=1}^{k_2}$ to construct two estimators of $[\nabla g(x^t)]^{-1}$:
		$$
		[\nabla \hat g^{A}(x^t)]^{-1} = \frac{K}{2L_\phi} \Pi_{i=1}^{k_1} \Big(I - \frac{\nabla \phi(x^t,\tilde\xi^{ti})}{2L_\phi}\Big); \ [\nabla \hat g^{B}(x^t)]^{-1} = \frac{K}{2L_\phi} \Pi_{j=1}^{k_2} \Big(I - \frac{\nabla \phi(x^t, \tilde\xi^{tj})}{2L_\phi}\Big).
		$$
		(By convention, let $\Pi_{i=1}^{0} \Big(I - \frac{\nabla \phi(x^t,\tilde\xi^{ti})}{2L_\phi}\Big)=I$.)
		\STATE Draw a sample $\xi^t$ from $\mathbb{P}(\xi)$ and construct a gradient estimator
		\begin{center}
		$
		v_F(x^t)=[\nabla \hat g^{A}(x^t)]^{-\top}[\nabla \hat g^{B}(x^t)]^{-\top} \nabla \phi(x^t,\xi^{t})^\top \nabla f(\phi(x^t,\xi^{t})) + \lambda x^t.
		$
		\end{center}
		\STATE Update $ x^{t+1}=\Pi_{{\mathcal{X}}}(x^t-\gamma_t v_F(x^t))$.
		\ENDFOR
		\ENSURE $\hat x^T$ is selected uniformly from $ \{x^t\}_{t=1}^T$.
	\end{algorithmic}
\end{algorithm}

Based on the above discussion, we formally describe MSG in Algorithm \ref{alg:MSG}. Line 2 and Line 3 in MSG are to build a stochastic gradient estimator of $[\nabla g(x)]^{-\top}[\nabla g(x)]^{-\top}\nabla F(x)$, 
where $[\nabla g(x)]^{-\top}[\nabla g(x)]^{-\top}$ acts as a preconditioning matrix that rescale the gradient of the nonconvex objective $F$. 
For $\phi(x,\xi) = x\wedge\xi$, the $i$-th diagonal entry of $[\nabla g(x)]^{-1}$ is $(1-H_i(x_i))^{-1}$, where $H_i$ is the cumulative distribution function of $\xi_i$. Thus the preconditioning parameter enlarges all coordinates of $\nabla F(x)$. To analyze the convergence of MSG, we need to characterize the second moment of $v_F(x)$. To avoid potential dependence issues, we use two independent sets of samples of $\xi$ to estimate the first and the second $[\nabla g(x)]^{-1}$ terms in Line 2 of MSG. Also note that in MSG, we use the matrix inverse estimator \eqref{eq:matrix_inverse_estimator} with $c=2$ for simplicity. In addition, we use an independent sample $\xi^t$ to build a gradient estimator of $F(x)$. The regularization term $\lambda x^t$ in line 3 of MSG is also used to avoid vanishing gradient issues in practice, as we did for RSG.

\section{Global Convergence and Complexities Bounds}
\label{sec:global_convergence}
In this section, we demonstrate the global convergence and the sample and gradient complexities of RSG and MSG to achieve an $\eps$-optimal solution. {For ease of reference, we summarize the assumptions needed for SAA+SG, RSG, and MSG and their sample complexities bounds in Table \ref{tab:algos_assumption} in Appendix \ref{appendix:discussion_on_assupmtions}. The table also includes the global convergence of MSG for NRM applications when $\xi$ satisfies a discrete distribution.} Note that the intuition of RSG and MSG discussed in Section \ref{sec:alg} builds upon the unconstrained setting. When extended to a constrained setting, the following lemma is the key property that we use to address the hardness brought in by projection.  We defer the proof to Appendix \ref{sec:switch_projection_transformation}.

\begin{lemma}
\label{lm:switch_transformation_projection}
Suppose that $\mathcal{X}$ is a box constraint and Assumption \ref{assumption:reformulation}(a)(b) holds. For any $x\in\mathcal{X}$, we have $g\Big(\Pi_\mathcal{X}(x)\Big) = \Pi_\mathcal{U}\Big(g(x)\Big).$
\end{lemma}
Lemma \ref{lm:switch_transformation_projection} says that one could exchange the projection operator and the transformation operator when both $\mathcal{X}$ and $\mathcal{U}$ are box constraints and $g$ is a component-wise non-decreasing function. In the proof of RSG, Lemma \ref{lm:switch_transformation_projection} plays a key role in establishing an upper bound of the gradient mapping of $G$ using the gradient mapping of $F$ (see \eqref{eq:gradient_mapping_FG}). In the proof of MSG, Lemma \ref{lm:switch_transformation_projection} enables us to conduct the analysis in a way similar to the unconstrained case. 
The following theorem establishes the global convergence of RSG. The proof is in Appendix \ref{sec:proof_of_RSG}. 
\begin{theorem}
\label{thm:RSG}
Suppose that  Assumptions \ref{assumption:reformulation}, \ref{assumption:general} and \ref{assumption:general_2}(a)(b) hold and that $\nabla F$ is $S_F$-Lipschitz continuous. \textcolor{black}{For RSG with stepsizes $\gamma_t = \gamma =T^{-1/2}$ and a regularization parameter $\lambda \geq 0$}, there exists  a constant $M>0$ such that the expected error of RSG is upper bounded by
$$
\textcolor{black}{
    \EE[F(\hat x^T) - F(x^*)] 
\leq 
    (2L_\phi D_\mathcal{X}+L_fL_\phi/(2{S_F}\mu_g)) \max\{\mu_g^{-1}, L_\phi\}\sqrt{2 MT^{-1/2} + 2\lambda^2 D_\mathcal{X}^2}.
}
$$
\end{theorem}
\textcolor{black}{In the above theorem, the term $M T^{-1/2}$ comes from a stationary convergence of RSG, where the constant $M=\cO((S_F+\lambda) (L_f^2 L_\phi^2+\lambda ^2D_\mathcal{X}^2))$ is explicitly given in the analysis, the term $\lambda^2D_\mathcal{X}^2$ comes from adding the regularization, and the remaining terms come from building a relationship between stationary convergence and global convergence utilizing convexity of the reformulated problem. 
Theorem \ref{thm:RSG} implies that setting $\lambda \in [0, D_\mathcal{X}^{-1}T^{-1/4}]$, and $T=\cO(D_\mathcal{X}^4\eps^{-4})$ for any $\eps\in(0,1)$, we have $\EE[F(\hat x^T) - F(x^*)] = \cO(\eps)$.
Since RSG uses one sample and computes one gradient of $f$ at each iteration, for $\hat x_T$ to be an $\eps$-optimal solution of $F$, the sample and gradient complexities of RSG are both $\cO(\eps^{-4})$ in terms of the dependence on the accuracy $\eps$.  } We point out that the complexity of RSG has a $D_\mathcal{X}^4$ dependence on the radius, which, in the worst case, is equivalent to a quadratic dependence on the dimension $d$ since $\mathcal{X}$ is a box constraint. We will show later that such dependence does not have a significant impact on numerical experiments. 

In terms of analysis, we first build up the stationary convergence of projected SGD on constrained smooth optimization measured by the norm of the gradient mapping~\citep{davis2018stochastic, drusvyatskiy2019efficiency} and then establish a relationship between the stationary convergence and the global convergence. 

The following theorem demonstrates the global convergence of MSG. We defer the proof to Appendix \ref{sec:proof_of_msg}. Unlike RSG, the analysis of MSG does not require $\nabla F$ to be Lipschitz continuous as it directly demonstrates global convergence.
\begin{theorem}
\label{thm:msg}
Suppose that  Assumptions \ref{assumption:reformulation}, \ref{assumption:general} and \ref{assumption:general_2} hold. For MSG with stepsizes $\gamma_t = \gamma $, and regularization parameter $\lambda\geq0$, the expected error of MSG is upper bounded by
\begin{equation}
\label{eq:expected_error_MSG}
 \begin{aligned}
    \EE[F(\hat x^T)-F(x^*)]
\leq &
    \frac{\|u^1 -u^*\|^2}{2\gamma T}+  \gamma(L_\phi^2  +2L_\phi D_\mathcal{X} S_g )\Big(\frac{K^4 L_f^2}{16L_\phi^2}+ \lambda^2 D_\mathcal{X}^2\Big) + 2 L_\phi^2 D_\mathcal{X}^2 \lambda \\
    & +
    L_\phi^2 D_\mathcal{X} \frac{K L_f+2L_f}{\mu_g}\Big(1- \frac{\mu_g}{2L_\phi}\Big)^K.
\end{aligned}
\end{equation}
\end{theorem}
In the right-hand-side of \eqref{eq:expected_error_MSG}, the first term coming from a  telescoping sum appears in SGD analysis~\citep{nemirovski2009robust}, the second terms comes from the variance of the  estimator $v_F$ and the approximation error of MSG to the virtual SGD update on $G$, the third term comes tje from regularization, and the fourth term comes from the bias of estimating matrix inverse $[\nabla g(x)]^{-1}$. 

Setting $\gamma =  (D_\mathcal{X}T)^{-1/2}$, $\lambda\in[0, (D_\mathcal{X}T)^{-1/2}]$, $K = \cO(\log(D_\mathcal{X}\eps^{-1}\log(\eps^{-1})))$, and $T = \tilde \cO(D_\mathcal{X}^2\eps^{-2})$, we have $\EE[F(\hat x^T) - F(x^*)] =\cO(\eps)$.
Since MSG uses at most $2K-1$ number of samples per-iteration, the sample and gradient complexities of MSG are both $\tilde \cO(\eps^{-2})$. In terms of the dependence on the accuracy $\eps$, the theorem implies that the nonconvex problem \eqref{problem:math_original_extension} under Assumption \ref{assumption:general} and Assumption \ref{assumption:general_2} is fundamentally no harder than the classical stochastic convex optimization. \textcolor{black}{Note that the iteration complexity $T$ also depends on the radius of the domain $\mathcal{X}$, i.e., $T\propto D_\mathcal{X}^2$. Since $\mathcal{X}$ is a box constraint, in the worst case, $T$ scales linearly in dimension $d$, which is still better than that of RSG. We are unaware of any method that could get rid of the dimension dependence, and we leave it for future investigation.}

Next, we discuss the efficiency of MSG via showing a lower bound for problem \eqref{problem:math_original_extension}. For this purpose, note that \citet{agarwal2009information} developed an $\cO(\eps^{-2})$ lower bounds on the gradient complexity of any black-box stochastic first-order algorithms for obtaining an $\eps$-optimal solution of $\min_{x\in\mathcal{X}}F(x)$, where $F$ is convex and Lipschitz continuous. Interestingly, the hard instance that they used to construct the lower bound happens to be a special case of problem \eqref{problem:math_original_extension} when $\phi(x,\xi)=x+\xi$, $F(x) =\EE [f(x+\xi)]$ and $\mathcal{X}=[-10,10]\subset\RR$. It is easy to verify that Assumptions \ref{assumption:reformulation} and \ref{assumption:general}(c-f) hold for $\phi(x,\xi)= x+\xi$. Though the hard instance in \citet{agarwal2009information} is constructed by $\EE [f(x+\xi)]$, \citet{agarwal2009information} considered lower bounds for black-box stochastic first-order algorithms, i.e., algorithms that uses a gradient estimator of $F$ that can be of any form as long as it is unbiased and has bounded variance but does not have to be $\nabla f(x,\xi)$. This is different from MSG, which additionally uses $\nabla \phi$ to build up a preconditioning matrix and thus is not a black-box algorithm.  Though the results of \citet{agarwal2009information} is not directly applicable to MSG, we could use their analysis to establish a $\cO(\eps^{-2})$ lower bound on the gradient
complexity of any black-box stochastic first-order algorithms for solving problem \eqref{problem:math_original_extension}.  In addition, such lower bounds imply that the sample and gradient complexities of MSG match the best possible black-box stochastic gradient methods for solving \eqref{problem:math_original_extension} in terms of accuracy $\eps$ if ignoring the logarithmic term.

\section{Applications}
\label{sec:application}
In this section, {we first discuss the board applicability of the studied stochastic nonconvex optimization in operations management. Then, we model the air-cargo NRM problem with random demand, two-dimensional capacity, consumption, and routing flexibility under booking limit control as a special case of Problem \eqref{problem:math_original_extension}. We further show the global convergence of MSG on the NRM problem. } Interested readers please refer to Appendix \ref{sec:appendix_NRM_alg} for the modeling of passenger NRM and refer to \citep{feng2015air,klein2020review} for a comprehensive review of air-cargo NRM. Our booking limit control adapts to two interesting extensions of managing uncertain capacity in \cite{previgliano2021managing}, see Appendix \ref{appendix:computation_details}.

\subsection{{Operations Management Applications}}
\label{sec:examples}
{Several operations management applications are special cases of the studied problem \eqref{problem:math_original_extension}. Dynamic multisourcing problems with random capacities~\citep{chen2018preservation, feng2019dynamic}, assemble-to-rrder system with a random capacity ~\citep{chen2018preservation, feng2019dynamic} and lost-sale inventory problems with random supply~\citep{chen2022learning} are special cases of problem \eqref{problem:math_original_extension} with $\phi(x,\xi)=x\wedge\xi$. Newsvendor with procurement from multiple suppliers~\citep{dada2007newsvendor} is a special cases of problem \eqref{problem:math_original_extension} with $\phi_i(x_i,\xi_i)=x_i\xi_i/(x_i+\alpha\xi_i^\kappa)$ and 
  random supply from one producer to multiple firms~\citep{tang2014pay} is a special case of problem \eqref{problem:math_original_extension} with $\phi_i(x_i,\xi_i)=(x_i/(x_i+\xi_i))k$. Interested readers may refer to other applications in~\citet{feng2018supply,feng2022applications}. We also give an example of assemble-to-order systems in Appendix \ref{appendix:application_ATO}.

}

\subsection{Modeling for Air-cargo Network Revenue Management}\label{sec:air_cargo_model}
From a temporal perspective, the air-cargo NRM problem consists of \textit{reservation stage} and \textit{service stage}. During the reservation stage, we have to decide whether to accept reservation requests. Then at the service stage, we aim to minimize the penalty of rejecting show-ups by accommodating show-ups within the limited random capacity and potential routing options. Thus, there
are four significant factors in air-cargo NRM problems, two-dimensional capacity (weight and volume), random capacity, random consumption, and routing flexibility (i.e., demand class with specified origin-destination pair can be shipped via any feasible route in the airline network). In this paper, we consider these four factors all at once. {\color{black} \citet{barz2016air} considered the same setting and proposed a DPD method. However, their DPD method only heuristically addressed the routing decision while we explicitly model the routing decisions as decision variables.} In what follows, we formulate the problem using booking limit control as a two-stage stochastic optimization problem.

At the start of the reservation stage, we decide the booking limits, denoted by decision vector $x=(x_1,x_2,\dots,x_d)^\top$, for $d$ demand classes. Under booking limit control, we accept new requests for a demand class $i$ unless the booking limit $x_i$ is reached. We assume that the aggregated demand during the whole reservation stage, denoted by vector $\tilde{D}=(\tilde{D}_1,\tilde{D}_2,\dots,\tilde{D}_d)^\top$, is random and component-wise independent. Note that each request comes with a random weight and a random volume that are independent of $x$ and $\tilde{D}$, and will reveal at the service stage. In total, we accept up to $x\wedge \tilde D$ reservations during the reservation stage. At the end of the reservation stage, cancellations and no-shows are realized.

We use the random vector $\tilde{Z}=(\tilde{Z}_1,\tilde{Z}_2,\dots,\tilde{Z}_d)^\top$ to represent the number of show-ups for the service. Thus, the number of show-ups can be written as a function of the booking limit and random demand, i.e., $\tilde{Z}=\tilde{Z}(x\wedge\tilde{D})=(\tilde{Z}_1(x_1\wedge\tilde{D}_1),\tilde{Z}_2(x_2\wedge\tilde{D}_2),\dots,\tilde{Z}_d(x_d\wedge\tilde{D}_d))^\top$. We assume that $\tilde{Z}_i(x_i),~i\in[d]$ follows a Poisson distribution with a coefficient $p_ix_i$ and that $\tilde{Z}(x)$ is component-wise independent. Without loss of generality, we assume that the no-shows or cancellations are not refundable. Note that all-show-up setting is a special case with $\tilde{Z}_i(x_i)=x_i$ and $p_i=1,~i\in[d]$. We focus on Poisson show-ups rather than the more practical binomial show-ups because Poisson is well suited to the continuous optimization framework, and the same justification can be found in \citet{karaesmen2004overbooking}. Moreover, we consider the continuous booking limit $x\in\mathbb{R}^d$, which allows fractional acceptance throughout the paper for the same reason. \textcolor{black}{We also want to highlight that our algorithms can be heuristically adapted to the discrete booking limit and binomial show-up setting with more detailed discussions in Appendix \ref{sec:appendix_NRM_alg}.}

At the service stage, there are $m$ inventory classes associated with a two-dimensional random capacity, where the first dimension is weight capacity  $\tilde{c}_w=(\tilde{c}_{w1},\dots,\tilde{c}_{wm})^\top$ and the second dimension is volume capacity $\tilde{c}_v=(\tilde{c}_{v1},\dots,\tilde{c}_{vm})^\top$.  Each accepted demand from class $i$ has random weight $\tilde{W}_i$ and volume $\tilde{V}_i$. We assume that both random weight $\tilde{W}=(\tilde{W}_1,\dots,\tilde{W}_m)^\top$ and volume $\tilde{V}=(\tilde{V}_1,\dots,\tilde{V}_m)^\top$ are realized at the beginning of the service stage and are independent of the booking limit $x$ and the aggregated demand $\tilde{D}$.  The revenue gained by accepting one unit reservation is a function of random weight and volume, denoted by $r(\tilde{W},\tilde{V})=(r_1(\tilde{W},\tilde{V}),\dots,r_d(\tilde{W},\tilde{V}))^\top$. In the air-cargo industry, a common practice is to charge $r_i(\tilde{W},\tilde{V})=\theta_1\max\{\tilde{W}_i,\tilde{V}_i/\theta_2\}$ for demand class $i\in[d]$, with some constants $\theta_1,\theta_2$ \citep{barz2016air}. With the similar structure, we define $l(\tilde{W},\tilde{V})=(l_1(\tilde{W},\tilde{V}),\dots,l_d(\tilde{W},\tilde{V}))^\top$ as the penalty of rejecting one unit reservation. We assume each demand class $i\in[d]$ can be satisfied by $K_i$ different routes and define the binary parameter $b_{ijk}\in\{0,1\}$ to represent whether the inventory class $j$ is required to satisfy the demand from the $k^{th}$ route of demand class $i$. During the service stage, the first decision is the amount of served show-ups $w=(w_1,\dots,w_d)^\top$ under limited capacities, and the second decision is the routing decision, where we use variable $y_{ik}$ to denote the amount of demand allocated to $k^{th}$ route of demand class $i\in[d]$. Let $\Gamma(z, W, V, c_w,c_v)$ denote the penalty of rejecting accepted demand during the service stage. Then the air-cargo NRM problem under booking limit control has the following mathematical formulation:
\begin{equation}\label{eqn:ari_cargo_two_dim}
	\begin{aligned}
		\max_{x\geq 0} ~~& \mathbb{E}_{\tilde{D}}[f(x\wedge\tilde{D})],
	\end{aligned}
\end{equation}
where $f(x)=\mathbb{E}_{\tilde{W}, \tilde{V}}[r(\tilde{W}, \tilde{V})^\top x]-\mathbb{E}_{\tilde{Z}(x),\tilde{W}, \tilde{V}, \tilde{c}_w,\tilde{c}_v}[\Gamma(\tilde{Z}(x),\tilde{W}, \tilde{V}, \tilde{c}_w,\tilde{c}_v)]$ and
\begin{equation*}
	\begin{aligned}
 		\Gamma(z, W, V, c_w,c_v)=&~\min_{y,w}~ l(W, V)^\top (z-w)\\
		\textit{s.t.}&~\sum_{i=1}^n\sum_{k=1}^{K_i} b_{ijk} y_{ik} W_i\leq {c_w}_j, \forall j\in[m];\quad\sum_{i=1}^n\sum_{k=1}^{K_i} b_{ijk} y_{ik} V_i\leq {c_v}_j, \forall j\in[m];\\
		 &~w_i=\sum_{k=1}^{K_i} y_{ik}, \forall i\in[d];\quad  0\leq w\leq z;\quad y\geq 0.\\
	\end{aligned}
\end{equation*}
In the above model, the first and the second constraints represent the weight and volume capacities constraints for all inventory classes $j\in[m]$ . The third constraint $w_i=\sum_{k=1}^{K_i} y_{ik}$ indicates that the total accepted demand $w_i$ of class $i$ is allocated over $K_i$ different routes. From a modeling perspective, to the best of our knowledge, we are the first to explicitly model the optimal routing decisions under the booking limit control. In comparison, the DPD method proposed in \citet{barz2016air} only heuristically splits the reservations of class $i$ with the same origin-destination equally into fixed $K_i$ sub-classes. 

{
\subsection{Theoretical Results for NRM Applications}\label{sec:convergence_NRM}
} 
{
In this subsection, we discuss the global convergence of the proposed algorithms in NRM applications and discuss what happens if there lacks Assumption \ref{assumption:general_2}. } The next lemma specifies the conditions needed for $\phi(x,\xi)=x\wedge \xi$ to ensure the assumptions needed for global convergence. We defer the proof to Appendix \ref{proof_of_verify_x_wedge_xi}. { We also specify conditions needed for the other three $\phi$ functions to ensure Assumption \ref{assumption:general_2} in Appendix \ref{appendix:verify_all_phi}. A summary of the conditions is in Table \ref{tab:phi_assumption}.}

\begin{lemma}
\label{lemma:verify_x_wedge_xi}
For $\phi(x,\xi) = x\wedge\xi$ with component-wise independent random vector $\xi$, if the CDF of $\xi_i$, $H_i(x_i)$, is $S_g$-Lipschitz continuous and $1-H(\bar X_i) = \PP(\xi_i\geq \bar X_i)\geq \mu_g$ for any $i\in [d]$, then all needed assumptions on $\phi$ and $\PP(\xi)$ to ensure global convergence of RSG and MSG hold.  
\end{lemma}
Note that the convexity of the objective function and the gradient calculation follow a similar derivation in \cite{karaesmen2004overbooking} and are reproduced in Appendix  \ref{sec:appendix_NRM_alg} for completeness. Specifically, in the all-show-up case (there does not exist cancellations and no-shows) in NRM problems, the function $f$ in \eqref{eqn:ari_cargo_two_dim} is concave as NRM is a maximization problem. On the other hand, obtaining the gradient of $f$ requires solving a linear program (LP). Under the condition specified in in Lemma \ref{lemma:verify_x_wedge_xi},  Theorems \ref{thm:RSG} and \ref{thm:msg} imply that for solving air-cargo NRM \eqref{eqn:ari_cargo_two_dim} under the all-show-up case, MSG requires solving $\tilde\cO(\eps^{-2})$ number of LPs while RSG needs to solve $\cO(\eps^{-4})$ number of LPs to achieve an $\eps$-optimal solution.

In what follows, we discuss the convergence of MSG when Assumption \ref{assumption:general_2}(c), the smoothness of the transformation function $g$,  is missing. To ensure Assumption \ref{assumption:general_2}(c) holds, it requires Lipschitz continuous CDF assumption on $\xi$, meaning that $\xi$ is a continuous random vector. However, in NRM applications, the distribution of $\xi$ can be discrete, like Possion or multinomial. The next theorem shows the approximate global convergence rate of MSG without assuming $\xi$ is continuous. We defer the analysis to Appendix \ref{appendix:a2.3_fails}

\begin{theorem}
\label{cor:msg_without_smoothness}
Suppose that  Assumptions \ref{assumption:reformulation}, \ref{assumption:general} and \ref{assumption:general_2}(b) hold. For MSG with stepsizes $\gamma_t = \gamma =\cO(1/\sqrt{T})$, and regularization parameter $\lambda=\cO(1/\sqrt{T})$, for a discrete distribution $\xi$ with a support over $\mathbb{Z}^d$, the expected error of MSG is upper bounded by
\begin{equation}
 \begin{aligned}
    \EE[F(\hat x^T)-F(x^*)]
= &
    \tilde\cO\Big(1/\sqrt{T} + \sqrt{d \max_{k\in\mathbb{Z},i\in[d]} \PP(\xi_i=k)}\Big).
\end{aligned}
\end{equation}
When $\xi$ is a Poisson random vector with  an arrival rate vector $\beta$ or a multinomial distributed random vector with $\beta$ number of trails, the expected error bound becomes 
\begin{equation}
\label{eq:msg_approxition}
 \begin{aligned}
    \EE[F(\hat x^T)-F(x^*)]
\approx &
    \tilde\cO(1/\sqrt{T} + \sqrt{d/\beta}).
\end{aligned}
\end{equation}
\end{theorem}
When $\beta$ is large, the approximation in \eqref{eq:msg_approxition} comes from Stirling's formula and the property of Poisson and multinomial distributions. The theorem implies that the global convergence of MSG still holds even without the smoothness of $g$. The reason is that both Poisson and multinomial distributions can be well approximated by continuous normal distribution when $\beta$ is large. {The result can further extend to more general nonsmooth function $g$ when one can bound $\big\|\EE[ g(x^t - \gamma v_F(x^t)) - g(x^t) +\gamma \nabla g(x^t)^\top v_F(x^t)\mid u^t]\big\|$. Theorem \ref{cor:msg_without_smoothness} exploits the concavity of $g(x) =\EE[x\wedge \xi]$ and the discrete nature of $\xi$ that its support is on $\mathbb{Z}^d$ to bound this term. For other nonsmooth $g$ functions, one may consider leveraging  smooth approximations to provide an upper bound on this term, which may require a case-by-case study.}

In Appendix \ref{appendix:a2.3_fails}, we further show that the performance of both RSG and MSG are also not influenced even if Assumption \ref{assumption:general_2}(b) is lacking. The reason follows a similar discussion in Example \ref{example:RSG_on_truncation} that adding regularization can address the issue.
Our numerical experiments also support such an observation that regularization is crucial to escape local solutions.

\section{Computational Experiments and Results}\label{sec:numerics}

For network revenue management problems, we conduct extensive numerical experiments implemented in Python using the Gurobi linear programming solver.
The implementation details and revenue evaluation of the proposed RSG (Algorithm \ref{alg:RSG}), MSG (Algorithm \ref{alg:MSG}), and SAA+SG (Algorithm \ref{alg:sgd_on_g} in Appendix \ref{sec:sgd_u}) and other benchmark strategies including DLP, DPDs \citep{erdelyi2010dynamic, barz2016air}, and VCBP \citep{previgliano2021managing}, can be found in Appendix \ref{appendix:computation_details}. {Discussions on assumptions in numerical studies are in Appendix \ref{appendix:numerical_assumption_fails}. Numerical convergence comparison of RSG, MSG, and SAA+SG on a specific instance is in Appendix \ref{appendix:convergence_numerical_NRM}.}

\subsection{Passenger Network Revenue Management with the Random Capacity}\label{sec:NRM_passenger}
\textbf{Experimental Setup.} We use test examples from \cite{erdelyi2010dynamic}. Among these examples, the reservation stage is divided into $T=240$ discrete periods with specified arrival probability for each demand class at each period. We label these test instances by tuple $(N, \kappa, \delta, \sigma, p, \rho, \gamma)$ with definitions given as follows. (1) $N$: the network contains one hub and $N\in\{4,8\}$ spokes (see Figure \ref{fig:network_structure}(a)); (2) $\kappa$: airline offers a high and a low fare itinerary in each origin-destination pair, where the high fare is $\kappa\in\{4,8\}$ times the price associated with low fare class. Thus, the number of inventory classes (flight legs) is $2N$, and the number of demand classes (itineraries) is $2N(N+1)$; (3) $(\delta, \sigma)$: penalty of rejecting one unit show-up from demand class $i$ is $l_i=\delta r_i+\sigma \max\{r_{\hat{l}}:\hat{l}=1,\dots,d\}$ with $(\delta,\sigma)\in\{(4,0),(8,0),(1,1)\}$; 
(4) $p$: show-up probability is given by $p\in\{0.90,0.95\}$, which follows binomial distribution and is the same for all demand classes; (5) $\rho$: load factor $\rho\in\{1.2,1.6\}$ is defined as total expected demand divided by total capacity; (6) $\gamma$: the random capacity follows the truncated Gaussian distribution with range $[0,\infty)$ and two different levels of coefficient of variation $\gamma\in\{0.1,0.5\}$. In total, there are 96 different test instances.

\noindent \textbf{Comparison to Other Control Policies.}
In the following, we compare proposed methods to other control policies in the existing literature in two aspects: expected revenue and computation time. We consider {an alternative} setting under binomial random show-ups as mentioned in the experimental setup. Thus the function $f$ is only component-wise convex rather than convex as required in our theory. We still apply our algorithms and report the results. 

\begin{table}[ht]
\centering \scriptsize
\caption{%\scriptsize
Revenue Comparison of MSG and Other Benchmarks for Different Sets of Test Instances in Passenger NRM}
\renewcommand\arraystretch{1.2}
\label{tab:s-h-rand-cap-summary}
\begin{tabular}{@{}lrrrrrrrrrrrrr@{}}
\\[-1.8ex]\hline 
\hline\\[-1.8ex]
\multirow{2}{*}{Benchmark Strategies} & \multicolumn{2}{c}{$N$} & \multicolumn{2}{c}{$\kappa$} & \multicolumn{3}{c}{$(\delta,\sigma)$} & \multicolumn{2}{c}{$p$} & \multicolumn{2}{c}{$\rho$} & \multicolumn{2}{c}{$\gamma$} \\ \cmidrule(l){2-3} \cmidrule(l){4-5} \cmidrule(l){6-8} \cmidrule(l){9-10} \cmidrule(l){11-12} \cmidrule(l){13-14} 

                                      & 4         & 8         & 4           & 8           & (4,0)     & (8,0)      & (1,1)     & 0.9       & 0.95      & 1.2            & 1.6            & 0.1        & 0.5        \\ 
\hline\\[-1.8ex]                                      
MSG v.s. DLP                         & 23.9\%    & 63.2\%    & 43.8\%      & 43.3\%      & 9.1\%     & 57.3\%     & 64.2\%    & 42.3\%    & 44.7\%    & 65.0\%         & 22.1\%         & 10.4\%     & 76.7\%     \\
MSG v.s. DPD                        & 3.0\%     & 13.7\%     & 11.0\%       & 5.6\%       & 3.0\%     & 14.9\%      & 7.1\%     & 8.1\%     & 8.6\%     & 10.5\%          & 6.1\%          & 0.7\%      & 16.0\%      \\
MSG v.s. VCBP                          & 4.4\%     & 5.3\%     & 5.7\%       & 3.9\%       & 4.4\%     & 6.9\%      & 3.2\%     & 6.0\%     & 3.6\%     & 4.0\%          & 5.6\%          & 4.0\%      & 5.7\%      \\ 

\\[-1.8ex]\hline 
\hline\\[-1.8ex]
\end{tabular}
\vskip -0.2in
\end{table}

For the comparison in expected revenue, Table \ref{tab:s-h-rand-cap} in Appendix \ref{appendix:complete_numerical_passenger_NRM} documents the complete numerical results for all passenger NRM instances. In summary, there is no significant difference in the expected revenue between RSG, MSG, and SAA+SG. Thus, we only compare MSG to other control policies in the following. Averaging over all instances, MSG gains higher revenue than DLP, DPD, and VCBP by 43.6\%, 8.3\%, and 4.8\%, respectively. It is not surprising that DLP performs worst among all control policies since it does not account for the variance in demands, show-ups, and capacities. There are some interesting observations when we fix one factor and average over all instances with the fixed factor. For instance,  we evaluate the influence of the capacity variance factor by averaging over 48 instances with $\gamma=0.1$ and the other 48 instances with $\gamma=0.5$. Table \ref{tab:s-h-rand-cap-summary} summarizes such results. We find that DLP and DPD perform significantly worse in high capacity variance case $\gamma=0.5$ than the low variance case $\gamma=0.1$, while our booking limit control and VCBP can deal with the random capacity setting much better. \cite{previgliano2021managing} report a similar result that VCBP performs better than DPD in high capacity variance cases. We point out that the implemented DPD \citep{erdelyi2010dynamic} method is designed for random show-ups with deterministic capacity. Although we extend their DPD method to incorporate the random capacity using the sample average to approximate the boundary value function, we admit there might exist other DPD methods specifically designed for the random capacity. For completeness, we compare our booking limit control to DPD under exactly the same 48 deterministic capacity instances (Table 1 and Table 2 in \cite{erdelyi2010dynamic}) and report that DPD performs better than booking limit control by 1.22\%. However, there is no significant revenue gap after we resolve our booking limit model 10 times. The test examples in two columns $(\delta,\sigma)=(4,0)$ and $(8,0)$ in Table \ref{tab:s-h-rand-cap-summary} have the increasing penalty of rejecting customers. The performance gap increases from the low penalty to the high penalty setting, indicating our booking limit control makes a better trade-off between the high-fare and low-fare classes.

\begin{table}[ht]
\centering %\scriptsize
\normalsize
\caption{Computation Time Comparison (CPU seconds)} 
\label{tab:comp_time} 
\begin{tabular}{cccccc}
\\[-1.8ex]\hline 
\hline\\[-1.8ex]
{ Number of Spokes } & RSG & MSG & SAA+SG & VCBP & DPD \\ \hline\\[-1.8ex]
$N=4 $               & 12  & 8   & 16     & 45   & 57  \\
                         $N=8 $            & 44  & 32  & 57     & 94   & 85 \\
\\[-1.8ex]\hline 
\hline
\end{tabular}
\vskip -0.1in
\end{table}

The comparison in computation time is summarized in Table \ref{tab:comp_time}. Since the number of spokes $N$ is the key parameter that affects the computation time, we report average CPU seconds averaged over all $N=4$ or $N=8$ instances. The stopping criteria of VCBP, RSG, MSG, and SAA+SG are specified in Appendix \ref{appendix:computation_details}. With different spokes $N$, we get different test instances with $n=2N(N+1)$ demand classes and $m=2N$ inventory classes. 
Next, we discuss the per-iteration computational costs.
VCBP solves one LP with $n$ decision variables and $2n+m$ constraints at each iteration and uses the backward path to get gradients with the computation cost of $O(mT)$, where $T$ is the number of total arrivals. Our proposed algorithms solve the same size LP with $n$ decision variables and $2n+m$ constraints at each iteration, and the computation cost of remaining arithmetic operations is mild compared to the LP solving. DPD solves $m$ single-leg dynamic programming, and the computation cost is bounded by $O(mT^2)$. 
Our results in Table \ref{tab:comp_time} show that MSG has the lowest computation cost at both $N=4$ and $N=8$. However, the computation cost of the DPD method scales better with respect to $N$. It is worth mentioning that the scalability with respect to $N$ of VCBP is the same as our algorithms as they all solve one LP of the same size at each iteration. 
In addition, although we only focus on fixed $T=240$ in our computation experiments and do not compare the scalability in $T$, the computation cost of our proposed algorithms for booking limits is independent of $T$ since we aggregate the reservation periods into a single stage.

\subsection{Air-cargo Network Revenue Management}\label{NRM_aircargo}

\textbf{Without Routing Flexibility Experimental Setup.} Since we do not have access to all the test instances for air-cargo NRM in \cite{barz2016air}, we construct similar instances based on parameters listed in Appendix \ref{appendix:air_cargo_param}, including the demand class label, average weight, average volume, the origin, the destination, and the per-unit-revenue. Note that the per-unit-revenue is the parameter $\theta_1$ in revenue $r_i=\theta_1\max\{W_i, V_i/\theta_2\}$ introduced in Section \ref{sec:air_cargo_model}. 

We adopt a similar setup as \citet{barz2016air} and set parameters as follows: $\theta_2=0.6$; the penalty is 2.4 times of revenue, e.g., $l=2.4r$; the coefficient of correlation between weight and volume consumption and the coefficient of correlation between the weight and volume capacity are both $0.8$; the planning horizon is $T=240$, which is consistent with the previous passenger network revenue case; we neither consider the no-show nor cancellation, which can be easily incorporated, to be consistent with the air-cargo DPD \citep{barz2016air} (ACDPD); all demand classes arrive with equal probability over the reservation stage; we consider two different levels of the coefficient of variation in the random consumption $CV_D\in\{0.1,0.4\}$, two levels of the coefficient of variation in the random capacity $CV_C\in\{0.1,0.4\}$, and two scenarios of the average load factor levels (i.e., $\mathbb{E}$[demand]/$\mathbb{E}$[capacity] with the fixed expected demand and varying expected capacity). The network structure is spoke-hub given in Figure \ref{fig:network_structure} (a). Since this network only contains one feasible route for any given origin-destination pair, there is no routing flexibility.

\begin{figure}[ht]
\centering
\subfigure[Spoke-and-hub Network without Routing Flexibility]{
    \begin{minipage}[t]{0.46\linewidth}
        \centering
        \includegraphics[width=0.5\linewidth]{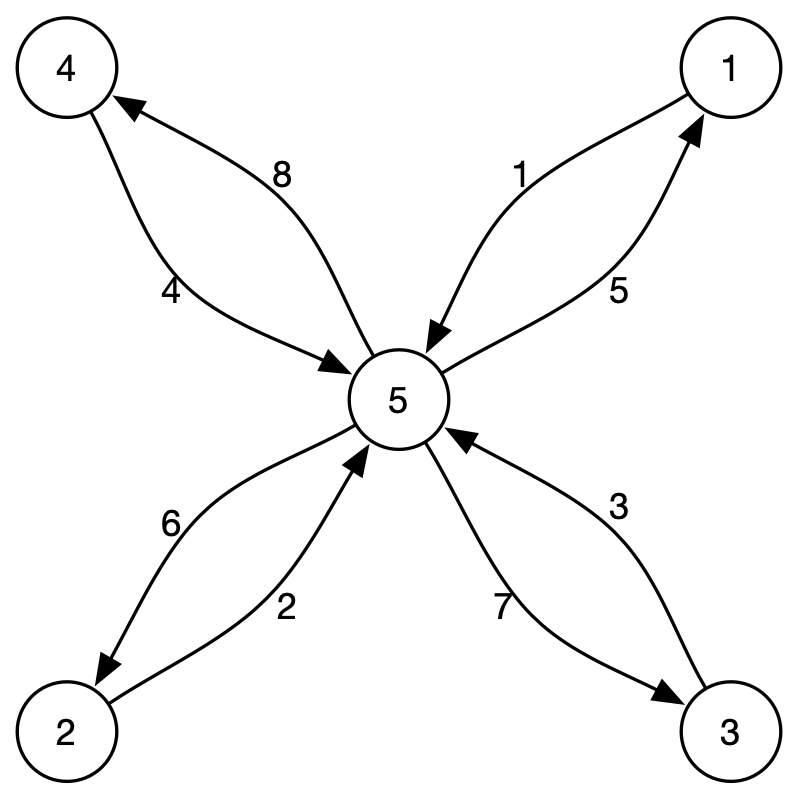}
    \end{minipage}
}%
\subfigure[Network with Routing Flexibility]{
    \begin{minipage}[t]{0.46\linewidth}
        \centering
        \includegraphics[width=0.5\linewidth]{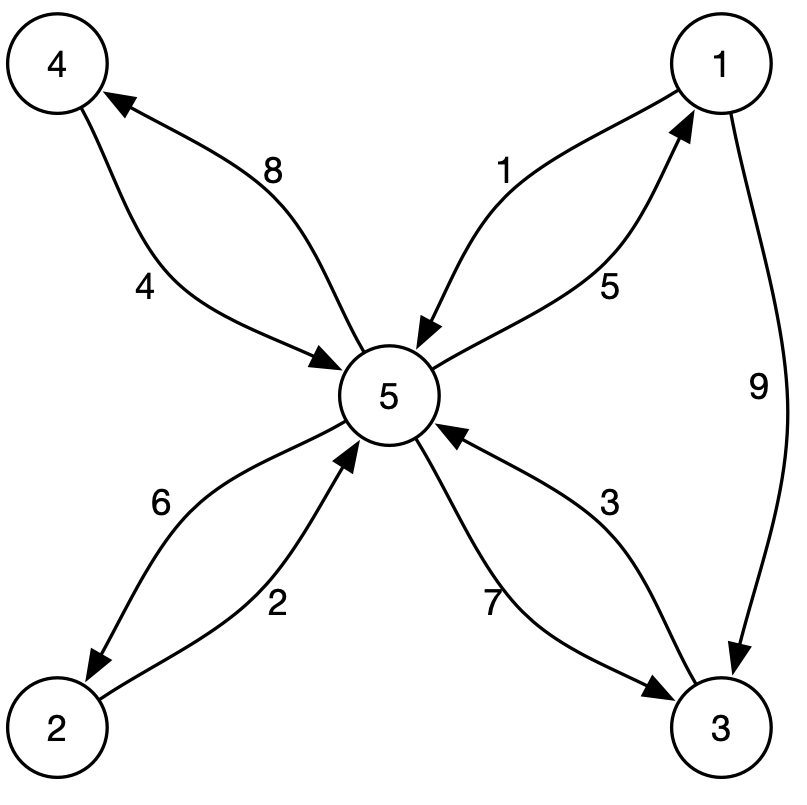}\\
    \end{minipage}%
}%
\caption{Flight Network Structure in Air-cargo NRM}
\label{fig:network_structure}
\vskip -0.1in
\end{figure}

As shown in Table \ref{tab:air_cargo_result} under ``Without Routing Flexibility'' columns, the booking limit control policy computed by MSG outperforms ACDPD by an average 12.86\% among all test instances. We observe a similar trend as in passenger network instances (see the $\gamma$ column in Table \ref{tab:s-h-rand-cap-summary}) that MSG outperforms ACDPD at a more significant level when the capacity and demand have higher variances, indicating that our MSG method accounts for randomness more effectively.

\begin{table}[!ht]
\centering \scriptsize 
%\normalsize
\caption{Revenue Comparison of MSG and ACDPD in Air-cargo NRM without/with Routing Flexibility}
\renewcommand\arraystretch{1.0}
\label{tab:air_cargo_result}
\begin{tabular}{ccc|rrr|rrr}
\\[-1.8ex]\hline 
\hline\\[-1.8ex]

\multicolumn{3}{c|}{Settings} &  \multicolumn{3}{c|}{Without Routing Flexibility} & \multicolumn{3}{c}{With Routing Flexibility}
\\\hline\\[-1.8ex]

$CV_D$                & $CV_C$                & $\mathbb{E}$[D]/$\mathbb{E}$[C] & ACDPD & MSG  & MSG v.s. ACDPD & ACDPD & MSG  & MSG v.s. ACDPD\\
\hline\\[-1.8ex]

\multirow{4}{*}{0.1} & \multirow{2}{*}{0.1} & 1.0                & 10,028        & 10,115 & $\odot$           & 9,284        & 9,328 & $\odot$        \\
                     &                      & 2.0             & 5,828        & 6,193 & 6.3\%     & 5,126        & 5,684 & 10.9\%                   \\
                     & \multirow{2}{*}{0.4} & 1.0                 & 5,975        & 6,989 & 17.0\%   & 5,197        & 6,262 & 20.5\%               \\
                     &                      & 2.0            & 3,758        & 4,213 & 12.1\%    & 3,282        & 4,193 & 27.8\%                    \\
\hline\\[-1.8ex]
                     
\multirow{4}{*}{0.4} & \multirow{2}{*}{0.1} & 1.0                  & 8,988        & 10,013 & 11.4\%         & 8,275        & 8,868 & 7.2\%              \\
                     &                      & 2.0             & 5,172        & 6,034 & 16.7\%       & 4,487        & 5,496 & 22.5\%            \\
                     & \multirow{2}{*}{0.4} & 1.0                  & 5,944        & 7,089 & 19.3\%       & 5,179        & 6,239 & 20.5\%                          \\
                     &                      & 2.0           & 3,506        & 4,216 & 20.2\%  & 3,172        & 4,076 & 28.5\%\\
\hline
\hline \\[-1.8ex]
\multicolumn{9}{l}{Columns ``ACDPD'' and ``MSG'' are expected revenue. ``MSG v.s. ACDPD'' is a relative revenue increase.}\\
\multicolumn{9}{l}{$\odot$ denotes there is no statistically significant difference between MSG and ACDPD at 95\% confidence level.}\\
\hline
\hline \\[-1.8ex]                      
\end{tabular}
% \begin{tablenotes}
% {\item %\scriptsize 
% Notes: Columns ``ACDPD'' and ``MSG'' are expected revenue. ``MSG v.s. ACDPD'' is a relative revenue increase at 95\% confidence level. $\odot$ denotes there is no statistically significant difference between MSG and ACDPD at 95\% confidence level.\par}
% \end{tablenotes}
\vskip -0.1in
\end{table}

\noindent\textbf{With Routing Flexibility Experimental Setup.} Figure \ref{fig:network_structure} (b) demonstrates a network structure with routing flexibility. Compared to Figure \ref{fig:network_structure} (a), there is an additional leg (link 9) from node 1 to node 3 on top of the spoke-hub network. With the additional link 9, the request from origin 1 to destination 3 can be served with two route options: 1) Route 1: leg 9; 2) Route 2: leg 1 from node 1 to node 5, then leg 7 from node 5 to node 3. We set the average capacity level (hence the total capacity level) the same way in the without-routing case by scaling the capacity levels of leg 1 to leg 8 by 8/9, and adding extra capacity to leg 9.

As shown in Table \ref{tab:air_cargo_result} under ``With Routing Flexibility'' columns, booking limit control outperforms ACDPD by an average 17.22\% among all test instances. An important observation is that booking limit control outperforms ACDPD even more with routing flexibility compared to fixed routes setting. It is not surprising because ACDPD only heuristically deals with the routing decisions by splitting reservations of class $i$ with the same origin and destination equally into fixed $K_i$ classes with different routes. For example, the requests from origin 1 to destination 3 are equally divided into two routes during the reservation stage. In addition, we still observe a similar trend that higher variance leads to a larger performance gap between MSG and ACDPD.

\section{Conclusion and Future Directions}
In this paper, we propose three gradient-based methods for solving a family of stochastic nonconvex optimization \eqref{problem:math_original_extension} to global optimality with non-asymptotic guarantees, in which the complexity of MSG matches the lower bounds. We model air-cargo NRM under booking limit control policy as two-stage stochastic optimization models and as special cases of the proposed model \eqref{problem:math_original_extension} and illustrate the superior performance of proposed algorithms theoretically and numerically.

Much remains open and requires further investigation. 1. When $\xi$ is a positively dependent random vector $\xi$~\citep{chen2019stochastic}, the convex reformulation \eqref{eqn:trans-model-2-extension} does not hold and there exists an infinitely-dimensional stochastic convex reformulation. There is a lack of efficient algorithms for solving such an infinite-dimensional problem despite some statistical results about the asymptotic performance of SAA \citep{deng2022smoothing} and \citep{singham2020sample}. {2. In this paper, our algorithms adjust the booking limit by leveraging the aggregated demand collected after each reservation period. It remains interesting to design an online booking limit control policy that adjusts the booking limit right after accepting or rejecting a demand request.}

\ACKNOWLEDGMENT{This research is partly supported by National Science Foundation Grants CMMI-1761699, CRII-1755829, ZJU-UIUC Institute Research Program, and NCCR Automation of SNSF in Switzerland. }

\bibliographystyle{informs2014}
\bibliography{main}

\begin{thebibliography}{61}
\providecommand{\natexlab}[1]{#1}
\providecommand{\url}[1]{\texttt{#1}}
\providecommand{\urlprefix}{URL }

\bibitem[{Agarwal et~al.(2020)Agarwal, Kakade, Lee, \protect\BIBand{} Mahajan}]{agarwal2020optimality}
Agarwal A, Kakade SM, Lee JD, Mahajan G (2020) Optimality and approximation with policy gradient methods in markov decision processes. \emph{Conference on Learning Theory}, 64--66 (PMLR).

\bibitem[{Agarwal et~al.(2009)Agarwal, Wainwright, Bartlett, \protect\BIBand{} Ravikumar}]{agarwal2009information}
Agarwal A, Wainwright MJ, Bartlett PL, Ravikumar PK (2009) Information-theoretic lower bounds on the oracle complexity of convex optimization. \emph{Advances in Neural Information Processing Systems}, 1--9.

\bibitem[{Agrawal \protect\BIBand{} Devanur(2014)}]{agrawal2014fast}
Agrawal S, Devanur NR (2014) Fast algorithms for online stochastic convex programming. \emph{Proceedings of the twenty-sixth annual ACM-SIAM symposium on Discrete algorithms}, 1405--1424 (SIAM).

\bibitem[{Anderson et~al.(2019)Anderson, Doyle, Low, \protect\BIBand{} Matni}]{anderson2019system}
Anderson J, Doyle JC, Low SH, Matni N (2019) System level synthesis. \emph{Annual Reviews in Control} 47:364--393.

\bibitem[{Balseiro et~al.(2023)Balseiro, Lu, \protect\BIBand{} Mirrokni}]{balseiro2023best}
Balseiro SR, Lu H, Mirrokni V (2023) The best of many worlds: Dual mirror descent for online allocation problems. \emph{Operations Research} 71(1):101--119.

\bibitem[{Barz \protect\BIBand{} Gartner(2016)}]{barz2016air}
Barz C, Gartner D (2016) Air cargo network revenue management. \emph{Transportation Science} 50(4):1206--1222.

\bibitem[{Ben-Tal \protect\BIBand{} Teboulle(1996)}]{ben1996hidden}
Ben-Tal A, Teboulle M (1996) Hidden convexity in some nonconvex quadratically constrained quadratic programming. \emph{Mathematical Programming} 72(1):51--63.

\bibitem[{Bertsekas(2009)}]{bertsekas2009convex}
Bertsekas D (2009) \emph{Convex optimization theory}, volume~1 (Athena Scientific).

\bibitem[{Bhandari \protect\BIBand{} Russo(2019)}]{bhandari2019global}
Bhandari J, Russo D (2019) Global optimality guarantees for policy gradient methods. \emph{arXiv preprint arXiv:1906.01786} .

\bibitem[{Bottou et~al.(2018)Bottou, Curtis, \protect\BIBand{} Nocedal}]{bottou2018optimization}
Bottou L, Curtis FE, Nocedal J (2018) Optimization methods for large-scale machine learning. \emph{Siam Review} 60(2):223--311.

\bibitem[{Brumelle \protect\BIBand{} McGill(1993)}]{brumelle1993airline}
Brumelle SL, McGill JI (1993) Airline seat allocation with multiple nested fare classes. \emph{Operations Research} 41(1):127--137.

\bibitem[{Chen et~al.(2022)Chen, Jiang, Zhang, \protect\BIBand{} Zhou}]{chen2022learning}
Chen B, Jiang J, Zhang J, Zhou Z (2022) Learning to order for inventory systems with lost sales and uncertain supplies. \emph{arXiv preprint arXiv:2207.04550} .

\bibitem[{Chen \protect\BIBand{} Gao(2019)}]{chen2019stochastic}
Chen X, Gao X (2019) Stochastic optimization with decisions truncated by positively dependent random variables. \emph{Operations Research} 67(5):1321--1327.

\bibitem[{Chen et~al.(2015)Chen, Gao, \protect\BIBand{} Hu}]{chen2015new}
Chen X, Gao X, Hu Z (2015) A new approach to two-location joint inventory and transshipment control via l$^\natural$-convexity. \emph{Operations Research Letters} 43(1):65--68.

\bibitem[{Chen et~al.(2018)Chen, Gao, \protect\BIBand{} Pang}]{chen2018preservation}
Chen X, Gao X, Pang Z (2018) Preservation of structural properties in optimization with decisions truncated by random variables and its applications. \emph{Operations Research} 66(2):340--357.

\bibitem[{Chen \protect\BIBand{} Shi(2023)}]{chen2019network}
Chen Y, Shi C (2023) Network revenue management with online inverse batch gradient descent method. \emph{Production and Operations Management} .

\bibitem[{Ciarallo et~al.(1994)Ciarallo, Akella, \protect\BIBand{} Morton}]{ciarallo1994periodic}
Ciarallo FW, Akella R, Morton TE (1994) A periodic review, production planning model with uncertain capacity and uncertain demand—optimality of extended myopic policies. \emph{Management Science} 40(3):320--332.

\bibitem[{Dada et~al.(2007)Dada, Petruzzi, \protect\BIBand{} Schwarz}]{dada2007newsvendor}
Dada M, Petruzzi NC, Schwarz LB (2007) A newsvendor’s procurement problem when suppliers are unreliable. \emph{Manufacturing \& Service Operations Management} 9(1):9--32.

\bibitem[{Davis \protect\BIBand{} Drusvyatskiy(2018)}]{davis2018stochastic}
Davis D, Drusvyatskiy D (2018) Stochastic subgradient method converges at the rate $ o (k^{-1/4}) $ on weakly convex functions. \emph{arXiv preprint arXiv:1802.02988} .

\bibitem[{Deng et~al.(2022)Deng, Xiong, Yang, \protect\BIBand{} Yang}]{deng2022smoothing}
Deng C, Xiong Y, Yang L, Yang Y (2022) A smoothing saa method for solving a nonconvex multisource supply chain stochastic optimization model. \emph{Mathematical Problems in Engineering} 2022.

\bibitem[{Drusvyatskiy \protect\BIBand{} Paquette(2019)}]{drusvyatskiy2019efficiency}
Drusvyatskiy D, Paquette C (2019) Efficiency of minimizing compositions of convex functions and smooth maps. \emph{Mathematical Programming} 178(1):503--558.

\bibitem[{Erdelyi \protect\BIBand{} Topaloglu(2009)}]{erdelyi2009separable}
Erdelyi A, Topaloglu H (2009) Separable approximations for joint capacity control and overbooking decisions in network revenue management. \emph{Journal of Revenue and Pricing Management} 8(1):3--20.

\bibitem[{Erdelyi \protect\BIBand{} Topaloglu(2010)}]{erdelyi2010dynamic}
Erdelyi A, Topaloglu H (2010) A dynamic programming decomposition method for making overbooking decisions over an airline network. \emph{INFORMS Journal on Computing} 22(3):443--456.

\bibitem[{Fatkhullin et~al.(2023)Fatkhullin, He, \protect\BIBand{} Hu}]{fatkhullin2023stochastic}
Fatkhullin I, He N, Hu Y (2023) Stochastic optimization under hidden convexity. \emph{arXiv preprint arXiv:2401.00108} .

\bibitem[{Feng et~al.(2015)Feng, Li, \protect\BIBand{} Shen}]{feng2015air}
Feng B, Li Y, Shen ZJM (2015) Air cargo operations: Literature review and comparison with practices. \emph{Transportation Research Part C: Emerging Technologies} 56:263--280.

\bibitem[{Feng et~al.(2019)Feng, Jia, \protect\BIBand{} Shanthikumar}]{feng2019dynamic}
Feng Q, Jia J, Shanthikumar JG (2019) Dynamic multisourcing with dependent supplies. \emph{Management Science} 65(6):2770--2786.

\bibitem[{Feng \protect\BIBand{} Shanthikumar(2018)}]{feng2018supply}
Feng Q, Shanthikumar JG (2018) Supply and demand functions in inventory models. \emph{Operations Research} 66(1):77--91.

\bibitem[{Feng \protect\BIBand{} Shanthikumar(2022)}]{feng2022applications}
Feng Q, Shanthikumar JG (2022) Applications of stochastic orders and stochastic functions in inventory and pricing problems. \emph{Production and Operations Management} 31(4):1433--1453.

\bibitem[{Ghadimi \protect\BIBand{} Lan(2013)}]{ghadimi2013stochastic}
Ghadimi S, Lan G (2013) Stochastic first-and zeroth-order methods for nonconvex stochastic programming. \emph{SIAM Journal on Optimization} 23(4):2341--2368.

\bibitem[{Ghadimi et~al.(2016)Ghadimi, Lan, \protect\BIBand{} Zhang}]{ghadimi2016mini}
Ghadimi S, Lan G, Zhang H (2016) Mini-batch stochastic approximation methods for nonconvex stochastic composite optimization. \emph{Mathematical Programming} 155(1):267--305.

\bibitem[{Ghai et~al.(2022)Ghai, Lu, \protect\BIBand{} Hazan}]{ghai2022non}
Ghai U, Lu Z, Hazan E (2022) Non-convex online learning via algorithmic equivalence. \emph{Advances in Neural Information Processing Systems} 35:22161--22172.

\bibitem[{Han et~al.(2020)Han, Zhou, \protect\BIBand{} Weissman}]{han2020optimal}
Han Y, Zhou Z, Weissman T (2020) Optimal no-regret learning in repeated first-price auctions. \emph{arXiv preprint arXiv:2003.09795} .

\bibitem[{Hong et~al.(2020)Hong, Wai, Wang, \protect\BIBand{} Yang}]{hong2020two}
Hong M, Wai HT, Wang Z, Yang Z (2020) A two-timescale framework for bilevel optimization: Complexity analysis and application to actor-critic. \emph{arXiv preprint arXiv:2007.05170} .

\bibitem[{Hu et~al.(2020{\natexlab{a}})Hu, Chen, \protect\BIBand{} He}]{hu2020sample}
Hu Y, Chen X, He N (2020{\natexlab{a}}) Sample complexity of sample average approximation for conditional stochastic optimization. \emph{SIAM Journal on Optimization} 30(3):2103--2133.

\bibitem[{Hu et~al.(2021)Hu, Chen, \protect\BIBand{} He}]{hu2021bias}
Hu Y, Chen X, He N (2021) On the bias-variance-cost tradeoff of stochastic optimization. \emph{Advances in Neural Information Processing Systems} 34.

\bibitem[{Hu et~al.(2024)Hu, Wang, Xie, Krause, \protect\BIBand{} Kuhn}]{hu2024contextual}
Hu Y, Wang J, Xie Y, Krause A, Kuhn D (2024) Contextual stochastic bilevel optimization. \emph{Advances in Neural Information Processing Systems} 36.

\bibitem[{Hu et~al.(2020{\natexlab{b}})Hu, Zhang, Chen, \protect\BIBand{} He}]{hu2020biased}
Hu Y, Zhang S, Chen X, He N (2020{\natexlab{b}}) Biased stochastic first-order methods for conditional stochastic optimization and applications in meta learning. \emph{Advances in Neural Information Processing Systems} 33:2759--2770.

\bibitem[{Karaesmen \protect\BIBand{} Van~Ryzin(2004)}]{karaesmen2004overbooking}
Karaesmen I, Van~Ryzin G (2004) Overbooking with substitutable inventory classes. \emph{Operations Research} 52(1):83--104.

\bibitem[{Karimi et~al.(2016)Karimi, Nutini, \protect\BIBand{} Schmidt}]{karimi2016linear}
Karimi H, Nutini J, Schmidt M (2016) Linear convergence of gradient and proximal-gradient methods under the polyak-{\l}ojasiewicz condition. \emph{Joint European Conference on Machine Learning and Knowledge Discovery in Databases}, 795--811 (Springer).

\bibitem[{Klein et~al.(2020)Klein, Koch, Steinhardt, \protect\BIBand{} Strauss}]{klein2020review}
Klein R, Koch S, Steinhardt C, Strauss AK (2020) A review of revenue management: Recent generalizations and advances in industry applications. \emph{European Journal of Operational Research} 284(2):397--412.

\bibitem[{Kleywegt et~al.(2002)Kleywegt, Shapiro, \protect\BIBand{} Homem-de Mello}]{kleywegt2002sample}
Kleywegt AJ, Shapiro A, Homem-de Mello T (2002) The sample average approximation method for stochastic discrete optimization. \emph{SIAM Journal on Optimization} 12(2):479--502.

\bibitem[{Kunnumkal \protect\BIBand{} Topaloglu(2008)}]{kunnumkal2008using}
Kunnumkal S, Topaloglu H (2008) Using stochastic approximation methods to compute optimal base-stock levels in inventory control problems. \emph{Operations Research} 56(3):646--664.

\bibitem[{Li \protect\BIBand{} Pang(2017)}]{li2017dynamic}
Li D, Pang Z (2017) Dynamic booking control for car rental revenue management: A decomposition approach. \emph{European Journal of Operational Research} 256(3):850--867.

\bibitem[{Li \protect\BIBand{} Ye(2022)}]{li2022online}
Li X, Ye Y (2022) Online linear programming: Dual convergence, new algorithms, and regret bounds. \emph{Operations Research} 70(5):2948--2966.

\bibitem[{Miao \protect\BIBand{} Wang(2021)}]{miao2021network}
Miao S, Wang Y (2021) Network revenue management with nonparametric demand learning: $\sqrt{T}$-regret and polynomial dimension dependency. \emph{Available at SSRN 3948140} .

\bibitem[{Miao \protect\BIBand{} Wang(2024)}]{miao2024demand}
Miao S, Wang Y (2024) Demand balancing in primal-dual optimization for blind network revenue management. \emph{arXiv preprint arXiv:2404.04467} .

\bibitem[{Nemirovski et~al.(2009)Nemirovski, Juditsky, Lan, \protect\BIBand{} Shapiro}]{nemirovski2009robust}
Nemirovski A, Juditsky A, Lan G, Shapiro A (2009) Robust stochastic approximation approach to stochastic programming. \emph{SIAM Journal on optimization} 19(4):1574--1609.

\bibitem[{Previgliano \protect\BIBand{} Vulcano(2021)}]{previgliano2021managing}
Previgliano F, Vulcano G (2021) Managing uncertain capacities for network revenue optimization. \emph{Manufacturing \& Service Operations Management} .

\bibitem[{Sakos et~al.(2024)Sakos, Vlatakis-Gkaragkounis, Mertikopoulos, \protect\BIBand{} Piliouras}]{sakos2024exploiting}
Sakos I, Vlatakis-Gkaragkounis EV, Mertikopoulos P, Piliouras G (2024) Exploiting hidden structures in non-convex games for convergence to nash equilibrium. \emph{Advances in Neural Information Processing Systems} 36.

\bibitem[{Shaked \protect\BIBand{} Shanthikumar(2007)}]{shaked2007stochastic}
Shaked M, Shanthikumar JG (2007) \emph{Stochastic orders} (Springer Science \& Business Media).

\bibitem[{Singham \protect\BIBand{} Lam(2020)}]{singham2020sample}
Singham DI, Lam H (2020) Sample average approximation for functional decisions under shape constraints. \emph{2020 Winter Simulation Conference (WSC)}, 2791--2799 (IEEE).

\bibitem[{Starr(1969)}]{starr1969quasi}
Starr RM (1969) Quasi-equilibria in markets with non-convex preferences. \emph{Econometrica: journal of the Econometric Society} 25--38.

\bibitem[{Sun(2021)}]{sun2023nonconvex}
Sun J (2021) Provable nonconvex methods/algorithms \urlprefix\url{https://sunju.org/research/nonconvex/}.

\bibitem[{Talluri \protect\BIBand{} Van~Ryzin(1998)}]{talluri1998analysis}
Talluri K, Van~Ryzin G (1998) An analysis of bid-price controls for network revenue management. \emph{Management Science} 44(11-part-1):1577--1593.

\bibitem[{Tang \protect\BIBand{} Kouvelis(2014)}]{tang2014pay}
Tang SY, Kouvelis P (2014) Pay-back-revenue-sharing contract in coordinating supply chains with random yield. \emph{Production and Operations Management} 23(12):2089--2102.

\bibitem[{Vapnik(1999)}]{vapnik1999nature}
Vapnik V (1999) \emph{The nature of statistical learning theory} (Springer science \& business media).

\bibitem[{Wang et~al.(2021)Wang, Meng, Wang, \protect\BIBand{} Qu}]{wang2021two}
Wang T, Meng Q, Wang S, Qu X (2021) A two-stage stochastic nonlinear integer-programming model for slot allocation of a liner container shipping service. \emph{Transportation Research Part B: Methodological} 150:143--160.

\bibitem[{Wang(2016)}]{wang2016optimal}
Wang X (2016) Optimal allocation of limited and random network resources to discrete stochastic demands for standardized cargo transportation networks. \emph{Transportation Research Part B: Methodological} 91:310--331.

\bibitem[{Ying et~al.(2023)Ying, Guo, Ding, Lavaei, \protect\BIBand{} Shen}]{ying2023policy}
Ying D, Guo MA, Ding Y, Lavaei J, Shen ZJ (2023) Policy-based primal-dual methods for convex constrained markov decision processes. \emph{Proceedings of the AAAI Conference on Artificial Intelligence}, volume~37, 10963--10971.

\bibitem[{Yuan et~al.(2021)Yuan, Luo, \protect\BIBand{} Shi}]{yuan2021marrying}
Yuan H, Luo Q, Shi C (2021) Marrying stochastic gradient descent with bandits: Learning algorithms for inventory systems with fixed costs. \emph{Management Science} 67(10):6089--6115.

\bibitem[{Zhang et~al.(2020)Zhang, Koppel, Bedi, Szepesvari, \protect\BIBand{} Wang}]{zhang2020variational}
Zhang J, Koppel A, Bedi AS, Szepesvari C, Wang M (2020) Variational policy gradient method for reinforcement learning with general utilities. \emph{Advances in Neural Information Processing Systems} 33:4572--4583.

\end{thebibliography}

\newpage

\renewcommand{\theHsection}{A\arabic{section}}

\begin{APPENDICES}

\begin{center}
\textbf{\Large Online Appendices}
\end{center}
~\\
\section*{Organization of Appendices}
The appendices are organized as follows.
In Appendix \ref{sec:technical_details}, we show the technical details on the proof of Lemmas \ref{lemma:gradient_advanced}, \ref{lm:properties_of_inverse_expectation_estimation}, \ref{lm:switch_transformation_projection}, \ref{lemma:verify_x_wedge_xi}, the second moment of the gradient estimators in MSG, and the auxiliary results related to stationary convergence of RSG. In Appendix \ref{sec:proof_of_RSG}, we demonstrate the analysis of the global convergence of RSG. In Appendix \ref{sec:proof_of_msg}, we demonstrate the analysis of the global convergence of MSG. In Appendix \ref{appendix:discussion_on_assupmtions}, we discuss the conditions required for operations management applications to ensure that the assumptions needed by the global convergence results of RSG and MSG. We also discuss situations when certain assumptions do not hold and what happens to the practical performance of the proposed algorithm. In Appendix \ref{sec:sgd_u}, we discuss the detailed construction of SAA+SG method, which builds an empirical convex reformulation via SAA and solves the empirical convex reformulation via SGD, and the sample and gradient complexities of SAA+SG. In Appendix \ref{appendix:application_ATO}, we give a model formulation of assemble-to-order systems as a special case of our nonconvex optimization problem. In Appendix \ref{sec:appendix_NRM_alg}, we further discuss the details of the NRM problem given in Section \ref{sec:application}, including the modeling of the passenger NRM, concavity of the NRM models, computing the stochastic gradient of $f$ in NRM problems, and discussions about integer booking limits and Poisson show-ups in NRM problems.  In Appendix \ref{appendix:numerics}, we discussed details of the numerical implementation and demonstrate the full numerical results.

\section{Technical Details}\label{sec:technical_details}
\subsection{Proof of Lemma \ref{lemma:gradient_advanced}}
\label{proof_of_closed_form_gradient}

By Assumption \ref{assumption:general}(d), for any $\xi\in\Xi$, the probability that $\phi(x,\xi)$ is non-differentiable in $x$ is zero. In addition, $\phi$ is Lipschitz continuous in $x\in\mathcal{X}$ for any given $\xi$. Without loss of generality, for a given $\xi\in\Xi$, we define
\begin{equation}
\label{eq:redefined_nabla_phi}   
\bar \nabla \phi(x,\xi) = 
\begin{cases}
\nabla \phi(x,\xi) & \text{ if } \phi(x,\xi)  \text{ is differentiable in } x,\\
0 & \text{ otherwise.}
\end{cases}
\end{equation}
For simplicity, we shall use $\nabla \phi(x,\xi)$ and $\bar \nabla \phi(x,\xi)$ indifferently. 

\begin{proof1}
{\color{black} 
Since $\phi$ is $L_\phi$-Lipschitz continuous, it holds for any $x\in\mathcal{X}$,  $\xi\in\Xi$, $i\in[d]$, and $h\not = 0$ that
$$
\Big\|\frac{\phi_i(x_i+h,\xi_i) - \phi_i(x_i+h,\xi_i)}{h} \Big\|\leq L_\phi.
$$
Since $\xi$ is component-wise independent and $\phi(x,\xi) = (\phi_1(x_1,\xi_1),...,\phi_d(x_d,\xi_d))^\top$, without loss of generality, let us consider the first coordinate $\nabla g_1(x_1)$. The other coordinates follow directly.
$$
\begin{aligned}
    \nabla g_1(x_1) 
= & 
    \nabla \EE_{\xi_1} [\phi_1(x_1,\xi_1)]   \\
= & 
    \lim_{h\rightarrow 0} \int_{t\in\Xi_1} \frac{\phi_1(x_1+h,t) - \phi_1(x_1,t)}{h} dH_1(t)\\
= &
    \lim_{h\rightarrow 0} \int_{t \in \Theta}\frac{ \phi_1(x_1+h,t)  - \phi_1(x_1,t)}{h} d H_1(t) + \lim_{h\rightarrow 0} \int_{t \in \Theta^c}\frac{ \phi_1(x_1+h,t)  - \phi_1(x_1,t)}{h} d H_1(t)\\
=   &
      \int_{t \in \Theta}\lim_{h\rightarrow 0} \frac{ \phi_1(x_1+h,t)  - \phi_1(x_1,t)}{h} d H_1(t) + \lim_{h\rightarrow 0} \int_{t \in \Theta^c}\frac{ \phi_1(x_1+h,t)  - \phi_1(x_1,t)}{h} d H_1(t)\\
=   &
    \EE_{\xi_1} [\nabla \phi_1(x_1,\xi_1)|\Theta] \PP(\Theta)+ \lim_{h\rightarrow 0} \int_{t \in \Theta^c}\frac{ \phi_1(x_1+h,t)  - \phi_1(x_1,t)}{h} d H_1(t)\\
= &
    \EE_{\xi_1} [\nabla \phi_1(x_1,\xi_1)|\Theta] \PP(\Theta)\\
= &
    \EE_{\xi_1} [ \nabla \phi_1(x_1,\xi_1)],
\end{aligned}
$$
where $\Xi_1$ denotes the support of $\xi_1$, the event $\Theta:=\{\xi_1\in\Xi_1\mid \phi_1(x_1,\xi_1) \text{ is differentiable in }  x_1 \}$ and $\Theta^c$ denotes the complement of $\Theta$, the second equality holds by definition of the derivative, the third equality holds naturally, the forth equality holds by  dominated convergence theorem and mean-value theorem that one could switch the order of limit and integration as $\phi$ is Lipschitz continuous, the fifth equality holds by the definitions of derivative and conditional expectation, and the sixth equality holds as for any given $h\not =0$,
$
 \frac{\phi_1(x_1+h,t) - \phi_1(x_1,t) }{h}
$
is uniformly upper and lower bounded by the  Lipschitz continuous parameter $L_\phi$ and $P(\Theta^c) = 0$ by Assumption \ref{assumption:general}(d). By \eqref{eq:redefined_nabla_phi}, the last equality holds.

Since $f$ is continuously differentiable and $L_f$-Lipschitz continuous,
by Assumption \ref{assumption:general}(b)(c), it holds that
$$
\|\nabla\phi(x,\xi)^\top \nabla f(\phi(x,\xi)) \|\leq \|\nabla\phi(x,\xi)\| \|\nabla f(\phi(x,\xi))\|\leq L_\phi L_f.
$$
Following a similar argument, we have 
$$
\begin{aligned}
    \nabla F(x) 
= & 
    \nabla_x  \EE_\xi [f(\phi(x,\xi))]
=  
     \EE_\xi \nabla_x [f(\phi(x,\xi))] 
= 
    \EE_{\xi} [\nabla \phi(x,\xi)^\top\nabla f(\phi(x, \xi))].
\end{aligned}
$$
}
By Assumption \ref{assumption:general_2}(b), we have $\nabla g(x)\succeq \mu_g I$ for any $x\in\mathcal{X}$. By the inverse function theorem, we have
$$
\nabla g^{-1}(u) =  [\nabla g(g^{-1}(u))]^{-1}.
$$
Since $G(u) = F(g^{-1}(x))$, by the chain rule, it holds that
\begin{equation}
\begin{aligned}
    \nabla G(u) 
& = 
    \nabla g^{-1}(u)^\top\nabla F(g^{-1}(u))\\
& = 
    \nabla g^{-1}(u)^\top\EE_{\xi} [\nabla \phi(g^{-1}(u),\xi)^\top\nabla f(\phi(g^{-1}(u), \xi)) ]\\
& = 
    [\nabla g(g^{-1}(u))]^{-\top}\EE_{\xi} [\nabla \phi(g^{-1}(u),\xi)^\top\nabla f(\phi(g^{-1}(u), \xi)) ]
\end{aligned}    
\end{equation}
which completes the proof.
\end{proof1}

\subsection{Proof of Lemma \ref{lm:properties_of_inverse_expectation_estimation}: Estimating Matrix Inverse}
\label{proof:inverse_estimation}

\textbf{Remark:}  
Note that the distribution of $k$ in Lemma \ref{lm:properties_of_inverse_expectation_estimation} is a uniform distribution over the support $\{0,\ldots,K-1\}$. One could use other distributions to build up estimators for matrix inverse. For instance, when using a geometric distribution with parameter $p$ over support $\{0,1,\ldots,\infty\}$,  the estimator is 
$$
[\nabla \tilde g(x)]^{-1} = 
\begin{cases}
\frac{1}{p^k(1-p) c L_\phi}\prod_{i=1}^{k} \Big(I - \frac{\nabla \phi(x,\xi^{i})}{cL_\phi}\Big) & \text{ if } k\geq 1;\\
\frac{1}{p^k (1-p) c L_\phi} I & \text{ if } k=0.
\end{cases}
$$
One can show that $[\nabla \tilde g(x)]^{-1}$ is unbiased, has a bounded second moment, and needs $\cO(1)$ number of samples in expectation to construct. However, with a small probability, the estimator $[\nabla \tilde g(x)]^{-1}$ could have very large entries. 
\begin{proof1}
We first bound the bias. To simplify notation, we let $\prod_{i=1}^{k} \Big(I - \frac{\nabla \phi(x,\xi^{i})}{cL_\phi}\Big)=I$ for $k=0$.
\begin{align*}
    \EE [\nabla \hat g(x)]^{-1}
= &   
    \EE_k \EE_{\{\xi^{i}\}_{i=1}^k}\frac{K}{cL_\phi} \prod_{i=1}^{k} \Big(I - \frac{\nabla \phi(x,\xi^{i})}{cL_\phi}\Big)
= 
    \EE_k \frac{K}{cL_\phi} \prod_{i=1}^{k} \Big(I - \frac{\EE_{\xi^{i}}\nabla \phi(x,\xi^{i})}{cL_\phi}\Big)\\
= &
    \EE_k \frac{K}{cL_\phi} \prod_{i=1}^{k} \Big(I - \frac{\EE_{\xi}\nabla \phi(x,\xi)}{cL_\phi}\Big)
= 
    \frac{K}{cL_\phi}\frac{1}{K}\sum_{k=0}^{K-1} \Big(I - \frac{\EE_{\xi}\nabla \phi(x,\xi)}{cL_\phi}\Big)^k \\
= &     
    \frac{1}{cL_\phi}\sum_{k=0}^{K-1} \Big(I - \frac{\nabla g(x)}{cL_\phi}\Big)^k,
\end{align*} 
where the last equality holds as $0\preceq\nabla\phi(x,\xi)\preceq L_\phi I$ and dominated convergence theorem guarantees interchange of expectation and gradient. 
On the other hand, since $I\succeq \frac{\nabla g(x)}{cL_\phi}\succeq \frac{\mu_g}{cL_\phi}$ for any $x$, we have
$$
[\nabla g(x)]^{-1}  = \frac{1}{cL_\phi}\sum_{k=0}^\infty \Big(I-\frac{\nabla g(x)}{cL_\phi}\Big)^k.
$$
As a result, the bias of the estimator is upper bounded.
$$
\begin{aligned}
    \|\EE [\nabla \hat g(x)]^{-1}-[\nabla g(x)]^{-1}\| 
= &
    \Big\|\frac{1}{cL_\phi}\sum_{k=K}^{\infty} \Big(I - \frac{\nabla g(x)}{cL_\phi}\Big)^k\Big\|
\leq 
    \frac{1}{cL_\phi}\sum_{k=K}^{\infty} \Big\|\Big(I - \frac{\nabla g(x)}{cL_\phi}\Big)^k\Big\|\\
\leq & 
    \frac{1}{cL_\phi}\Big\|\Big(I - \frac{\nabla g(x)}{cL_\phi}\Big)^K\Big\|~~\sum_{k=0}^{\infty} \Big\|\Big(I - \frac{\nabla g(x)}{cL_\phi}\Big)^k\Big\|\\
\leq &
    \frac{1}{cL_\phi}\Big\|I - \frac{\nabla g(x)}{cL_\phi}\Big\|^K~~\sum_{k=0}^{\infty} \Big\|I - \frac{\nabla g(x)}{cL_\phi}\Big\|^k\\
= &
    \frac{1}{cL_\phi}\Big\|I - \frac{\nabla g(x)}{cL_\phi}\Big\|^K \frac{1}{1 -\Big\|I - \frac{\nabla g(x)}{cL_\phi}\Big\| }  
\leq 
    \frac{1}{\mu_g}\Big(1- \frac{\mu_g}{cL_\phi}\Big)^K,
\end{aligned}
$$
where the first inequality holds by triangle inequality, the second and third inequality holds by spectral norm, the second equality holds as $0\prec I - \frac{\nabla g(x)}{cL_\phi}\prec I$,  and the last inequality holds by Assumptions \ref{assumption:general_2}(b). 
As for the second moment, since $0\preceq\nabla\phi(x,\xi)\preceq L_\phi I$, we have 
$$
\begin{aligned}
    \EE \|[\nabla \hat g(x)]^{-1}\|^2
\leq &
    \EE_k  \EE_{\{\xi^{i}\}_{i=1}^k}\Big[ \frac{K^2}{c^2L_\phi^2} \prod_{i=1}^{k} \Big\|I - \frac{\nabla \phi(x,\xi^{i})}{cL_\phi}\Big\|^2\Big]
\leq 
     \frac{K^2}{c^2L_\phi^2} \EE_k \prod_{i=1}^{k} \|I\|^2\\
= & 
    \frac{K^2}{c^2L_\phi^2},
\end{aligned}
$$
where the first inequality holds by spectral norm, and the second inequality holds as $c>1$.

The average number of samples used to construct the estimator is 
$$
\EE_k k = \frac{1}{K}\sum_{k=0}^{K-1} k = \frac{(K-1)}{2}.
$$
\end{proof1}

\subsection{Proof of Lemma \ref{lm:switch_transformation_projection}: Switching Projection and Transformation}
\label{sec:switch_projection_transformation}
\begin{proof1}
By definition, $\mathcal{X}  = [\munderbar{X}_1,\bar{X}_1]\times\ldots\times [\munderbar{X}_d,\bar{X}_d]$, $\mathcal{U}  = [\EE[\phi_1(\munderbar{X}_1,\xi_1)], \EE[\phi_1(\bar{X}_1,\xi_1)]]\times\ldots\times [\EE[\phi_d(\munderbar{X}_d,\xi_d)], \EE[\phi_d(\bar{X}_d,\xi_d)]]:=[\munderbar{U}_1,\bar{U}_1]\times\ldots\times [\munderbar{U}_d,\bar{U}_d]$.

It suffices to show the one-dimensional case because $\phi(x,\xi) = (\phi_1(x_1,\xi_1),\ldots, \phi_d(x_d,\xi_d))^\top$, $\xi$ is component-wise independent, and both $\mathcal{X}$ and $\mathcal{U}$ are box constraints. Without loss of generality, we denote $\mathcal{X} = [\munderbar{X},\bar{X}]$ and $\mathcal{U} = [\munderbar{U},\bar{U}]$.

\textbf{Case I}: if $x\in\mathcal{X}$, then $g(x)\in\mathcal{U}$. It holds that
$$
g\Big(\Pi_\mathcal{X}(x)\Big) = g(x) = \Pi_\mathcal{U}\Big(g(x)\Big).
$$

\textbf{Case II}: if $x\leq \munderbar{X}$ It holds that
$$
g\Big(\Pi_\mathcal{X}(x)\Big) = \EE[\phi(\Pi_\mathcal{X}(x),\xi)]=\EE[\phi(\munderbar{X},\xi)]=\munderbar{U}.
$$
Since $\Pi_\mathcal{U}\Big(g(x)\Big)\in\mathcal{U}$, it holds that
$$
\Pi_\mathcal{U}\Big(g(x)\Big) \geq \munderbar{U}.
$$
On the other hand, since projection from $\RR$ to an interval is a  non-decreasing  function and $\phi(x,\xi)$ is also non-decreasing for any $\xi$, we have the following,
$$
\Pi_\mathcal{U}\Big(g(x)\Big) = \Pi_\mathcal{U}\Big(\EE[\phi(x,\xi)]\Big) \leq \Pi_\mathcal{U}\Big(\EE[\phi(\munderbar{X},\xi)]\Big) = \Pi_\mathcal{U}\Big(\munderbar{U}\Big)=\munderbar{U}.
$$
As a result, it holds that, $g\Big(\Pi_\mathcal{X}(x)\Big) = \munderbar{U} = \Pi_\mathcal{U}\Big(g(x)\Big)$

\textbf{Case III}: if $x\geq \bar{X}$. It holds that
$$
g\Big(\Pi_\mathcal{X}(x)\Big) = \EE[\phi(\Pi_\mathcal{X}(x),\xi)]=\EE[\phi(\bar{X},\xi)]=\bar{U}.
$$
Since $\Pi_\mathcal{U}\Big(g(x)\Big)\in\mathcal{U}$, we have the following inequality,
$$
\Pi_\mathcal{U}\Big(g(x)\Big) \leq \bar{U}.
$$
On the other hand, due to the non-decreasing property of box projection and $\phi(x,\xi)$ for any $\xi$, we have the following,
$$
\Pi_\mathcal{U}\Big(g(x)\Big) = \Pi_\mathcal{U}\Big(\EE[\phi(x,\xi)]\Big) \geq \Pi_\mathcal{U}\Big(\EE[\phi(\bar{X},\xi)]\Big) = \Pi_\mathcal{U}\Big(\bar{U}\Big)=\bar{U}.
$$
Thus, $g\Big(\Pi_\mathcal{X}(x)\Big) = \bar{U} = \Pi_\mathcal{U}\Big(g(x)\Big)$.
Summarizing all the cases, we obtain the desired result.
\end{proof1}

\subsection{Second Moments of Gradient Estimators $v_F$ and $v_G$ in MSG}
The following lemma characterizes the second moments of gradient estimators $v_F(x)$ and $v_G(u)$ used in the analysis of MSG.
\begin{lemma}
\label{lm:bounded_second_moment}
Under Assumption \ref{assumption:general},  the second moment of $v_G(u)$ and $v_F(x)$ are bounded with
$$
\EE \|v_G(x)\|^2 \leq  \frac{K^2L_f^2}{4} .
$$
$$
\EE \|v_F(u)\|^2\leq \frac{K^4 L_f^2}{8L_\phi^2}+2\lambda^2D_{\mathcal{X}}^2.
$$
\end{lemma}

\begin{proof1}
By Lemma \ref{lm:properties_of_inverse_expectation_estimation}, we have $\EE \|[\nabla \hat g^{j}(g^{-1}(u))]^{-1}\|^2\leq \frac{K^2}{4L_\phi^2}$ for $j=A,B$ and $c=2$. It holds that
\begin{align*}
    \EE \|v_G(u)\|^2
= &
    \EE \|[\nabla \hat g^{B}(g^{-1}(u))]^{-\top} \nabla \phi(g^{-1}(u), \xi)^\top \nabla f( \phi(g^{-1}(u), \xi))\|^2 \\
\leq &
    \EE \|[\nabla \hat g^{B}(g^{-1}(u))]^{-1}\|^2 ~\EE [ \| \nabla \phi(g^{-1}(u), \xi)\|^2 \|\nabla f( \phi(g^{-1}(u), \xi))\|^2 ]\\
\leq &
     \frac{K^2}{4L_\phi^2} L_\phi^2 L_f^2\\
= &
    \frac{K^2 L_f^2}{4},
\end{align*}
where the first inequality uses the Cauchy-Schwarz inequality and the fact that  $[\nabla \hat g^{B}(g^{-1}(u))]^{-1}$ is independent of $\nabla \phi(g^{-1}(u), \xi)^\top \nabla f( \phi(g^{-1}(u), \xi))$, and the second inequality holds by Lemma \ref{lm:properties_of_inverse_expectation_estimation} and Lipschitz continuity of $\phi$ and $f$.

As for $v_F(x)$, we have
\begin{align*}
    \EE \|v_F(x)\|^2
= &
    \EE \|[\nabla \hat g^{A}(x)]^{-\top}[\nabla \hat g^{B}(x)]^{-\top} \nabla \phi(x, \xi)^\top \nabla f( \phi(x, \xi))+\lambda x\|^2 \\
\leq &
    2\EE \|[\nabla \hat g^{A}(x)]^{-\top}[\nabla \hat g^{B}(x)]^{-\top} \nabla \phi(x, \xi)^\top \nabla f( \phi(x, \xi))\|^2+2\EE\|\lambda x\|^2\\
\leq &
     2\EE \|[\nabla \hat g^{A}(x)]^{-1}\|^2 \EE \|[\nabla \hat g^{B}(x)]^{-1}\|^2 \EE \| \nabla \phi(x, \xi)^\top \nabla f( \phi(x, \xi))\|^2 + 2\lambda^2 \EE \| x\|^2\\
\leq &
    2\frac{K^4}{16L_\phi^4}L_\phi^2 L_f^2 + 2\lambda^2 D_\mathcal{X}^2\\
= &
    \frac{K^4 L_f^2}{8L_\phi^2}+ 2\lambda^2 D_\mathcal{X}^2,
\end{align*}
where the first inequality uses the Cauchy-Schwarz inequality, and the second inequality uses the Cauchy-Schwarz inequality and the fact that $[\nabla \hat g^{A}(x)]^{-1}$, $[\nabla \hat g^{B}(x)]^{-1}$, and $\nabla \phi(x, \xi)^\top \nabla f( \phi(x, \xi))$ are independent. 
\end{proof1}

\subsection{Auxiliary Results Related to Gradient Mapping}
We first restate the definition of gradient mapping for constrained optimization problems. For a smooth objective $F$ over a convex domain $\mathcal{X}$, define $\tilde x:= \Pi_{\mathcal{X}}(x-\alpha\nabla F(x))$ for some $\alpha>0$. The definition of gradient mapping of $F$ is given as
\begin{equation}
\label{eq:gradient_mapping}
\tilde \nabla F_\alpha( x ) := \frac{x -\tilde x}{\alpha}.    
\end{equation}

The following lemmas are used in the proof of Theorem \ref{thm:RSG}. In particularly, Lemma \ref{lm:gradient_generalized_gradient} establishes the optimality gap and  the gradient mapping. Lemma \ref{lm:generalied_gradient_stationary_convergence} characterizes the convergence rate, measured in terms of the norm of gradient mapping, of projected SGD on weakly convex objectives. 

\begin{lemma}
\label{lm:gradient_generalized_gradient}
For a convex function $G$ over a convex domain $\mathcal{U}$ with $u\in\mathcal{U}$ and
$\tilde u=\Pi_\mathcal{U}(u - \alpha\nabla G(u))$ for any $\alpha>0$, it holds for any $u^*\in\mathcal{U}$  that 
$$
G(u) - G(u^*)\leq 
(u -  u^*)^\top\frac{u-\tilde u}{\alpha} +   \nabla G(u)^\top (u-\tilde u).
$$
\end{lemma}
\begin{proof1}
Define 
$$
h_\mathcal{U}(u) = 
\begin{cases}
0 & \text{ if } u\in\mathcal{U}, \\
\infty & \text{ otherwise}.
\end{cases}
$$ 
Equivalently, we may rewrite
$$
\tilde u = \argmin_{u^\prime\in\RR^d} \frac{1}{2\alpha}\|u^\prime - (u - \alpha\nabla G(u))\|^2 + h_\mathcal{U}(u^\prime).
$$
By the first-order optimality condition, it holds that
$$
\frac{u-\tilde u}{\alpha} -  \nabla G(u) \in \partial h_\mathcal{U}(\tilde u),
$$
where $\partial$ denotes the subdifferential set of the convex function $h_\mathcal{U}$ at $\tilde u$. 
By definition of $h_\mathcal{U}(u)$, $h_\mathcal{U}(\tilde u) = h_\mathcal{U}(u^*)=0$, it holds that
\begin{align*}
&
    G(u) - G(u^*)\\
= &
    G(u)  - G(u^*) + h_\mathcal{U}(\tilde u) - h_\mathcal{U}(u^*)\\
\leq &
    \nabla G(u)^\top (u-u^*) + \Big(\frac{u-\tilde u}{\alpha} -  \nabla G(u)\Big)^\top(\tilde u-u^*)\\
= &
    \nabla G(u)^\top (u-u^*) + \Big(\frac{u-\tilde u}{\alpha} -  \nabla G(u)\Big)^\top\Big(u - \alpha\frac{u-\tilde u}{\alpha} - u^*\Big)\\
= &
    (u -  u^*)^\top\frac{u-\tilde u}{\alpha} - \alpha\Big\|\frac{u-\tilde u}{\alpha}\Big\|^2 + \alpha\nabla G(u)^\top \frac{u-\tilde u}{\alpha} \\
\leq & 
    (u -  u^*)^\top\frac{u-\tilde u}{\alpha} +  \nabla G(u)^\top (u-\tilde u).
\end{align*}
where the first inequality uses convexity of $G$ and $h_\mathcal{U}(u)$
\end{proof1}

Consider the general stochastic optimization problem:
$$
\min_{x\in\mathcal{X}} \varphi(x):= \EE_\xi [ \Phi(x,\xi)],
$$
where $\mathcal{X}$ is a convex set.
Recall the projected SGD updates with a independent random sample $\xi^t$ and stepsize $\gamma$: 
$$
x^{t+1} = \Pi_\mathcal{X}(x^t-\gamma \nabla \Phi(x^t,\xi^t)).
$$
Let $\hat x^T$ be uniformly selected from $\{x^t\}_{t=1}^T$.
Denote 
$$
\tilde \varphi_\alpha(x):= \min_{y\in\mathcal{X}} \{\varphi(y)+\frac{1}{2\alpha} \|y-x\|^2\},
$$
$$
\mathrm{prox}_{\alpha \varphi}(x) := \argmin_{y\in\mathcal{X}} \{\varphi(y)+\frac{1}{2\alpha} \|y-x\|^2\}.
$$
Function $\tilde \varphi_\alpha$ is the Moreau envelop of $\varphi$ and is widely used in stationary convergence of nonconvex functions~\citep{davis2018stochastic,hu2020biased,drusvyatskiy2019efficiency}. By \citet{davis2018stochastic}, the gradient of the Moreau envelop is given by
$$
\nabla \tilde \varphi_\alpha(x) = \frac{x-\mathrm{prox}_{\alpha \varphi}(x)}{\alpha}.
$$

\begin{lemma}
\label{lm:generalied_gradient_stationary_convergence}
If $\Phi(x,\xi)$ is $L$-Lipschitz continuous in $x$ for any given $\xi$ and $ \nabla \varphi(x)$ is $S$-Lipschitz continuous in $x\in\mathcal{X}$,  the output of projected SGD with stepsize $\gamma=1/\sqrt{T}$ satisfies the following inequality:
$$
    \EE\|\tilde \nabla \varphi_{1/S} (x) \|^2
\leq
    \frac{9}{2}\Big(1+\frac{1}{\sqrt{2}}\Big)^2\frac{(\varphi_{(1/2S)}(x_1)-\min_{x\in\mathcal{X} }\varphi(x))+S L^2}{\sqrt{T}}.    
$$
\end{lemma}
\begin{proof1}
Let $\alpha=1/S$.
Theorem 4.5 and equation (4.9) in \citet{drusvyatskiy2019efficiency} showed that
$$
\frac{1}{4}\Big\|\frac{x-\mathrm{prox}_{(\alpha/2) \varphi}(x)}{\alpha/2}\Big\|\leq \|\tilde \nabla \varphi_\alpha (x) \|\leq \frac{3}{2}\Big(1+\frac{1}{\sqrt{2}}\Big) \Big\|\frac{x-\mathrm{prox}_{(\alpha/2) \varphi}(x)}{\alpha/2}\Big\|.
$$
In addition, Corollary 2.2 in  \citet{davis2018stochastic} showed that the output of projected SGD with stepsize $\gamma = \frac{1}{\sqrt{T}}$ on $\varphi$  satisfies
$$
    \EE \Big\|\frac{\hat x^T-\mathrm{prox}_{ (\alpha/2)\varphi}(\hat x^T)}{\alpha/2}\Big\|^2
\leq
    2\frac{(\varphi_{(\alpha/2)}(x^1)-\min_{x\in\mathcal{X} }\varphi(x))+S L^2}{\sqrt{T}}.
$$
Combining the above two inequalities, we have
\begin{equation}
\label{eq:stationary_convergence_PSGD}
\begin{aligned}
    \EE\|\tilde \nabla \varphi_{1/S} (\hat x^T) \|^2
\leq &
    \frac{9}{4}\Big(1+\frac{1}{\sqrt{2}}\Big)^2 \EE \Big\|\frac{\hat x^T-\mathrm{prox}_{(\alpha/2) \varphi}(\hat x^T)}{\alpha/2}\Big\|^2\\
\leq &
    \frac{9}{2}\Big(1+\frac{1}{\sqrt{2}}\Big)^2\frac{(\varphi_{(1/2S)}(x^1)-\min_{x\in\mathcal{X} }\varphi(x))+S L^2}{\sqrt{T}}.  
\end{aligned}
\end{equation}
\end{proof1}

\section{Proof of Theorem \ref{thm:RSG}: Global Convergence of RSG}
\label{sec:proof_of_RSG}
\begin{proof1}
Recall that $\hat x^T$ is the output of RSG. By definition of gradient mapping given in \eqref{eq:gradient_mapping}, we have
$$
\tilde \nabla F_\alpha(\hat x^T) = \frac{\hat x^T-\Pi_{\mathcal{X}}(\hat x^T-\alpha\nabla F(\hat x^T) )}{\alpha}.
$$
RSG is equivalent to projected SGD on the regularized objective $ F^\lambda(x) = F(x)+\frac{\lambda}{2}\|x\|^2$. Since $F^\lambda$ is $(L_\phi L_f+\lambda D_\mathcal{X})$-Lipschitz continuous and $(S_F+\lambda)$-weakly convex,  Lemma \ref{lm:generalied_gradient_stationary_convergence} implies that $\hat x^T$, the output of RSG with stepsize $\gamma=1/\sqrt{T}$, satisfies
$$
    \EE \|\tilde \nabla F_{\alpha}^\lambda(\hat x^T)\|^2
\leq 
   \frac{9}{2}\Big(1+\frac{1}{\sqrt{2}}\Big)^2\frac{(\varphi_{(1/2S_F)}(x^1)-\min_{x\in\mathcal{X} }F(x))+(S_F+\lambda) (L_f^2 L_\phi^2+\lambda ^2D_\mathcal{X}^2)}{\sqrt{T}}. 
$$
Denote 
\begin{equation}
\label{eq:constant_M}
M:=\frac{9}{2}\Big(1+\frac{1}{\sqrt{2}}\Big)^2[(\varphi_{(1/2S_F)}(x^1)-\min_{x\in\mathcal{X} }F(x))+(S_F+\lambda) (L_f^2 L_\phi^2+\lambda ^2D_\mathcal{X}^2)].
\end{equation}
The gradient mapping of $F$ satisfies the following inequality.
\begin{equation}
\label{eq:F_gradient_and_generalized_gradient}
\begin{aligned}
    \EE\|\tilde \nabla F_\alpha( \hat x^T )\|^2 
\leq & 
    2\EE\|\tilde \nabla F^\lambda_\alpha( \hat x^T )\|^2 + 2\EE\|\tilde \nabla_\alpha F( \hat x^T )- \tilde \nabla F^\lambda_\alpha( \hat x^T )\|^2\\    
= & 
    2\EE\|\tilde \nabla F^\lambda_\alpha( \hat x^T )\|^2 + 2\EE\Big\|\frac{\Pi_{\mathcal{X}}(\hat x^T-\alpha\nabla F(\hat x^T))-\Pi_{\mathcal{X}}(\hat x^T-\alpha\nabla F(\hat x^T) - \alpha\lambda \hat x^T)}{\alpha}\Big\|^2\\
\leq & 
    2\EE\|\tilde \nabla F^\lambda_\alpha( \hat x^T )\|^2 + 2\EE\|\lambda \hat x^T\|^2\\
\leq & 
    2M T^{-1/2} + 2\lambda^2 D_\mathcal{X}^2,
\end{aligned}
\end{equation}
where the first inequality uses the Cauchy-Schwarz inequality, the second inequality uses the non-expansiveness property of the projection operator~\citep{bertsekas2009convex}, i.e., $\|\Pi_\mathcal{X}(x) - \Pi_\mathcal{X}(x^\prime)\|\leq \|x-x^\prime\|$ for a convex closed set $\mathcal{X}\subseteq\RR^d$ and any $x,x^\prime\in\RR^d$, and the third inequality utilizes the fact that $\mathcal{X}$ is compact with radius $\mathcal{D}_\mathcal{X}$. In what follows, we establish the relationship between optimality gap and gradient mapping convergence.

For $u=g(x)$ and $x= g^{-1}(u)$, recall that $\tilde x= \Pi_{\mathcal{X}}(x-\alpha\nabla F(x))$. The following inequality holds
\begin{equation}
\label{eq:RSG_key}
\begin{aligned}
    F(x)-F(x^*) 
\overset{(a)}{=} &
    G(u) -G(u^*) 
\overset{(b)}{\leq} 
    ({u}-u^*)^\top \tilde \nabla_\alpha G(u) + \alpha\nabla G(u)^\top \tilde \nabla_\alpha G(u)  \\
\overset{(c)}{\leq} &
    \|u-u^*\|~\|\tilde \nabla_\alpha G({u})\|+ \alpha\|\nabla G(u)\| \|\tilde \nabla_\alpha G(u)\|\\
\overset{(d)}{=} & 
    (\|g(x)-g(x^*)\|+\alpha\|\nabla G(u)\|)~\Big\|\frac{u-\Pi_\mathcal{U}(u-\alpha\nabla G(u))}{\alpha}\Big\|\\
\overset{(e)}{\leq} & 
    (2L_\phi D_\mathcal{X}+\alpha L_fL_\phi \mu_g^{-1}) \Big\|\frac{u-\Pi_\mathcal{U}(u-\alpha\nabla G(u))}{\alpha}\Big\|\\
\overset{(f)}{=} & 
    (2L_\phi D_\mathcal{X}+ L_fL_\phi \mu_g^{-1}/2S_F)\sqrt{\sum_{i=1}^d\Big(\frac{u_i-\Pi_{\mathcal{U}_i}(u_i-\alpha[\nabla G(u)]_i)}{\alpha}\Big)^2},
\end{aligned}
\end{equation}
where (a) holds as $G(u)=F(g^{-1}(u))$, $x=g^{-1}(u)$, and $x^*=g^{-1}(u^*)$; (b) holds according to Lemma \ref{lm:gradient_generalized_gradient}; (c) holds by the Cauchy-Schwarz inequality; (d) uses the definition of gradient mapping and the fact that $g(x)=u$, $g(x^*)=u^*$; (e) uses Lipschitz continuity of  $\phi$, $f$ and $g^{-1}$ (since $\nabla g\succeq \mu_g I$) and the fact that  $\mathcal{X}$ is compact; (f) holds as $\mathcal{U}$ is a box constraint with the $i$-th coordinate interval being $\mathcal{U}_i = [\EE [\phi_i(\munderbar{X}_i,\xi_i), \phi_i(\bar{X}_i,\xi_i)]$. 

For coordinate $i\in[d]$, we divide the following analysis into two cases: 1) $x_i - \alpha [\nabla F(x)]_i\in \mathcal{X}_i = [\munderbar{X}_i, \bar X_i]$; 2) $x_i - \alpha [\nabla F(x)]_i \not \in \mathcal{X}_i$. For the first case,  we have
$$
\begin{aligned}
&
    \Big(\frac{u_i-\Pi_{\mathcal{U}_i}(u_i-\alpha[\nabla G(u)]_i)}{\alpha}\Big)^2 
\leq 
    [\nabla G(u)]_i^2 = [\nabla g(x)]_{ii}^{-2} [\nabla F(x)]_i^2\\
= &
    [\nabla g(x)]_{ii}^{-2} \Big(\frac{x_i-\Pi_{\mathcal{X}_i}(x_i-\alpha[\nabla F(x)]_i)}{\alpha}\Big)^2
\leq 
    \mu_g^{-2} \Big(\frac{x_i-\Pi_{\mathcal{X}_i}(x_i-\alpha[\nabla F(x)]_i)}{\alpha}\Big)^2,
\end{aligned}
$$
where the first inequality utilizes the fact that $u_i\in\mathcal{U}_i$, $\Pi_{\mathcal{U}_i}(u_i)=u_i$, and the  non-expansiveness of projection operator, the first equality holds as $\nabla G(u) = [\nabla g(x)]^{-\top} \nabla F(x)$ and $\nabla g(x)$ is a diagonal matrix, the second equality holds as $x_i - \alpha [\nabla F(x)]_i\in \mathcal{X}_i =[\munderbar{X}_i, \bar X_i]$, and the second inequality holds by Assumption \ref{assumption:general_2}(b). For the second case, we have 
\begin{equation}
\label{eq:gradient_mapping_FG}
\begin{aligned}
&
    \Big(\frac{u_i-\Pi_{\mathcal{U}_i}(u_i-\alpha[\nabla G(u)]_i)}{\alpha}\Big)^2
= 
    \Big(\frac{g_i(x_i)-g_i(g_i^{-1}(\Pi_{\mathcal{U}_i}(u_i-\alpha[\nabla G(u)]_i)))}{\alpha}\Big)^2 \\
\leq &
    L_\phi^2 \Big(\frac{x_i-g_i^{-1}(\Pi_{\mathcal{U}_i}(u_i-\alpha[\nabla G(u)]_i))}{\alpha}\Big)^2 
= 
    L_\phi^2 \Big(\frac{x_i-\Pi_{\mathcal{X}_i}(g_i^{-1}(u_i-\alpha[\nabla G(u)]_i))}{\alpha}\Big)^2,
\end{aligned}    
\end{equation}
where the first equality holds as $g_i$ is a bijective mapping under Assumption \ref{assumption:general}, the inequality holds as $g$ is $L_\phi$-Lipschitz continuous, and the second equality holds by Lemma \ref{lm:switch_transformation_projection}. 

If $\Pi_{\mathcal{X}_i}(x_i - \alpha [\nabla F(x)]_i) = \bar X_i$, it means $[\nabla F(x)]_i\leq 0$. Since $\nabla g$ is a diagonal positive definite matrix  and $\nabla G(u) = [\nabla g(x)]^{-1} \nabla F(x)$, it holds that $[\nabla G(u)]_i\leq 0$. As a result,  we have $u_i-\alpha[\nabla G(u)]_i\geq u_i$ and thus $g_i^{-1}(u_i-\alpha[\nabla G(u)]_i)\geq x_i$. Hence, it holds that 
$$
| x_i-\Pi_{\mathcal{X}_i}(g_i^{-1}(u_i-\alpha[\nabla G(u)]_i)) | \leq |x_i - \bar X_i| = |x_i - \Pi_{\mathcal{X}_i}(x_i - \alpha [\nabla F(x)]_i)|.
$$
A similar argument holds when $\Pi_{\mathcal{X}_i}(x_i - \alpha [\nabla F(x)]_i) = \munderbar{X}_i$. As a result, for the second case, we have
$$
     \Big(\frac{u_i-\Pi_{\mathcal{U}_i}(u_i-\alpha[\nabla G(u)]_i)}{\alpha}\Big)^2 \leq L_\phi^2 \Big(\frac{x_i - \Pi_{\mathcal{X}_i}(x_i - \alpha [\nabla F(x)]_i)}{\alpha}\Big)^2.
$$
Summarizing the two cases, we have
$$
\Big(\frac{u_i-\Pi_{\mathcal{U}_i}(u_i-\alpha[\nabla G(u)]_i)}{\alpha}\Big)^2 \leq \max\{\mu_g^{-2}, L_\phi^2\}\Big(\frac{x_i - \Pi_{\mathcal{X}_i}(x_i - \alpha [\nabla F(x)]_i)}{\alpha}\Big)^2. 
$$
Setting $x = \hat x^T$ in \eqref{eq:RSG_key} and taking full expectation, we have
$$
\begin{aligned}
    \EE [F(\hat x^T) - F(x^*)] 
\leq & 
    (2L_\phi D_\mathcal{X}+L_fL_\phi \mu_g^{-1}/2S_F)\EE \sqrt{\max\{\mu_g^{-2}, L_\phi^2\} \sum_{i=1}^d\Big(\frac{\hat x_i^T - \Pi_{\mathcal{X}_i}(\hat x_i^T- \alpha [\nabla F(\hat x^T)]_i)}{\alpha}\Big)^2}\\
= &
    (2L_\phi D_\mathcal{X}+L_fL_\phi \mu_g^{-1}/2S_F)\max\{\mu_g^{-1}, L_\phi\} \EE\Big\|\frac{\hat x^T - \Pi_{\mathcal{X}}(\hat x^T - \alpha [\nabla F(\hat x^T)])}{\alpha}\Big\|\\
= & 
    (2L_\phi D_\mathcal{X}+L_fL_\phi \mu_g^{-1}/2S_F)\max\{\mu_g^{-1}, L_\phi\} \EE\|\tilde \nabla F_\alpha(\hat x^T)\|\\
\leq &
    (2L_\phi D_\mathcal{X}+L_fL_\phi \mu_g^{-1}/2S_F) \max\{\mu_g^{-1}, L_\phi\}\sqrt{2 M T^{-1/2} + 2\lambda^2 D_\mathcal{X}^2},
\end{aligned}
$$
where the last inequality holds by \eqref{eq:F_gradient_and_generalized_gradient}.
\end{proof1}

\section{Proof of Theorem \ref{thm:msg}: Global Convergence of MSG}
\label{sec:proof_of_msg}
\begin{proof1}
Denote $u^t := g(x^t)$, $v_G(u^t):= [\nabla \hat g^{B}(x^t)]^{-\top} \nabla \phi(x^t,\xi^{t})^\top \nabla f(\phi(x^t,\xi^{t}))$. We have
\begin{equation}
\label{eq:gradient_F_and_G}
v_F(x^t) = [\nabla \hat g^{A}(x^t)]^{-\top} v_G(u^t)+\lambda x^t.
\end{equation}
We first establish an upper bound on the objective value $F(x^t)$ to the optimal objective value $F(x^*)$. For this purpose, first note that
\begin{align*}
& 
    \EE[\| u^{t+1} -u^*\|^2 \mid u^t] -\|u^t -u^*\|^2 \\
\overset{(a)}{=} ~ &
    \EE[\| g(\Pi_\mathcal{X}(x^t- \gamma v_F(x^t))) -u^*\|^2 \mid u^t]   -\|u^t -u^*\|^2 \\
\overset{(b)}{=} ~ &
    \EE[\| \Pi_\mathcal{U}(g(x^t- \gamma v_F(x^t))) -u^*\|^2 \mid u^t]  - \|u^t -u^*\|^2 \\
\overset{(c)}{\leq} ~ &
    \EE[\| g(x^t- \gamma v_F(x^t)) -u^*\|^2 \mid u^t]  - \|u^t -u^*\|^2 \\
\overset{}{=} ~ &
    \EE[\| g(x^t- \gamma v_F(x^t)) - u^t + u^t -u^*\|^2 \mid u^t]  - \|u^t -u^*\|^2 \\
\overset{(d)}{=} ~ &
    \EE[\| g(x^t- \gamma v_F(x^t)) - u^t\|^2 \mid u^t] +2\EE[ (u^t -u^*)^\top (g(x^t- \gamma v_F(x^t)) - u^t ) \mid u^t]  \\
\overset{(e)}{=} ~ &
    \EE[\| g(x^t- \gamma v_F(x^t)) - g(x^t)\|^2 \mid u^t] -2\gamma\EE[ (u^t -u^*)^\top \nabla G(u^t)\mid u^t] \\
    & + 2\EE[ (u^t -u^*)^\top [ g(x^t- \gamma v_F(x^t)) - g(x^t)  + \gamma\nabla G(u^t) \mid u^t] \\
\overset{(f)}{\leq} ~ &
    \EE[\| g(x^t- \gamma v_F(x^t)) - g(x^t)\|^2 \mid u^t] -2\gamma (G(u^t)-G(u^*)) \\
    & +2(u^t -u^*)^\top \EE[ g(x^t- \gamma v_F(x^t)) - g(x^t) + \gamma\nabla G(u^t) \mid u^t],
\end{align*}
where (a) uses the fact that $u^{t+1} = g(x^{t+1})$ and the definition of $x^{t+1}$ specified by MSG, (b) follows from Lemma \ref{lm:switch_transformation_projection}, which is the key step for handling the constraints, (c) utilizes the fact that projection operator is non-expansive and $u^* = g(x^*) =\Pi_\mathcal{U}(u^*)\in \mathcal{U}$,  (d) holds by expanding $\| g(x^t- \gamma v_F(x^t)) - u^t + u^t -u^*\|^2$, (e) follows from the definition $u^t = g(x^t)$, and (f) follows from Theorem \ref{thm:equivalent-trans-extension} that $G(u)$ is convex.
After rearranging terms and taking full expectation, we have
\begin{equation}
\label{eq:msg_per_iteration}
\begin{aligned}
    2\gamma( \EE ( G(u^t) -G(u^*)) 
\leq &
    \underbrace{\EE \|u^t -u^*\|^2 - \EE\| u^{t+1} -u^*\|^2}_{:=A_t}
    +
    \underbrace{\EE[\| g(x^t- \gamma v_F(x^t)) - g(x^t)\|^2]}_{:=B_t}\\
    & + 
    \underbrace{2\EE (u^t -u^*)^\top [ g(x^t- \gamma v_F(x^t)) - g(x^t) + \gamma\nabla G(u^t)]}_{:= C_{t}}.\\
\end{aligned}
\end{equation}
Summing up \eqref{eq:msg_per_iteration} from $t=1$ to $t=T$ and dividing $2\gamma$ on both sides, we have
\begin{equation}
\label{eq:results}
    \EE ( F(\hat x^T) -F(x^*))
=
    \frac{1}{T}\sum_{t=1}^T \EE ( F(x^t) -F(x^*)) 
=
    \frac{1}{T}\sum_{t=1}^T \EE ( G(u^t) -G(u^*)) 
\leq
    \frac{1}{T}\sum_{t=1}^T \frac{A_t+B_t+C_t}{2\gamma},
\end{equation}
where the first equality holds as $\hat x^T$ is selected uniformly from $\{x^t\}_{t=1}^T$, and the second equality holds as $F(x^t) = G(u^t)$ and $F(x^*)= G(u^*)$. It remains to upper bound the right-hand-side of \eqref{eq:results}.
Note that sum of $\{A_t\}_{t=1}^T$ forms a telescoping sum that is widely used in derivation of gradient-based methods~\citep{nemirovski2009robust}.
\begin{equation}
\label{eq:error_A}
\frac{1}{T}\sum_{t=1}^T A_t = \|u^1-u^*\|^2 -\EE \|u^{T+1}-u^*\|^2 \leq  \|u^1-u^*\|^2.
\end{equation}
Next we establish an upper bound on $B_t$. 
By Assumption \ref{assumption:general} that $\phi(x,\xi)$ is $L_\phi$-Lipschitz continuous in $x$ for any $\xi$, $g(x) = \EE \phi(x,\xi)$  is $L_\phi$-Lipschitz continuous. As a result, we have
\begin{equation}
\label{eq:error_B}
\begin{aligned}
    B_t 
= &
    \EE[\| g(x^t- \gamma v_F(x^t)) - g(x^t)\|^2]
\leq  
    L_\phi^2 \gamma^2 \EE[\|v_F(x^t)\|^2]
\leq 
    L_\phi^2 \gamma^2 \Big(\frac{K^4 L_f^2}{8L_\phi^2}+ 2\lambda^2 D_\mathcal{X}^2\Big),
\end{aligned}
\end{equation}
where the second inequality holds  by Lemma \ref{lm:bounded_second_moment} about the second moment of $v_F$. 

Upper bounding $C_t$ is another key step of the analysis. By definition, we have
\begin{align}
    C_t 
= &
    2\EE(u^t -u^*)^\top [ g(x^t- \gamma v_F(x^t)) - g(x^t) + \gamma\nabla G(u^t)] \label{eq:error_C}\\
= & 
    2\EE (u^t -u^*)^\top [ g(x^t- \gamma v_F(x^t)) - g(x^t) + \gamma v_G(u^t)] + 2\EE (u^t -u^*)^\top [ \gamma \nabla G(u^t)-\gamma v_G(u^t)]\nonumber\\ 
\leq &
    2\EE\{ \|u^t -u^*\| \| \EE [g(x^t- \gamma v_F(x^t)) - g(x^t) + \gamma v_G(u^t)\mid u^t]\|\} + 2\gamma \EE\{ \|u^t -u^*\| \|\EE [ \nabla G(u^t)-v_G(u^t)\mid u^t]\|\}\nonumber\\
\leq &
    4 L_\phi D_\mathcal{X}\EE\{\underbrace{\| \EE [ g(x^t- \gamma v_F(x^t)) - g(x^t) + \gamma v_G(u^t)\mid u^t]\|}_{:=C_{t,1}}\} + 4\gamma L_\phi D_\mathcal{X}\EE\{\underbrace{\|\EE [ \nabla G(u^t)-v_G(u^t)\mid u^t]\|}_{:=C_{t,2}}\},\nonumber
\end{align}
where the first inequality holds by the tower property and the Cauchy-Schwarz inequality, the second inequality holds as $\|u^t-u^*\| = \| g(x^t)-g(x^*)\|\leq L_\phi \|x^t-x^*\|\leq 2L_\phi D_\mathcal{X}$. Note that $g(x^t - \gamma v_F(x^t))- g(x^t)$ can be interpreted as  a ``gradient estimator'' in $u$ space which corresponds to the gradient estimator $v_F$ in $x$ space, $C_{t,1}$ reflects the approximation error between $g(x^t - \gamma v_F(x^t))- g(x^t)$ and the gradient estimator $v_G(u^t)$, and $C_{t,2}$ controls the bias of $v_G(u^t)$.
It remains to upper bound $C_{t,1}$ and $C_{t,2}$. 
Since $u^t=g(x^t)$, it holds that
\begin{align}
    C_{t,2} =& \|\EE [ \nabla G(u^t)-v_G(u^t)\mid x^t]\|
=  
    \|\EE [[\nabla g(x^t)]^{-\top} \nabla F(x^t) -[\nabla \hat g^{B}(x^t)]^{-\top} \nabla F(x^t)\mid x^t \|\nonumber\\
= &
    \|\EE [([\nabla g(x^t)]^{-1} -[\nabla \hat g^{B}(x^t)]^{-1} )^\top \nabla F(x^t)\mid x^t \|
\leq 
    \|\EE [\nabla g(x^t)]^{-1} -[\nabla \hat g^{B}(x^t)]^{-1} \mid x^t ] \|~\| \nabla F(x^t)\|\nonumber\\
\leq &
    \frac{1}{\mu_g}\Big(1- \frac{\mu_g}{2L_\phi}\Big)^K L_\phi L_f, \label{eq:error_C_2}
\end{align}
where the first inequality uses  the Cauchy-Schwarz inequality, and the second inequality holds by Lemma \ref{lm:properties_of_inverse_expectation_estimation} and the fact that $\|\nabla F(x)\|\leq L_\phi L_f$.
Next we bound $C_{t,1}$. 
\begin{align}
&
    C_{t,1} = \Big\|\EE\Big[ g(x^t - \gamma v_F(x^t)) - u^t + \gamma v_G(u^t)\mid u^t\Big]\Big\|\nonumber\\
= & 
    \Big\|\EE\Big[ g(x^t - \gamma v_F(x^t)) - g(x^t) +\gamma \nabla g(x^t)^\top v_F(x^t)  -\gamma \nabla g(x^t)^\top v_F(x^t) + \gamma v_G(u^t)\mid u^t\Big]\Big\|\nonumber\\
\leq & 
    \Big\|\EE\Big[ g(x^t - \gamma v_F(x^t)) - g(x^t) +\gamma \nabla g(x^t)^\top v_F(x^t)\mid u^t\Big]\Big\| + \Big\| \EE \Big[\gamma v_G(u^t) -\gamma \nabla g(x^t)^\top v_F(x^t) \mid u^t\Big]\Big\|\nonumber\\
\leq &
    \frac{\gamma^2 S_g}{2}\EE[   \|v_F(x^t)\|^2  \mid u^t]+\gamma\Big\|\EE \Big[v_G(u^t)-  \nabla g(x^t)^\top v_F(x^t)\mid u^t\Big] \Big\|, \label{eq:error_C_1}
\end{align}
where the equality uses the fact that $u^t =g(x^t)$, the first inequality uses the triangle inequality, and the second inequality uses the Lipschitz continuity of $\nabla g$, i.e., $g$ is smooth.

For the first term, with Lemma \ref{lm:bounded_second_moment}, we have
\begin{equation}
\label{eq:second-order-bounds}
    \frac{\gamma^2 S_g}{2}\EE[\| v_F(x^t)\|^2\mid u^t] 
\leq
    \frac{\gamma^2 S_g}{2}\Big(\frac{K^4 L_f^2}{8L_\phi^2}+ 2\lambda^2 D_\mathcal{X}^2\Big).
\end{equation}

For the second term,  by \eqref{eq:gradient_F_and_G}, we have $v_F(x^t) = [\nabla \hat g^{A}(x^t)]^{-1}v_G(u^t)+\lambda x^t$. It holds that
\begin{align*}
&    
    \gamma\Big\|\EE \Big[v_G(u^t)-  \nabla g(x^t)^\top v_F(x^t)\mid u^t\Big] \Big\|\\
= & 
    \gamma\Big\|\EE \Big[v_G(u^t)-   \nabla g(x^t)^\top[\nabla \hat g^{A}(x^t)]^{-\top} v_G(u^t)-\lambda \nabla g(x^t)^\top x^t\mid u^t\Big] \Big\|\\
\leq &  
    \gamma\Big\|\nabla g(x^t)^\top \EE \Big[ (\nabla g(x^t)^{-\top} -   [\nabla \hat g^{A}(x^t)]^{-\top}) v_G(u^t)\mid u^t\Big] \Big\| + \gamma \|\lambda \nabla g(x^t)^\top x^t\|\\
\leq & 
    \gamma\Big\|\nabla g(x^t)^\top  \Big[\EE_{k_1, \{\xi^{ti}\}_{i=1}^{k_1}}(\nabla g(x^t)^{-\top} -   [\nabla \hat g^{A}(x^t)]^{-\top})\EE_{\xi^t, \{\xi^{tj}\}_{j=1}^{k_2}, k_2} v_G(u^t) \Big] \Big\| + \gamma \lambda\| \nabla g(x^t)\|\|x^t\|\\
\leq &
    \gamma\|\nabla g(x^t)\|\Big[  \Big\|\EE_{k_1, \{\xi^{ti}\}_{i=1}^{k_1}}(\nabla g(x^t)^{-1} -   [\nabla \hat g^{A}(x^t)]^{-1})\Big\| ~\EE_{\xi^t, \{\xi^{tj}\}_{j=1}^{k_2}, k_2} \|v_G(u^t)\|\Big] + \gamma \lambda\| \nabla g(x^t)\|\|x^t\|\\
\leq &    
    \gamma L_\phi  \Big\|\EE_{k_1, \{\xi^{ti}\}_{i=1}^{k_1}}([\nabla g(x^t)]^{-1} -   [\nabla \hat g^{A}(x^t)]^{-1})\Big\| \EE_{\xi^t, \{\xi^{tj}\}_{j=1}^{k_2}, k_2} \|v_G(u^t)\|    + \gamma L_\phi\lambda D_\mathcal{X},
\end{align*}    
where the first inequality uses the triangle inequality, the second inequality uses the tower property for conditional expectation where we specify each expectation with respect to what randomness, and the last inequality holds by the Cauchy-Schwarz inequality and the fact that $\|\nabla g(x^t)\|\leq L_\phi$. Using Lemma \ref{lm:properties_of_inverse_expectation_estimation} about the bias of matrix inverse estimator and the first moment of $v_G$ derived via Jensen's inequality from Lemma \ref{lm:bounded_second_moment} about second moment of $v_G$, we have
\begin{equation}
\label{eq:first-order-bounds}
    \gamma\Big\|\EE \Big[v_G(u^t)-  \nabla g(x^t) v_F(x^t)\mid u^t\Big] \Big\|
\leq 
    \gamma \frac{L_\phi K L_f}{2\mu_g}\Big(1- \frac{\mu_g}{2L_\phi}\Big)^K + \gamma L_\phi\lambda D_\mathcal{X}.
\end{equation}
Plugging \eqref{eq:second-order-bounds}, \eqref{eq:first-order-bounds} into \eqref{eq:error_C_1}, we have
$$
C_{t,1}\leq \gamma \frac{L_\phi K L_f}{2\mu_g}\Big(1- \frac{\mu_g}{2L_\phi}\Big)^K + \gamma L_\phi\lambda D_\mathcal{X} + \frac{\gamma^2 S_g}{2}\Big(\frac{K^4 L_f^2}{8L_\phi^2}+ 2\lambda^2 D_\mathcal{X}^2\Big).
$$
Combining with \eqref{eq:error_C_2} and \eqref{eq:error_C}, we have  
$$
C_t \leq 2 L_\phi^2 D_\mathcal{X} \gamma \frac{K L_f+2L_f}{\mu_g}\Big(1- \frac{\mu_g}{2L_\phi}\Big)^K + 4 L_\phi^2 D_\mathcal{X}^2 \gamma \lambda  + 2L_\phi D_\mathcal{X}\gamma^2 S_g \Big(\frac{K^4 L_f^2}{8L_\phi^2}+ 2\lambda^2 D_\mathcal{X}^2\Big).
$$
Together with \eqref{eq:error_B}, \eqref{eq:error_A}, and \eqref{eq:results},  we have
\begin{align*}
    \EE[F(\hat x^T)-F(x^*)]
\leq &
    \frac{\|u^1 -u^*\|^2}{2\gamma T}+  (L_\phi^2 \gamma +2L_\phi D_\mathcal{X}\gamma S_g )\Big(\frac{K^4 L_f^2}{16L_\phi^2}+ \lambda^2 D_\mathcal{X}^2\Big) + 2 L_\phi^2 D_\mathcal{X}^2 \lambda \\
    & +
    L_\phi^2 D_\mathcal{X} \frac{K L_f+2L_f}{\mu_g}\Big(1- \frac{\mu_g}{2L_\phi}\Big)^K.
\end{align*}
Plugging in $\gamma =c_1T^{-1/2}$ and $\lambda \in [0, c_2T^{-1/2}]$ obtains the desired result. 
\end{proof1}

\section{Discussions on Assumptions}
\label{appendix:discussion_on_assupmtions}
In this section, we discuss the conditions required for operations management applications to ensure the assumptions needed by the global convergence results of RSG and MSG. We also discuss situations when certain assumptions do not hold and what happens to the practical performance of the proposed algorithm. 

The following Table \ref{tab:algos_assumption} summarizes the assumptions needed for global convergence of RSG, MSG, and SAA+SG algorithms. 

{
\begin{table}[t]\centering \footnotesize
    \caption{Summary of Assumptions and Complexity of Global Convergence} 
    \renewcommand\arraystretch{1.2}
    \label{tab:algos_assumption}
    \begin{tabular}{@{\extracolsep{1pt}}cccc}
    \\\hline
    \hline
    \multirow{2}{*}{Algorithm}  & \multirow{2}{*}{Assumptions} &{Sample} & {Gradient} \\ 
    {} & {} & {Complexity} & {Complexity}\\
    \hline\\[-1.8ex]

    SAA+SG   &  Assumptions \ref{assumption:reformulation} and \ref{assumption:general} & \multirow{2}{*}{$\tilde\cO(d\eps^{-2})$ }   & \multirow{2}{*}{$\tilde\cO(d^2\eps^{-4})$}         \\
    (Algorithm \ref{alg:sgd_on_g} Theorem \ref{lemma:saa+sgd})  & $f(\phi(x;\xi))$ is sub-Gaussian & &  \\ 
    \\[-1.8ex]\hline\\[-1.8ex]
    
    RSG   &  Assumptions \ref{assumption:reformulation} and \ref{assumption:general} & \multirow{2}{*}{$\cO(\eps^{-4})$ }   & \multirow{2}{*}{$\cO(\eps^{-4})$}         \\
    (Algorithm \ref{alg:RSG} Theorem \ref{thm:RSG})  & Assumption \ref{assumption:general_2} & &  \\ 
    \\[-1.8ex]\hline\\[-1.8ex]
    
    MSG   &  Assumptions \ref{assumption:reformulation} and \ref{assumption:general} & \multirow{2}{*}{$\tilde \cO(\eps^{-2})$}    & \multirow{2}{*}{$\tilde \cO(\eps^{-2})$}        \\
    (Algorithm \ref{alg:MSG} Theorem \ref{thm:msg})  & Assumption \ref{assumption:general_2} & &  \\
    \\[-1.8ex]\hline\\[-1.8ex]

    MSG   &  Assumptions \ref{assumption:reformulation} and \ref{assumption:general} & \multirow{2}{*}{$\tilde \cO(\eps^{-2})$}    & \multirow{2}{*}{$\tilde \cO(\eps^{-2})$}        \\
    (Algorithm \ref{alg:MSG} Theorem \ref{cor:msg_without_smoothness})  & $\xi\sim\mathrm{Poisson}(\beta)$ with a large $\beta$  & &  \\
    \\[-1.8ex]\hline\\[-1.8ex]

     MSG   &  Assumptions \ref{assumption:reformulation} and \ref{assumption:general} & \multirow{2}{*}{$\tilde \cO(\eps^{-2})$}    & \multirow{2}{*}{$\tilde \cO(\eps^{-2})$}        \\
    (Algorithm \ref{alg:MSG} Theorem \ref{cor:msg_without_smoothness})  & $\xi\sim\mathrm{Multinomial}(n)$ with a large $n$& &  \\
    \\[-1.8ex]\hline
    \hline\\[-1.8ex]
    \end{tabular}
    \end{table}
}

\subsection{Conditions of $\phi$ and $\PP(\xi)$ to Ensure Assumption \ref{assumption:general_2}}
\label{appendix:verify_all_phi}
One could question that Assumption \ref{assumption:general_2} might be hard to satisfy for some $\phi$ function and distribution $\PP(\xi)$ that appears in applications.
Below we list two sets of combinations of conditions on $\phi$ and $\PP(\xi)$ to ensure Assumption \ref{assumption:general_2}.

\begin{lemma}
\label{lemma:conditions_for_smoothness}
To ensure that Assumption \ref{assumption:general_2} holds, it suffices to have either one of the two conditions:
\begin{itemize}
    \item [(i.)] For any $x\in\mathcal{X}$, function $\nabla \phi(x,\xi)$ is $S_g$-Lipschitz continuous and  $\phi(x,\xi)$ satisfies  $\mu_g I\preceq \nabla \phi(x,\xi)$ for any realization of $\xi$ within the support of $\PP(\xi)$.
    \item [(ii.)] For any $x\in\mathcal{X}$, function   $\phi(x,\xi)$ satisfies  $\mu_\phi(x,\xi) I\preceq \nabla \phi(x,\xi)$ for some $\mu_\phi(x,\xi)\geq0$ and any realization of $\xi$ within the support of $\PP(\xi)$. In addition, it holds that $\EE \mu_\phi(x,\xi)\geq\mu_g>0$ for any $x\in\mathcal{X}$. Function $\nabla \phi(x,\xi)$ is not Lipschitz continuous yet $\EE \nabla \phi(x,\xi)$ is $S_g$-Lipschitz continuous.
\end{itemize}
\end{lemma}

The proof of the lemma is obvious and thus omitted. Next, we show that the four $\phi$ functions listed in Section \ref{sec:alg}, i.e., $\phi(x,\xi)=x\xi$, $\phi(x,\xi)=x\xi/(x+\alpha\xi^\kappa)$, $\phi(x,\xi)=(x/(x+\xi))k$, and $\phi(x,\xi) = x\wedge\xi$, all satisfy one of the conditions listed in the lemma above. 

\begin{lemma}
\label{lemma:verify_all_phi}
We have the following results.
\begin{itemize}
    \item For $\phi(x,\xi)=x\xi$, $\phi(x,\xi)=x\xi/(x+\alpha\xi^\kappa)$, and $\phi(x,\xi)=(x/(x+\xi))k$, suppose that the domain $\mathcal{X}\subseteq \RR^d_+$ is nonnegative and compact, and the support of the distribution $\PP(\xi)$ is nonnegative and bounded, then the first condition in Lemma \ref{lemma:conditions_for_smoothness} holds.
    \item For $\phi(x,\xi) = x \wedge\xi$, when $\PP(\xi_i\geq\bar X_i)\geq \mu_g$ for all $i=1,...,d$ and the CDF of the distribution $\PP(\xi)$ is $S_g$-Lipschitz continuous, then the second condition in Lemma \ref{lemma:conditions_for_smoothness} holds.
\end{itemize}
\end{lemma}
The following Table \ref{tab:phi_assumption} summarize the conditions on $\phi$ and $\PP(\xi)$ to ensure that Assumption \ref{assumption:general_2} holds. The assumption needed for $\phi(x,\xi)=x\wedge\xi$ implies that $\xi$ is a continuously distributed random vector and that $\bar X_i<\mathrm{ess}\sup \xi_i$ for all $i=1,...,d$. We shall discuss in the next subsection what if these conditions do not hold so that Assumption \ref{assumption:general_2} fails. Next, we show the proof of Lemma \ref{lemma:verify_all_phi}.
{
\begin{table}[t]\centering \scriptsize
    \caption{Conditions on $\phi$ and $\PP(\xi)$ to Ensure Assumption \ref{assumption:general_2}} 
    \renewcommand\arraystretch{1.2}
    \label{tab:phi_assumption}
    \begin{tabular}{@{\extracolsep{1pt}}cccc}
    \\\hline
    \hline
    \multirow{2}{*}{Function $\phi$ }  & \multirow{2}{*}{Conditions Needed} &\multirow{1}{*}{Assumption \ref{assumption:general_2}(a)} & \multirow{1}{*}{Assumption \ref{assumption:general_2}(b)} \\ 
    {} & {} & $\nabla g(x)\succeq \mu_g I$  & $\nabla g$ is Lipschitz continuous \\
    \hline\\[-1.8ex]

    {$\phi(x,\xi)=x\xi$}   &  \multirow{2}{*}{Nonnegative and compact domain $\mathcal{X}\subset \RR^d_+$}  & \multirow{4}{*}{\checkmark }   & \multirow{4}{*}{\checkmark }        \\
    \multirow{2}{*}{$\phi(x,\xi)=x\xi/(x+\alpha\xi^\kappa)$}  &  & &  \\ 
      & \multirow{2}{*}{Nonnegative and bounded support of $\PP(\xi)$} & &  \\ 
     $\phi(x,\xi)=(x/(x+\xi))k$ & & &  \\ 
    \\[-1.8ex]\hline\\[-1.8ex]
    
    \multirow{4}{*}{$\phi(x,\xi)=x\wedge\xi$}   &  $\PP(\xi_i\geq\bar X_i)\geq \mu_g$ for all $i\in[d]$ & \multirow{2}{*}{\checkmark}    &        \\
    &(implying that $\bar X_i<\mathrm{ess}\sup \xi_i$) && \\[4pt]
    \cline{2-4}
     & CDF of $\xi$ is $S_g$-Lipschitz continuous for all $i\in[d]$ & &\multirow{2}{*}{\checkmark}   \\
     &(implying that $\xi$ is continuous distributed) && \\
    \\[-1.8ex]\hline
    \hline\\[-1.8ex]
    \end{tabular}

    \end{table}
} 
\begin{proof1}
For ease of demonstration, we consider the case when $d=1$. It can be easily generalized to higher dimensions as $\phi(x,\xi) = (\phi_1(x_1,\xi_1),\ldots,\phi_d(x_d,\xi_d))^\top$ is separable. 

We first show that $\phi(x,\xi)=x\xi$, $\phi(x,\xi)=x\xi/(x+\alpha\xi^\kappa)$, and $\phi(x,\xi)=(x/(x+\xi))k$ satisfy the first condition in Lemma \ref{lemma:conditions_for_smoothness}. By the asssumption on the nonnegative domain $\mathcal{X}$ and the nonnegative support of the distribution $\PP(\xi)$, without loss of generality, we assume that $0< \munderbar{\xi} \leq \xi\leq \bar \xi$ and $0<\munderbar{X}\leq x\leq \bar{X}$. 
It is easy to see that these three $\phi$ functions are continuously differentiable, Lipschitz continuous, and strictly increasing.
\begin{itemize}
    \item In example $\phi(x,\xi)=x\xi$, we have $\nabla \phi(x,\xi)=\xi$. Thus $\mu_g=\munderbar{\xi}$. In addition, $\nabla^2 \phi(x,\xi) = 0$, thus $S_g$ can be any positive number. Thus $\phi(x,\xi)=x\xi$ satisfies the first condition in Lemma \ref{lemma:conditions_for_smoothness}.
    \item In example $\phi(x,\xi)=x\xi/(x+\alpha\xi^\kappa)$, we have $\nabla \phi(x,\xi)=(\alpha\xi^{\kappa+1})/{(x+\alpha\xi^\kappa)^2}$, which is monotonically decreasing in $x\in[0,\infty)$. As a result, together with the boundedness of $\xi$ and $x$, one can easily verify that $\mu_g=\alpha\munderbar{\xi}^{\kappa+1}/(\bar{X}+\alpha^2\bar{\xi}^{2\kappa})^2$. On the other hand, $\nabla^2\phi(x,\xi) = -2 (\alpha\xi^{\kappa+1})/{(x+\alpha\xi^\kappa)^3}$, which is monotonically increasing in $x\in[\munderbar{X},\bar X]$. Thus $S_g = -2 (\alpha\munderbar{\xi}^{\kappa+1})/{(\bar X+\alpha\bar \xi^\kappa)^3}$. As a result, $\phi(x,\xi)=x\xi/(x+\alpha\xi^\kappa)$  satisfies the first condition in Lemma \ref{lemma:conditions_for_smoothness}.
    \item In example $\phi(x,\xi)=(x/(x+\xi))k$, we have $\nabla \phi(x,\xi) = (\xi k )/(x+\xi)^2$, which is monotonically decreasing in $x\in[\munderbar{X},\bar X]$. Together with the boundedness of $\xi$ and $x$, we have $\mu_g=k\munderbar{\xi}/(\bar{X}+\bar{\xi})^2$. In addition, $\nabla^2 \phi(x,\xi) = -2(\xi k )/(x+\xi)^3$, which is monotonically increasing in $x$. Thus $S_g = -2(\munderbar{\xi} k )/(\bar x+\bar \xi)^3$.  As a result, $\phi(x,\xi)=(x/(x+\xi))k$  satisfies the first condition in Lemma \ref{lemma:conditions_for_smoothness}.
\end{itemize}  

Next, we show that $\phi(x,\xi)=x\wedge\xi$ satisfies the second condition in Lemma \ref{lemma:conditions_for_smoothness} under the specified assumptions. For any given $x\in\mathcal{X}$, notice that $\nabla \phi(x,\xi) = \mathbb{I}(x\leq \xi)$.  Thus $\mu_\phi(x,\xi) = 0$ if $x>\xi$ and $\mu_\phi(x,\xi)=1$ if $x\leq \xi$. Thus
$\EE \mu_\phi(x,\xi)= \PP(\xi\geq x) \geq \PP(\xi\geq \bar X)\geq \mu_g$. Note that $\nabla g\succeq \mu_g$ can also be shown via
$$
    \nabla g(x) = 1-H(x)=\PP(\xi\geq x)\geq\PP(\xi\geq \bar X)\geq \mu_g,
    $$
where $H$ is the CDF of $\xi$ and the inequality holds by the specified assumption. On the other hand, function $\phi(x,\xi)=\mathbb{I}(x\leq \xi)$ is not Lipschitz continuous. However, we have for any $x,y\in\mathcal{X}$ that
    $$
    |\nabla g(x) -\nabla g(y)| = |H(y) - H(x)|\leq S_g |x - y|,
    $$
where the last inequality holds as the CDF of the distribution $\PP(\xi)$ is $S_g$-Lipschitz continuous.
\end{proof1}

{\color{black}
\subsection{Proof of Lemma \ref{lemma:verify_x_wedge_xi}}
\label{proof_of_verify_x_wedge_xi}
\begin{proof1}
For $\phi(x,\xi)= x\wedge \xi$, when  $\xi$ is component-wise independent random vector and the CDF of $\xi_i$, $H_i(x_i)$, is $S_g$-Lipschitz continuous and $1-H(\bar X_i) = \PP(\xi_i\geq \bar X_i)\geq \mu_g$ for any $i\in [d]$, we verify Assumption \ref{assumption:reformulation}(b)(c), Assumption \ref{assumption:general}(c)(d) and Assumption \ref{assumption:general_2} and that $F$ is smooth. 
Note that Assumption \ref{assumption:general_2} is also verified in Lemma \ref{lemma:verify_all_phi}.
\begin{enumerate}

    \item[\textit{Verification of Assumption \ref{assumption:reformulation}(b).}] 
    It is obvious that $x\wedge\xi$ is component-wise non-decreasing in $x$ for any given $\xi$.
    
    \item[\textit{Verification of Assumption \ref{assumption:reformulation}(c).}] 
    \citet{feng2018supply} has shown that $x\wedge\xi$ is stochastic linear in mid-point when $d=1$. The extension to high-dimensional cases follows as $\xi$ is component-wise independent and $x\wedge \xi = (x_1\wedge\xi_1,\ldots,x_d\wedge\xi_d)$.
    
    \item[\textit{Verification of Assumption \ref{assumption:general}(c).}] 
    $$
    \|\phi(x,\xi) - \phi(y,\xi)\| = \|x\wedge\xi - y\wedge\xi\|\leq \|x - y\|.
    $$
    Thus $\phi(x,\xi)$ is $1$-Lipschitz continuous in $x$ for any given $\xi$. 
    \item[\textit{Verification of Assumption \ref{assumption:general}(d).}] 
    
    Since $\phi(x,\xi)=x\wedge \xi$, the only non-differetiable points are within $\{x\mid x_i = \xi_i \text{ for some } i\in[d] \}$, which forms a zero-measure set. Thus $\phi(x,\xi)$ is almost everywhere differentiable in $x\in\mathcal{X}$ for any given $\xi$.

    To show that $\phi(x,\xi)$ is almost surely differentiable for $x\in\mathcal{X}$,  equivalently, we need to show  for $x\in\mathcal{X}$ that
    $$
    \PP(\xi\mid x\wedge\xi \text{ is differentiable in } x) =1.
    $$
    It is equivalent to 
    $$
    \PP(\xi\mid \xi_i \not = x_i \text{ for any } i\in[d])=1.
    $$
    Let us assume that $\PP(\xi\mid \xi_i \not = x_i \text{ for any } i\in[d])<1$, i.e., there exists a $x^0$ such that $\PP(\xi\mid \xi_i  = x_i^0 \text{ for some } i\in[d])>0$.
    It contradicts the fact that the CDF of $\xi$ is $L_H$-Lipschitz continuous. Therefore, we obtain the desired result. 
    {\item[\textit{Verification of Assumption \ref{assumption:general_2}(a).}] 

    Since $\xi$ is component-wise independent, and $g(x) = \EE[x\wedge\xi] $, we consider the $i$-th coordinate. 
    $$
    \lim_{\Delta x\rightarrow 0} \frac{g_i(x_i+\Delta_x) - g_i(x_i)}{\Delta_x} =  1-H_i(x_i).
    $$
    Since $H$ is continuously differentiable, we know that $\nabla g$ exists and that $g$ is continuously differentiable. 
    }
    \item[\textit{Verification of Assumption \ref{assumption:general_2}(b).}] 
    
    Since $\xi$ is component-wise independent, $\nabla g(x)$ is a diagonal matrix.  
    Since $\PP(\xi_i\geq \bar X_i)\geq \tilde \mu_g$ for any $i\in [d]$,  we have 
    $$
    \nabla_i g_i(x_i) = 1-H_i(x_i)=\PP(\xi_i\geq \bar X_i)\geq \tilde \mu_g.
    $$
    Therefore $\nabla g(x)\succeq \tilde \mu_g I$.
    \item[\textit{Verification of Assumption \ref{assumption:general_2}(c).}] 
    Since $\xi$ is component-wise independent, $\nabla g(x)$ is a diagonal matrix. In addition, it holds that $\nabla_i g_i(x_i) = 1-H_i(x_i)$. Since $H_i(x_i)$ is $L_H$-Lipschitz continuous for any $i\in[d]$, we have 
    $$
    |\nabla_i g_i(x_i) -\nabla_i g_i(y_i)| = |H_i(y_i) - H_i(x_i)|\leq L_H |x_i - y_i|,
    $$
    where $x_i$ and $y_i$ are the $i$-th coordinate of $x, y\in\mathcal{X}$. As a result, $\nabla g(x)$ is $L_H$-Lipschitz continuous. 
    
\item[\textit{Verification of Lipschitz continuity of $\nabla F$}]
Without loss of generality, we consider the case when $d=1$. The extension to a higher-dimensional case is straightforward. For $\phi(x,\xi)=x\wedge\xi$, we have 
$$
\begin{aligned}
    \|\nabla F(x) -\nabla F(y) \|
= &
    \|\EE \mathbf{1}(x\leq \xi) \nabla f(x\wedge\xi) -\mathbf{1}(y\leq \xi) \nabla f(y\wedge\xi) \| \\
\leq &
    \|\EE \mathbf{1}(x\leq \xi) [\nabla f(x\wedge\xi) - \nabla f(y\wedge\xi)] \| + \|\EE [\mathbf{1}(x\leq \xi) - \mathbf{1}(y\leq \xi)]\nabla f(y\wedge\xi)\|\\
\leq &
    \EE\|\nabla f(x\wedge\xi) - \nabla f(y\wedge\xi) \|+ \Big|\int_{t\in[\min(x,y),\max(x,y)]} \nabla f(y\wedge t) d H(t)\Big|\\
\leq &
    S_f |x-y|+ L_f\int_{t\in[\min(x,y),\max(x,y)]} d H(t)\\
\leq &
    (S_f + L_H)|x-y|,
\end{aligned}
$$
where the third inequality uses smoothness of $F$ and the fourth inequality uses Lipschitz continuity of $H$. It implies that $\nabla F$ is Lipschitz continuous.
\end{enumerate}
\end{proof1}
}

\subsection{Performance of RSG and MSG when Assumption \ref{assumption:general_2}(b) Fails.}
\label{appendix:a2.3a_fails}
In the following two subsections, we discuss what happens to RSG and MSG when these requirements to ensure Assumption \ref{assumption:general_2} are not satisfied when $\phi(x,\xi)=x\wedge\xi$. Note that the SAA+SG method does not require Assumption \ref{assumption:general_2}, and thus SAA+SG is not influenced. We first consider when Assumption \ref{assumption:general_2}(a) fails.

Recall that it requires $\PP(\xi_i\geq \bar X_i)\geq\mu_g$ for all $i\in[d]$ to ensure that  $\nabla g(x)\succeq\mu_g I$ for any $x\in\mathcal{X}$. Without loss of generality, we consider the one-dimensional case, i.e., $d=1$. Suppose that $\PP(\xi\geq \bar X)\geq\mu_g$ does not hold for any $\mu_g>0$. By the analysis of Lemma \ref{lemma:verify_x_wedge_xi}, we know that there exists a $\hat x\in\mathcal{X}$ such that $\PP(\xi\geq \hat x) = 0$. In other words, $\mathrm{ess}\sup\xi<\hat x\leq \bar X$. Recall the discussion in Example \ref{example:RSG_on_truncation}. If RSG and MSG encounter such a point $\hat x$ at iterate $t$, i.e., $x^t=\hat x$, we have $\nabla [x\wedge \xi]=\mathbb{I}(\hat x\leq \xi)=0$ with probability $1$. Thus RSG and MSG algorithms with a regularization parameter $\lambda=0$ could get stuck at a local point within $\mathcal{X}_\mathrm{local}^*$ and fail to converge globally. However, since we deliberately use a non-zero regularization term in the gradient estimator design, as discussed in Example \ref{example:RSG_on_truncation}, RSG and MSG will automatically shrink $\hat x$ such that the next update is $x^{t+1}=(1-\lambda\gamma)x^t =(1-\lambda\gamma)\hat x$. Such a shrinking update ensures that we will find a $\tilde x$ such that $\tilde x<\mathrm{ess}\sup\xi$. As a result, even if $\bar X_i>\mathrm{ess}\sup \xi$ so that $\PP(\xi\geq\bar X)\geq \mu_g$ does not hold, RSG and MSG automatically avoid the trivial local solutions, i.e., large $x$ and create an ``effective" upper bound $\tilde X$  such that the algorithm is sure to converge below it. In addition, it holds that $\tilde X_i<\mathrm{ess}\sup \xi_i<\bar X_i$ and that $\PP(\xi_i\geq\tilde X_i)\geq \mu_g$ for some $\mu_g>0$. Then RSG and MSG will search for the optimal solution over the ``effective domain" $\tilde{\mathcal{X}} = [\munderbar{X},\tilde X]$. Compared to optimizing over the original domain $\mathcal{X}$, optimizing over the effective domain only rules out the local solutions, i.e., $\mathcal{X}_\mathrm{local}^*$ as defined in Example \ref{example:RSG_on_truncation}. In addition, the needed condition $\PP(\xi\geq \tilde X)\geq\mu$ holds on the effective domain. Thus
the conclusion is that without this assumption, the practical performance of RSG and MSG is not much influenced. This assumption is only needed for demonstrating rigorous analysis. 

\subsection{Performance of MSG when Assumption \ref{assumption:general_2}(c) Fails}
\label{appendix:a2.3_fails}

In what follows, we discuss the case when $g(x)=\EE [x\wedge\xi]$ is not smooth, i.e., Assumption \ref{assumption:general_2}(c) fails. Note that this assumption is only used in the analysis of MSG. Thus we investigate how MSG behaves without such an assumption. In particular, we investigate the NRM case when the distribution of $\xi$ follows a Poisson or a multinomial distribution. As the smoothness of $g$ requires the CDF of $\xi$ to be Lipschitz continuous, it is clearly not satisfied when $\xi$ is a discrete random vector. This paragraph serves as a 

Recall the global convergence analysis of MSG. The smoothness of $g$ is only used in the analysis of inequality \eqref{eq:error_C_1} to obtain
\begin{align*}
&
    \Big|\EE\Big[ g(x^t - \gamma v_F(x^t)) - g(x^t) +\gamma \nabla g(x^t)^\top v_F(x^t)\mid u^t\Big]\Big|\leq \frac{\gamma^2 S_g}{2}\EE[   \|v_F(x^t)\|^2  \mid u^t].
\end{align*}

We particularly utilize the fact that the right-hand side depends quadratically on the stepsize $\gamma$ while $\EE[\|v_F(x^t)\|^2\mid u^t]$ is treated as $\tilde\cO(1)$. When translating into the final convergence rate, this error term needs to divide $\gamma$, and thus an $\cO(\gamma)$ error appears, which is of order $\cO(T^{-1/2})$ when picking $\gamma = T^{-1/2}$.

Now we try to bound this term without smoothness of $g$. We first consider the one-dimensional setting as $d=1$. {By definition $g(x)=\EE[x\wedge\xi]$, when $\xi$ follows a discrete distribution on $\mathbb{Z}$, $g$ is differentiable and smooth on $\mathcal{X}/\mathbb{Z}$. The non-differentiable points of $g$ on $\mathcal{X}$ is of measure zero. Thus without loss of generality, we consider when $x^t$ and $x^t - \gamma v_F(x^t)$ are differentiable points of $g$. } Notice that $g(x)$ is a concave function in $x$ as $x\wedge\xi$ is concave. Therefore, we can bound the desired term utilizing the concavity:
\begin{align*}
    &\Big\|\EE\Big[ g(x^t - \gamma v_F(x^t)) - g(x^t) +\gamma \nabla g(x^t)^\top v_F(x^t)\mid u^t\Big]\Big\|\\
    =&\EE\Big[g(x^t)-\gamma\nabla g(x^t)^\top v_F(x^t) - g(x^t - \gamma v_F(x^t)) \mid u^t \Big]\\
    \leq&\EE\Big[\gamma\nabla g(x^t - \gamma v_F(x^t))^\top v_F(x^t)-\gamma\nabla g(x^t)^\top v_F(x^t) \mid u^t \Big]\\
    \leq &\gamma \EE\Big[\|\nabla g(x^t - \gamma v_F(x^t))-\nabla g(x^t)\| \|v_F(x^t)\|\mid u^t \Big],
\end{align*}
where the first equality uses the concavity of $g(x)=\EE[x\wedge\xi]$, i.e., $g(x^t - \gamma v_F(x^t)) - g(x^t)+\gamma\nabla g(x^t)^\top v_F(x^t)\leq 0$, to get rid of the norm, the second inequality uses again the concavity of $g$ such that $g(x^t)-g(x^t - \gamma v_F(x^t)) \leq \gamma\nabla g(x^t - \gamma v_F(x^t))^\top v_F(x^t)$, and the third inequality holds by Cauchy-Schwarz inequality. By the definition of $v_F(x^t)$, we know that $\|v_F(x^t)\|= \cO(K^2L_f+\lambda \bar X)$ admits an uniform upper bound. Without loss of generality, we denote such an upper bound as $M_F$, i.e., $v_F(x)\leq M_F$. Further dividing $\gamma$, this error term leads to a $\cO\Big(\frac{1}{T}\sum_{t=1}^T\EE\Big[\|\nabla g(x^t - \gamma v_F(x^t))-\nabla g(x^t)\|\Big]\Big)$ error in the final global convergence rate.

Next, we derive a upper bound on $\|\nabla g(x - \gamma v_F(x))-\nabla g(x)\|$ for any $x\in\mathcal{X}$. Note that $\xi$ has a support $\mathbb{Z}$, $\PP(\xi = k) = p_{k}$, and $k\in \mathbb{Z}$. Since $\nabla g(x) = 1-H(x)$ where $H$ denotes the CDF of $\PP(\xi)$, we know  
$$
|\nabla g(x)-\nabla g(x - \gamma v_F(x))| = |H(x - \gamma v_F(x)) - H(x)|.
$$  
Thus the difference reflects the cumulative distribution function at two points.
Pick a stepsize $\gamma = \cO(1/\sqrt{T})$.  Note that when we want to achieve a high accuracy $\eps$, the number of iterations $T$ has to be large. For large $T$, without loss of generality, we know $|\gamma v_F(x)|< 1$. Notice that $\xi$ takes values in $\mathbb{Z}$. This means that within the interval of $[x - \gamma v_F(x), x]$ or $[x, x - \gamma v_F(x)]$, there is only one integer number which $\xi$ can take value. As a result, the difference between $\nabla g(x)$ and $\nabla g(x - \gamma v_F(x))$, equivalently the difference between $H(x - \gamma v_F(x))$ and $H(x)$, is at most the probability mass of $\PP(\xi)$ at one integer point in $\mathbb{Z}$. Equivalently, we have 
$$
|\nabla g(x - \gamma v_F(x))-\nabla g(x)|\leq = |H(x - \gamma v_F(x)) - H(x)|\leq \max_{k\in\mathbb{Z}}p_k.
$$

In OM practice, the error term $\max_{k\in\mathbb{Z}}p_{k}$ can be very small for certain distribution $\PP(\xi)$.
\begin{itemize}
    \item When $\xi$ takes Poisson distribution with a large arrival rate $\beta$. Note that 
$$
\max_{k\in \mathbb{Z}}\frac{\beta^k e^{-\beta}}{k!}\approx \frac {e^{-\beta}\beta^{\beta}}{\beta!}\approx \frac{1}{\sqrt{2\pi \beta}}.
$$ 
The first approximation holds as $p_k={\beta}/{k}p_{k-1}$ for a Poisson distribution, i.e., $p_k$ increases in $k$ when $k<\beta$. The second approximation follows Stirling's formula that $\beta! \sim \sqrt{2\pi\beta}(\beta/e)^\beta$. 
\item When $\xi$ follows multinomial distribution with $\beta$ trials, following the similar argument, we have 
$$
\max_{k\in\mathbb{Z}}p_{k}\approx \frac{1}{\sqrt{\beta\pi/2}}.
$$
\end{itemize}
It means that when $\xi$ satisfies a Poisson distribution with a large arrival rate $\beta$ or when $\xi$ satisfies a multinomial distribution with a large number of trials $\beta$. The error term is of order $\cO(1/\sqrt{\beta})$.

For the $d$-dimensional case, the final error is multiplied by $\sqrt{d}$ as 
$$
\|\nabla g(x - \gamma v_F(x))-\nabla g(x)\|\leq \sqrt{d} \max_{i\in[d],k\in\mathbb{Z}}p_{ik} .
$$
{Via calculation, we have
\begin{equation}
\begin{aligned}
    \EE[F(\hat x^T)-F(x^*)]
\leq &
    \frac{\|u^1 -u^*\|^2}{2\gamma T}+  \gamma L_\phi^2 \Big(\frac{K^4 L_f^2}{16L_\phi^2}+ \lambda^2 D_\mathcal{X}^2\Big)  +4L_\phi D_\mathcal{X} \sqrt{d} \max_{i\in[d],k\in\mathbb{Z}}p_{ik} \Big(\frac{K^2 L_f}{8L_\phi}+ 2\lambda D_\mathcal{X}\Big) \\ &+ 2 L_\phi^2 D_\mathcal{X}^2 \lambda 
     +
    L_\phi^2 D_\mathcal{X} \frac{K L_f+2L_f}{\mu_g}\Big(1- \frac{\mu_g}{2L_\phi}\Big)^K.
\end{aligned}
\end{equation}
Picking $\lambda=\cO(1/\sqrt{T})$ and $\gamma = \cO(1/\sqrt{T})$, } it concludes the results shown in Theorem \ref{cor:msg_without_smoothness}.

\section{A Stochastic Gradient Method for Finite-dimension Convex Reformulation and Convergence Analysis}\label{sec:sgd_u}

In this section, we discuss how to solve the convex reformulation \eqref{eqn:trans-model-2-extension} via SAA and SGD. We first illustrate the key difference between SAA+SG and RSG/MSG.

\subsection{Algorithmic Design Difference between SAA+SG and RSG/MSG}
\label{appendix:algorithmic_difference}
Figure \ref{fig:RSG} illustrates the updating procedure of RSG and MSG. The gradient estimator $v$ of $F$ is constructed differently for RSG and MSG. Only arrows are executed in the algorithm while the dashed line between $x^t$ and $u^t$ represents the relationship $u^t=g(x^t)$ that is only used in the analysis. The update of SAA+SG is given in Figure \ref{fig:SAASGD} as a comparison.

\begin{figure}[t]
    \centering
    \includegraphics[width=0.9\linewidth]{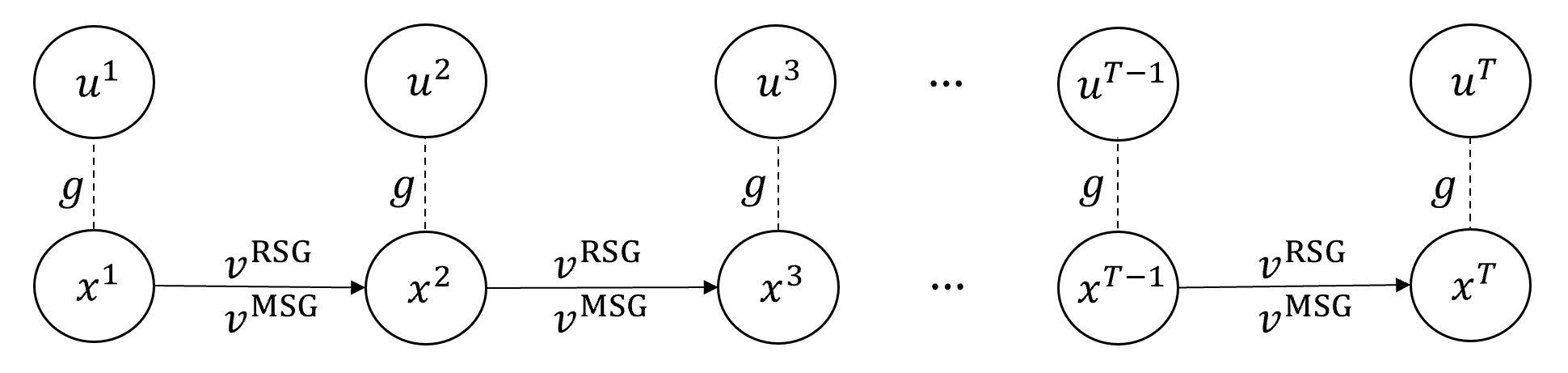}
    \caption{Illustration of RSG and MSG for Solving $F$}
    \label{fig:RSG}
    \vskip -0.2in
\end{figure}

\begin{figure}[t]
    \centering
    \includegraphics[width=0.9\linewidth]{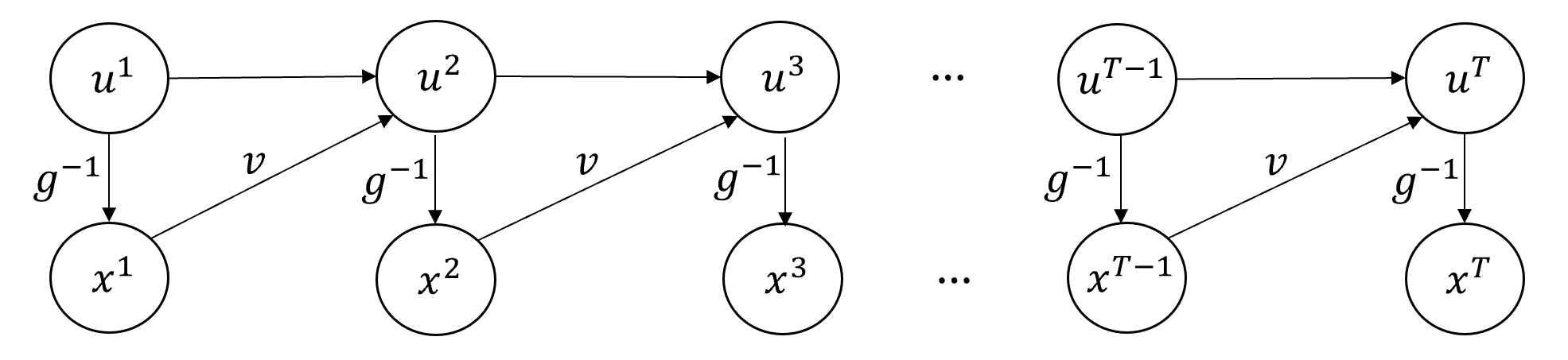}
    \caption{Illustration of SAA+SG for Solving $F$}
    \label{fig:SAASGD}
    \vskip -0.2in
\end{figure}

\subsection{Motivation of SAA+SG}
To perform projected stochastic gradient descent on $G$, based on Lemma \ref{lemma:gradient_advanced}, one needs to know $g^{-1}(u)$ to compute stochastic gradient estimator of $\EE[f(\phi(g^{-1}(u),\xi))]$ and needs to know the closed-form of $\mathcal{U}$ so as to perform projection onto $\mathcal{U}$. However, both $g^{-1}(u)$ and $\mathcal{U}$ involve unknown distribution $\PP(\xi)$.

A straightforward idea is to leverage SAA on the convex reformulation \eqref{eqn:trans-model-2-extension}, $\min_{u\in\mathcal{U}} G(u) = F(g^{-1}(u))$. Hence, one needs to build sample average estimators for the following three terms
\begin{itemize}
    \item $F(x) = \EE[f(\phi(x,\xi))]$;
    \item $g^{-1}(u)$ where $g(x) = \EE[\phi(x,\xi)]$;
    \item $\mathcal{U}=\{u\mid  \mathbb{E}[\phi_i(\munderbar{X}_i, \xi_i)]\leq u_i \leq \mathbb{E}[\phi_i(\bar{X}_i, \xi_i)], \text{ for any } i\in[d]\}$.
\end{itemize}
\begin{figure}
  \begin{center}
    \includegraphics[width=0.5\linewidth]{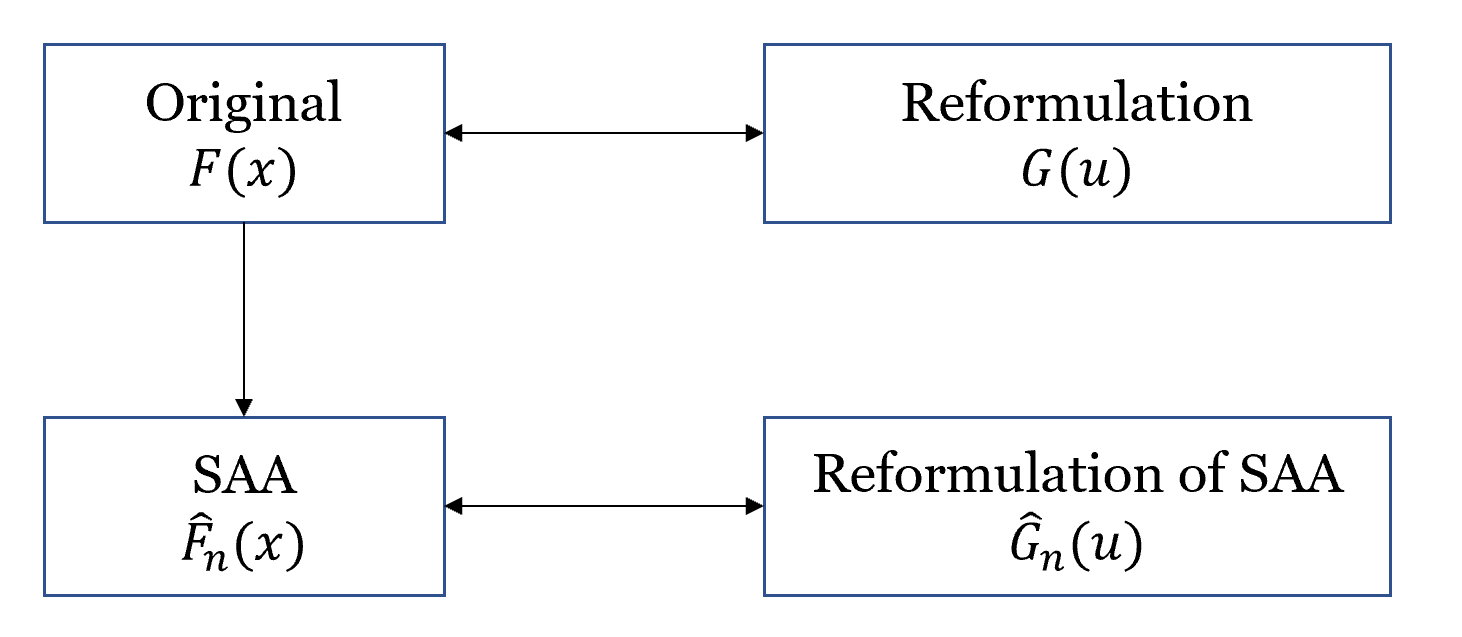}
  \end{center}
  \caption{Illustration of SAA+SG}
    \label{fig:idea_SAA}
\end{figure}
However, it is unclear whether we should 1) use the same set of samples to estimate these three terms, which might introduce undesired correlation when performing SGD to solve the empirical objective; or 2) use different sets of samples to estimate these three terms, which might lead to a potential nonconvex empirical objective. 

Instead, we follow a more principled way to construct a convex empirical objective. We use SAA to form an empirical objective $\hat F_n(x)$ for the original objective \eqref{problem:math_original_extension}. Then we utilize Proposition \ref{thm:equivalent-trans-extension} to form an equivalent convex reformulation $\hat G_n(u)$ of the empirical objective $\hat F_n(x)$. Next we solve $\hat G_n(u)$ using projected SGD. Figure \ref{fig:idea_SAA}  illustrates the key idea of the procedure. As a result, projected SGD is implementable on $\min_{\hat{\mathcal{U}}} \hat G_n(u)$ and as $n$ goes to infinity, $\hat F_n(x)$ is a good approximation of $F(x)$ according to law of large numbers. The formal definitions are in the following paragraph. We point out that such procedure coincidentally corresponds to using the same set of samples to estimate $F$, $g^{-1}$, and $\mathcal{U}$ and construct SAA for the convex reformulation as mentioned in the last paragraph.

The empirical optimization objective of \eqref{problem:math_original_extension} constructed via SAA is:
\begin{equation}
\label{problem:SAA}
    \min_{x\in\mathcal{X}} \hat F_n(x) \coloneqq \frac{1}{n}\sum_{j=1}^n f(\phi(x, \xi^j)),
\end{equation}
where $\{\xi^j\}_{j=1}^n$ are independent and identically distributed (i.i.d.) samples from $\PP(\xi)$.
Notice that the SAA problem \eqref{problem:SAA} can be interpreted as \eqref{problem:math_original_extension} with a uniform discrete distribution over $\{\xi^j\}_{j=1}^n$. Correspondingly, the SAA problem \eqref{problem:SAA} has a finite-dimensional convex reformulation by Proposition \ref{thm:equivalent-trans-extension}: 
\begin{equation}
\label{problem:SAA_finite_reformulation}
\min_{u\in\hat{\mathcal{U}}} \hat G_n (u) \coloneqq \frac{1}{n}\sum_{j=1}^n f( \phi(\hat g^{-1}(u), \xi^j)),
\end{equation}
where $\hat g(x) = \frac{1}{n}\sum_{j=1}^n \phi(x, \xi^j)$, $\hat{\mathcal{U}}=\{u\mid\frac{1}{n}\sum_{j=1}^n \phi(\munderbar{X}_i,\xi_i^j) \leq u_i\leq \frac{1}{n}\sum_{j=1}^n \phi(\bar{X}_i,\xi_i^j) \text{ for all } i\in[d]\}$, $\hat g^{-1}(u)=(\hat g^{-1}_1(u_1),\dots,\hat g^{-1}_d(u_d))^\top$ with $\hat g^{-1}_i(u_i)=\inf_{x\in[ \munderbar{X}_i,\bar{X}_i]} \{x\mid \hat g_i(x)\geq u_i\} \text{ for } i\in[d]$.

By Proposition \ref{thm:equivalent-trans-extension},  we know that $\hat G_n(u)$ is convex  and \eqref{problem:SAA_finite_reformulation} is equivalent to \eqref{problem:SAA}. 
From the classical SAA theory, to ensure an $\eps$- approximation error between SAA and the original objective, it requires a large number of samples $n=\tilde \cO(d\eps^{-2})$~\citep{kleywegt2002sample}.  Thus performing full-batch gradient descent on $\hat G_n(u)$ might not be efficient. Specifically, in NRM applications discussed in Section \ref{sec:application}, gradient descent requires solving $n$ linear programs at each iteration.  Instead, we perform stochastic gradient descent in the $u$-space on the empirical objective $\hat G_n(u)$. We denote such method as SAA+SG and the details are in Algorithm \ref{alg:sgd_on_g}.

In comparison, for classical stochastic optimization, it is generally unnecessary to first perform SAA then perform SGD for two reasons: 1) one can directly apply SGD; 2) the sample complexity of SGD  is better than that of SAA by a factor of $d$ in the convex setting, see a comparison between \citet{kleywegt2002sample} and \citet{nemirovski2009robust}.

\begin{algorithm}[ht]
	\caption{The Stochastic Gradient Method for Convex Reformulation (SAA+SG)}
	\label{alg:sgd_on_g}
	\begin{algorithmic}[1]
		\REQUIRE  Number of iterations $ T $, stepsizes $\{\gamma_t\}_{t=1}^{T}$, initialization point $u^1$, radius parameter $\delta_0$.
		\STATE Generate $n$ i.i.d. samples $\{\xi^j\}_{j=1}^n$ from $\PP(\xi)$.
		\STATE Set radius parameter $\delta = \min\{\delta_0, \frac{1}{2}\min_{i\in[d]} \frac{1}{n}\sum_{j=1}^n [\phi(\bar{X}_i,\xi_i^j) -\phi(\munderbar{X}_i,\xi_i^j)] \}$
		\FOR{$t=1$ to $T$ \do} 
		\STATE For given $u^t$, find $x^t\in\mathcal{X}$ such that $x^t=\hat g^{-1}(u^t)$. 
		\STATE Take a sample ${\xi^t}^\prime$ uniformly from $\{\xi^j\}_{j=1}^n$ and construct a gradient estimator
		\begin{center}
		$
        v(u^t) =\nabla \hat g(x^t)^{-\top} \nabla \phi(x^t,{\xi^t}^\prime)^\top\nabla f(\phi(x^t, {\xi^t}^\prime)).
		$
		\end{center}
		\STATE Update $ u^{t+1}=\Pi_{\hat{\mathcal{U}}_\delta}( u^t-\gamma_t v(u^t))$, where 
		\begin{center}
		$
        \hat{\mathcal{U}}_\delta = \{u\mid\frac{1}{n}\sum_{j=1}^n \phi(\munderbar{X}_i,\xi_i^j)+\delta \leq u_i\leq \frac{1}{n}\sum_{j=1}^n \phi(\bar{X}_i,\xi_i^j)-\delta \text{ for all } i\in[d]\}.
		$
		\end{center}
		\ENDFOR
		\ENSURE $\hat u^T$ and $\hat x^T$ where $\hat u^T = \frac{1}{T} \sum_{t=1}^T u^t$ and $\hat x^T=\hat g^{-1}(\hat u^T)$.
	\end{algorithmic}
\end{algorithm}

SAA+SG requires finding $x^t$ for a given $u^t$ at each iteration. Since $\phi(x,\xi)$ is a component-wise non-decreasing function in $x$ for any $\xi$, it is not very costly to find the corresponding $x^t$ for a given $u^t$.  Note that when updating $u^{t+1}$, we perform projection onto $\hat{\mathcal{U}}_\delta$ instead of $\hat{\mathcal{U}}$. This is to ensure that $[\nabla \hat g(x^t)]^{-1}$ is well-defined and we explain via the following example.

\begin{example}[Example of SAA+SG when $\phi(x,\xi)=x\wedge\xi$]
When $\xi$ is component-wise independent, consider the example when $\phi(x,\xi) = x\wedge \xi$ and $g(x) = \EE[ x\wedge \xi]$.  It is easy to verify that $\nabla_i g_i(x_i) = 1-\hat H_i(x_i)$, where $\hat H_i(\cdot)$ is the empirical CDF of the $i$-th coordinate of $n$ samples $\{\xi_i^j\}_{j=1}^n$. Suppose for some $t\in[T]$ and $i\in[d]$ that $u_i^t = \frac{1}{n}\sum_{j=1}^n \xi_i^j$ due to projection onto $\hat{\mathcal{U}}$. Hence, $x_i^t=\hat g_i^{-1}(u_i^t) = \max_{j\in[n]} \xi_i^j$. As a result, $\nabla \hat g_i(x_i^t) = 0$. Since $\nabla \hat g$ is a diagonal matrix, $[\nabla \hat g(x^t)]^{-1}$ is not well-defined.
\end{example}

Denote $ x_\mathrm{SAA}^*$ as the optimal solution of \eqref{problem:SAA};  $ u_\mathrm{SAA}^*$ as the optimal solution of \eqref{problem:SAA_finite_reformulation}; and $ u_\mathrm{SAA}^\delta$ as the optimal solution of $\min_{u\in\hat{\mathcal{U}}_\delta} \hat G_n (u) \coloneqq \frac{1}{n}\sum_{j=1}^n f( \phi(\hat g^{-1}(u), \xi^j))$.
The following theorem characterizes the approximation error of SAA on $F$ and expected error of projected SGD on $\hat G_n$.

\begin{theorem}
\label{lemma:saa+sgd}
The expected error of SAA+SG satisfies
$$
\EE [F(\hat x^T) - F(x^*)] \leq \EE [F(\hat x^T)-\hat F_n(\hat x^T)] + \EE [\hat G_n(\hat u^T) - \hat G_n(u^*_\mathrm{SAA})].
$$
Suppose Assumptions \ref{assumption:reformulation} and \ref{assumption:general} hold and $\phi_i(x_i,\xi_i)$ has left and right derivative in $x_i\in\mathcal{X}_i$ for any realization of $\xi_i$. 
If $f(\phi(x,\xi))$ is sub-Gaussian with a variance proxy $\sigma^2$, i.e., $\EE [\exp(t (f(\phi(x,\xi))-\EE f(\phi(x,\xi))))]\leq \exp\Big(t^2\sigma^2/2\Big)$ for any $x\in\mathcal{X}$, the approximation error of SAA satisfies
\begin{equation}
\label{eq:SAA+SG_SAA_error} 
\EE [F(\hat x^T)-\hat F_n(\hat x^T)]\leq \cO\Big(2\sqrt{\frac{d \log(D_\mathcal{X}\sqrt{n}) \sigma^2}{2n}} \Big)+\frac{2L_\phi L_f}{\sqrt{n}}.
\end{equation}
If $\hat g^{-1}$ is $L_{g^{-1}}$-Lipschitz continuous on $\hat g^{-1}(\hat{\mathcal{U}}_\delta)$,   letting $\gamma_t=\gamma$ and $\delta_0 = \frac{1}{\sqrt{dT}}$, the expected error of projected SGD on \eqref{problem:SAA_finite_reformulation} satisfies
\begin{equation}
\label{eq:SAA+SG_SGD_error}
\EE [\hat G_n(\hat u^T) - \hat G_n(u^*_\mathrm{SAA})]\leq \frac{\|u^1-u_\mathrm{SAA}^\delta\|^2}{\gamma T}+ \frac{\gamma }{T}\sum_{t=1}^T\EE\|v(u^t)\|^2 +  \frac{L_\phi L_f L_{g^{-1}}}{\sqrt{T}}.
\end{equation}
\end{theorem}
Note that the sub-Gaussian random function assumption is standard for SAA \citep{kleywegt2002sample}.  We point out that even if Assumptions \ref{assumption:general_2}(a) and (c) hold for $\phi(x,\xi)$ with $\xi$ under distribution $\PP(\xi)$, they may not hold for $\phi(x,\xi)$ with $\xi$ under the empirical distribution.
Note that \eqref{eq:SAA+SG_SAA_error} adopts from \citet{hu2020sample}. 
The first two terms in the right-hand-side of \eqref{eq:SAA+SG_SGD_error} also appear in classic projected SGD analysis~\citep{nemirovski2009robust} while the third terms comes from projection onto $\hat{\mathcal{U}}_\delta$ instead of $\hat{\mathcal{U}}$. We point out that in classical SGD analysis, one generally assumes that the gradient estimator $v(u^t)$ has $\cO(1)$ second moment for any $t\in[T]$. As a result, the sample and gradient complexity of classical SGD is $\cO(\eps^{-2})$ (by setting $\gamma = T^{-1/2}$ and $T = \cO(\eps^{-2})$) for convex objectives. Differently, such bounded $\cO(1)$ second moment condition might not hold for the gradient estimator $v(u)$ of SAA+SG for certain $\phi(x,\xi)$ that appears in supply chain and NRM applications, for instance when $\phi(x,\xi)=x\wedge\xi$. Proposition \ref{prop:saa+sgd} characterizes the second moment of the gradient estimator $v(x)$ when $\phi(x,\xi) = x\wedge \xi$ and demonstrates the corresponding sample and the gradient complexity of SAA+SG.

\begin{proposition}\label{prop:saa+sgd}
For $\phi(x,\xi) = x\wedge \xi$, under all conditions in Theorem \ref{lemma:saa+sgd}, we have
$\EE \|v(u)\|^2\leq ndL_f^2$  for any $u\in\mathcal{U}$. Setting $\gamma = (ndT)^{-1/2}$, then sample complexity $n$ of Algorithm \ref{alg:sgd_on_g} is $\tilde \cO(d\eps^{-2})$ and the gradient complexity is $\tilde \cO(d^2\eps^{-4})$.
\end{proposition}

The proposition shows that the second moment of the gradient estimator used in SAA+SG can be much larger than what classical SGD analysis normally assumes. Thus SAA+SGD method has a large gradient complexity meaning that the method takes a longer time to converge to a global optimal solution.   Such large second moment comes from estimating matrix inverse $[\nabla \hat g(x)]^{-1}$ via sample average. Note that one may not impose a variant of Assumption \ref{assumption:general_2}(b) that $\nabla \hat g(x)\succeq \mu_g I$ for any $x\in\mathcal{X}$ to control the second moment as the empirical distribution depends on generated samples. We point out that the upper bounds of the second moment derived in Proposition \ref{prop:saa+sgd} is based on  the worst $u\in\mathcal{U}$. For some $u\in\mathcal{U}$, the second moment $\EE\|v(u)\|^2$ could be bounded by $\cO(1)$. We leave the probabilistic characterization of the second moment of $\{v(u^t)\}_{t=1}^T$ for future investigation. In numerical experiments, we do observe that SAA+SG converges much slower than RSG and MSG, see e.g., Figure \ref{fig:comparison_of_algorithms}(a).

A natural question is whether we can design some alternative gradient estimator with a smaller second moment for $\phi(x,\xi)=x\wedge \xi$. The answer is yes. Utilizing the structure of $x\wedge\xi$, one can show that 
$$
[\nabla F(x)]_i =  (1-H_i(x_i))\EE_{\xi_{[-i]}} [\nabla f(x_i, x_{[-i]}\wedge {\xi}^\prime_{[-i]})]_i,
$$
where $[-i]$ denotes an index set $\{1,...,i-1,i+1,...,d\}$.
Therefore, for $x=g^{-1}(u)$, using the fact that $[\nabla g(x)]^{-\top}\nabla F(x) = \nabla G(u)$, we have
$$
[\nabla G(u)]_i = (1-H_i(x_i))^{-1} \EE_{\xi_{[-i]}} [(1-H_i(x_i))[\nabla f(x_i, x_{[-i]}\wedge {\xi}_{[-i]})]_i = \EE_{\xi_{[-i]}}  [\nabla f(x_i, x_{[-i]}\wedge {\xi}_{[-i]})]_i.
$$
where $(x_i, x_{[-i]}\wedge {\xi}^\prime_{[-i]}) = (x_1\wedge \xi_1,\ldots,x_{-1}\wedge \xi_{i-1},  x_i, x_{i+1}\wedge \xi_{i+1},\ldots,x_{d}^t\wedge \xi_{d})$.
Thus, one may construct a gradient estimator $\tilde v(u)$ with the $i$-th coordinate being
$$
[\tilde v(u)]_i = [\nabla f(x_i, x_{[-i]}\wedge {\xi}^\prime_{[-i]})]_i.
$$
The advantage of $\tilde v(u)$ is that 1) it does not need to know any information about $g$ to build a gradient estimator of $G$; 2) $\tilde v(u)$ has bounded second moment $dL_f^2=\cO(d)$ since $f$ is $L_f$-Lipschitz continuous. Thus the gradient complexity reduces to $\cO(d^2\eps^{-2})$. Note that with $\tilde v(u)$, we still need to first use SAA otherwise we cannot perform projection onto $\mathcal{U}$. The reduction in the gradient complexity via using $\tilde v(u)$ is not a free lunch. As the $i$-th coordinate of $\tilde v$ requires taking gradient of $f$ on the $i$-th input $(x_i, x_{[-i]}\wedge {\xi}^\prime_{[-i]})$. Therefore, to build  such an $\tilde v(u)$, it requires compute $\nabla f$ at $d$ different points $\{(x_i, x_{[-i]}\wedge {\xi}^\prime_{[-i]})\}_{i=1}^d$. Since estimating the gradient of $f$ in our NRM applications requires solving a linear program, it means that computing $\tilde v(u)$ require solving $d$ linear programs at each iteration which is much larger than solving $1$ linear program as required by SAA+SG. Hence, we do not intend to use the new estimator in practice.

\subsection{Proof of Theorem \ref{lemma:saa+sgd}}
\begin{proof1}
We decompose the expected error as follows:
\begin{equation}
\label{eq:error_decomposition_SAA}
\begin{aligned}
    &
\EE [F(\hat x^T) - F(x^*)]\\
    = &
\EE [F(\hat x^T)-\hat F_n(\hat x^T)+\hat F_n(\hat x^T)- \hat F_n(x_\mathrm{SAA}^*) + \hat F_n(x_\mathrm{SAA}^*) - \hat F_n(x^*) + \hat F_n(x^*) -F(x^*)]\\
    = &
\EE [F(\hat x^T)-\hat F_n(\hat x^T)+\hat G_n(\hat u^T)- \hat G_n(u_\mathrm{SAA}^*) + \hat F_n(x_\mathrm{SAA}^*) - \hat F_n(x^*) + \hat F_n(x^*) -F(x^*)]\\
    \leq &
 \EE [F(\hat x^T)-\hat F_n(\hat x^T)]+\EE_{\{\xi^j\}_{j=1}^n}[\EE[\hat G_n(\hat u^T)- \hat G_n(u_\mathrm{SAA}^*)|\{\xi^j\}_{j=1}^n]],
\end{aligned}
\end{equation}
where the second equality holds as $\hat G_n(u_\mathrm{SAA}^*) = \hat F_n(x_\mathrm{SAA}^*)$ by Proposition \ref{thm:equivalent-trans-extension}, and the inequality holds as $\hat F_n(x_\mathrm{SAA}^*) - \hat F_n(x^*)\leq 0$ and $\EE [\hat F_n(x^*) -F(x^*)]=0$. Note that $ \EE [F(\hat x^T)-\hat F_n(\hat x^T)]$ characterizes the approximation error of SAA and $\EE[\hat G_n(\hat u^T)- \hat G_n(u_\mathrm{SAA}^*)|\{\xi^j\}_{j=1}^n]]$ characterizes the error of SGD on $\hat G_n$.

\noindent \textbf{Approximation error of SAA:}
we first prove the approximation error of SAA using uniform convergence. Take a $\upsilon$-net $\{\tilde x^k\}_{k=1}^Q$ over $\mathcal{X}$ such that for any $x\in\mathcal{X}$, there exists a $k\in[Q]$ such that $\|\tilde x^k - x\|\leq \upsilon$. Such $\upsilon$-net exists when $Q=\cO\Big(\Big(\frac{D_\mathcal{X}}{\upsilon}\Big)^d\Big)$~\citep{kleywegt2002sample}. Denote $\bar x = \argmax_{x\in\mathcal{X}} [F(x)-\hat F_n(x)] $ and let $k_0\in [Q]$ be such that $\|\tilde x^{k_0}-\bar x\|\leq \upsilon $. We have the following result:
\begin{equation}
\label{eq:SAA_uniform}
\begin{aligned}
&
    \EE [F(\hat x^T)-\hat F_n(\hat x^T)] 
\leq 
    \EE \max_{x\in\mathcal{X}}[F(x)-\hat F_n(x)] 
= 
    \EE [F(\bar x)-\hat F_n(\bar x)]\\
= &
    \EE [F(\bar x)-F(\tilde x^{k_0})] + [F(\tilde x^{k_0})-\hat F_n(\tilde x^{k_0})] + [\hat F_n(\tilde x^{k_0}) - \hat F_n(\bar x)] \\
\leq &
    \EE [F(\tilde x^{k_0})-\hat F_n(\tilde x^{k_0})] + 2 L_\phi L_f \upsilon
\leq 
    \EE \max_{k_0\in[Q]} [F(\tilde x^{k_0})-\hat F_n(\tilde x^{k_0})] + 2 L_\phi L_f \upsilon,
\end{aligned}
\end{equation}
where the first inequality holds naturally, $\tilde x^{k_0}$ is the closest point in the $\upsilon$-net to $\bar x$, and the second inequality holds as $F(x)$ and $\hat F_n(x)$ are both $L_f L_\phi$-Lipschitz continuous and $\|\tilde x^{k_0} - \bar x\|\leq \upsilon$. Note that $\tilde x^{{k_0}}$ depends on the samples $\{\xi^j\}_{j=1}^n$. Thus $\tilde x^{{k_0}}$ is correlated with $\hat F_n$. To get rid of such dependence, we utilize the following argument for any $s>0$:
\begin{equation}
\label{eq:SAA_max}
\begin{aligned}
&
    \EE \max_{{k_0}\in[Q]} [F(\tilde x^{{k_0}})-\hat F_n(\tilde x^{{k_0}})]
= 
    \frac{1}{s}\log\Big(\exp\Big(s \EE \max_{{k_0}\in[Q]} [F(\tilde x^{{k_0}})-\hat F_n(\tilde x^{{k_0}})]\Big)\Big) \\
\leq &
    \frac{1}{s}\log\Big(\EE \exp\Big(s \max_{{k_0}\in[Q]} [F(\tilde x^{{k_0}})-\hat F_n(\tilde x^{{k_0}})]\Big)\Big)
=
    \frac{1}{s}\log\Big(\EE\max_{{k_0}\in[Q]} \exp\Big(s  [F(\tilde x^{{k_0}})-\hat F_n(\tilde x^{{k_0}})]\Big)\Big)\\
\leq &
    \frac{1}{s}\log\Big(\sum_{k=1}^Q\EE \exp\Big(s  [F(\tilde x^k)-\hat F_n(\tilde x^k)]\Big)\Big),
\end{aligned}
\end{equation}
where the first equality holds by definition, the first inequality holds by Jessen's inequality and the fact that exponential function is convex, the second equality holds as exponential function is strictly increasing, and the last inequality holds since exponential function is non-negative. After taking summation over $k\in[Q]$, each $\tilde x^k$ is from the $\upsilon$-net and is independent from $\hat F_n$.

By definition, we have $F(\tilde x^k)-\hat F_n(\tilde x^k)  = \frac{1}{n} \sum_{j=1}^n [\EE_\xi f(\phi(\tilde x^k,\xi)) - f(\phi(\tilde x^k,\xi^j))] $. Since each $\tilde x^k$ is independent of $\hat F_n$, we have
$
\EE_{\xi^j} [\EE_\xi f(\phi(\tilde x^k,\xi)) - f(\phi(\tilde x^k,\xi^j))] = 0.$ Utilizing the fact that $f(\phi(x,\xi))$ is sub-Gaussian for any $x\in\mathcal{X}$, we know that $\EE_\xi f(\phi(\tilde x^k,\xi)) - f(\phi(\tilde x^k,\xi^j))$ is a  zero-mean sub-Gaussian random variable. Therefore, it holds that
$$
\EE \exp\Big(s  [F(\tilde x^k)-\hat F_n(\tilde x^k)]\Big) \leq \exp\Big(\frac{s^2 \sigma^2}{2n}\Big) \ \text{ for any } k\in[Q].
$$
Combined with \eqref{eq:SAA_uniform} and \eqref{eq:SAA_max}, we have
$$
\begin{aligned}
&    
    \EE [F(\hat x^T)-\hat F_n(\hat x^T)] 
\leq 
   \frac{1}{s}\log\Big( Q\exp\Big(\frac{s^2 \sigma^2}{2n}\Big)\Big) + 2 L_\phi L_f \upsilon 
= 
    \frac{\log(Q)}{s} + \frac{s \sigma^2}{2n} + 2 L_\phi L_f \upsilon \\
= & 
    2\sqrt{\frac{\log(Q) \sigma^2}{2n}} + 2\frac{L_\phi L_f}{\sqrt{n}}
= 
    \cO\Big(\sqrt{\frac{d \log(D_\mathcal{X}\sqrt{n}) \sigma^2}{2n}} \Big)+2\frac{L_\phi L_f}{\sqrt{n}},
\end{aligned}
$$
where the second equality holds by setting $s = \sqrt{2\log(Q)n/\sigma^2}$ and $\upsilon= n^{-1/2}$, and the third equality uses the fact that $Q = \cO\Big(\Big(\frac{D_\mathcal{X}}{\upsilon}\Big)^d\Big)$. 

\textbf{Error of projected SGD on $\hat G_n(u)$}: next we demonstrate expected error of performing projected SGD on $\hat G_n(u)$. Since  $u^\delta_\mathrm{SAA}\in\hat{ \mathcal{U}}_\delta$, we have
\begin{equation*}
\begin{aligned}
&
    \EE \|u^{t+1} - u^\delta_\mathrm{SAA}\|^2 \\
= &
    \EE \|\Pi_{\hat{\mathcal{U}}_\delta}( u^t-\gamma v(u^t)) - \Pi_{\hat{\mathcal{U}}_\delta}(u^\delta_\mathrm{SAA})\|^2 \\
\leq &
    \EE \| u^t-\gamma v(u^t) - u^\delta_\mathrm{SAA}\|^2 \\
= &
    \EE \| u^t- u^\delta_\mathrm{SAA}\|^2 + \gamma^2 \EE \| v(u^t)\|^2 -2\gamma \EE(u^t-u^\delta_\mathrm{SAA})^\top v(u^t) \\
= &
    \EE \| u^t- u^\delta_\mathrm{SAA}\|^2 + \gamma^2 \EE \| v(u^t)\|^2 - 2\gamma \EE(u^t-u^\delta_\mathrm{SAA})^\top \nabla \hat G_n(u^t)\\
\leq &
    \EE \| u^t- u^\delta_\mathrm{SAA}\|^2 + \gamma^2 \EE \| v(u^t)\|^2 - 2\gamma \EE(\hat G_n(u^t)-\hat G_n(u^\delta_\mathrm{SAA})),
\end{aligned}
\end{equation*}
where the first inequality uses the fact that projection operator is non-expansive, the third equality uses the fact that $\EE [v(u^t)\mid u^t] = \nabla \hat G_n(u^t)$, and the second inequality uses convexity of $\hat G_n$. Rearranging terms and dividing $2\gamma$ on both sides, we have
$$
\EE[\hat G_n(u^t)-\hat G_n(u^\delta_\mathrm{SAA})]\leq \frac{\EE \| u^t- u^\delta_\mathrm{SAA}\|^2 - \EE \|u^{t+1} - u^\delta_\mathrm{SAA}\|^2}{2\gamma} + \frac{\gamma \EE \| v(u^t)\|^2}{2}.
$$
Summing up from $t=1$ to $t=T$ and dividing $T$ on both sides, we have
$$
\EE[\hat G_n(\hat u^T)-\hat G_n(u^\delta_\mathrm{SAA})] \leq \frac{1}{T}\sum_{t=1}^T \EE[\hat G_n(u^t)-\hat G_n(u^\delta_\mathrm{SAA})] \leq \frac{ \| u^1- u^\delta_\mathrm{SAA}\|^2}{2\gamma T} + \frac{1}{T}\sum_{t=1}^T \frac{\gamma \EE \| v(u^t)\|^2}{2},
$$
where the first inequality uses the definition of $\hat u^T$, the convexity of $\hat G_n(u)$, and Jensen's inequality. On the other hand, we have
$$
\begin{aligned}
&
    \EE [\hat G_n(u^\delta_\mathrm{SAA})-\hat G_n(u^*_\mathrm{SAA})]  
=  
    \EE [\hat G_n(u^\delta_\mathrm{SAA})-\hat G_n(\Pi_{\hat{\mathcal{U}}_\delta} (u^*_\mathrm{SAA}))] + \EE [\hat G_n(\Pi_{\hat{\mathcal{U}}_\delta} (u^*_\mathrm{SAA}))-\hat G_n(u^*_\mathrm{SAA})]  \\
\leq &
    \EE [\hat G_n(\Pi_{\hat{\mathcal{U}}_\delta} (u^*_\mathrm{SAA}))-\hat G_n(u^*_\mathrm{SAA})]
\leq 
    L_\phi L_f L_{g^{-1}} \EE \|\Pi_{\hat{\mathcal{U}}_\delta} (u^*_\mathrm{SAA}) - u^*_\mathrm{SAA}\|\\
\leq &
    L_\phi L_f L_{g^{-1}} \EE \delta \sqrt{d}
\leq 
    L_\phi L_f L_{g^{-1}} \delta_0 \sqrt{d}
= 
    L_\phi L_f L_{g^{-1}} \frac{1}{\sqrt{T}}.
\end{aligned}
$$
where the first inequality holds by optimality of $u^\delta_\mathrm{SAA}$,  the second inequality holds by Lipschitz continuity of $\phi$, $f$, and $\hat g^{-1}$, and the third inequality holds by definition of $\hat{\mathcal{U}}_\delta$ and $\delta_0 = \frac{1}{\sqrt{d T}}$.
\end{proof1}

\subsection{Proof of Proposition \ref{prop:saa+sgd}}

\begin{proof1}
When $\phi(x,\xi) = x\wedge \xi$, it holds that $
\nabla_i \hat g_i(x_i)  = 1- \hat H_i(x_i).
$
We further have
$$
\begin{aligned}
    \EE[\|v(u)\|^2|\{\xi^j\}_{j=1}^n] 
= &
    \EE_{{\xi}^\prime}\Big[\|\nabla \hat g(x)^{-\top}\nabla \phi(x,{\xi}^\prime)^\top\nabla f(\phi(x, {\xi}^\prime))\|^2  |\{\xi^j\}_{j=1}^n\Big]\\
= & 
    \EE_{{\xi }^\prime}\Big[\sum_{i=1}^d (1-\hat H_i(x_i))^{-2} \mathbf{1}(x_i\leq {\xi_i}^\prime)  [\nabla f(x \wedge {\xi}^\prime)]_i^2  |\{\xi^j\}_{j=1}^n\Big],
\end{aligned}
$$
where $\hat H_i$ is the empirical CDF of $\{\xi_i^j\}_{j=1}^n$ and the second equality holds by the definition of $[\nabla \hat g(x)]^{-1}$ and the fact that $\nabla \hat g(x)$ is a diagonal matrix.
Without loss of generality, assume that the inequality $\xi_i^1<\xi_i^2<\ldots<\xi_i^n$ holds for some $i\in[d]$. 
When $\xi_i^{j}\leq x_i<\xi_i^{j+1}$ for $j=1,2,\ldots,n-1$, it holds that
$$
\begin{aligned}
    \EE_{{\xi}^\prime}\Big[[\nabla f(x \wedge {\xi}^\prime)]_i^2 \mathbf{1}(x_i\leq {\xi_i}^\prime) (1-\hat H_i(x_i))^{-2} |\{\xi^j\}_{j=1}^n\Big] 
= &
    \EE_{\xi^\prime_{[-i]}}[\nabla_i f(x_i, x_{[-i]}\wedge\xi^\prime_{[-i]})]^2\frac{n-j}{n}(1-\hat H_i(x_i))^{-2}  \\
= &
    \frac{n}{n-j} \EE_{\xi^\prime_{[-i]}}[\nabla_i f(x_i, x_{[-i]}\wedge\xi^\prime_{[-i]})]^2,
\end{aligned}
$$
where $\xi_{[-i]}^\prime$ denotes  $\xi^\prime$ excluding the $i$-th coordinate, and $[\nabla_i f(x_i, x^\prime_{[-i]}\wedge\xi^\prime_{[-i]})]$ denotes the $i$-th coordinate of the gradient of $f$ on point $(x_{1}\wedge\xi^\prime_{1},...,x_{i-1}\wedge\xi^\prime_{i-1}, x_i, x_{i+1}\wedge\xi^\prime_{i+1},...,x_{d}\wedge\xi^\prime_{d})$.  The first equality holds as ${\xi_i}^{\prime}$ is selected uniformly from $\{\xi_i^j\}_{j=1}^n$ and the second equality holds as $(1-\hat H_i(x_i))^{-1} = \frac{n}{n-j}$. As a result, we have
$$
\EE[\|v(u)\|^2|\{\xi^j\}_{j=1}^n] \leq   \sum_{i=1}^d n \EE_{\xi^\prime_{[-i]}}[\nabla_i f(x_i, x_{[-i]}\wedge\xi^\prime_{[-i]})]^2 \leq n dL_f^2.
$$
Taking full expectation, we have $\EE\|v(u)\|^2=\EE_{\{\xi^j\}_{j=1}^n} \EE[\|v(u)\|^2|\{\xi^j\}_{j=1}^n] \leq nd L_f^2$. Together with Lemma \ref{lemma:saa+sgd}, we have
$$
\EE [F(\hat x^T) - F(x^*)]\leq \cO\Big(2\sqrt{\frac{d \log(D_\mathcal{X}\sqrt{n}) \sigma^2}{2n}} \Big)+\frac{2L_\phi L_f}{\sqrt{n}}+\frac{\|u^1-u_\mathrm{SAA}^\delta\|^2}{\gamma T}+ \gamma n dL_f^2 +  \frac{L_\phi L_f L_{g^{-1}}}{\sqrt{T}}.
$$
Setting $\gamma = (ndT)^{-1/2}$, $n = \tilde \cO(d\eps^{-2})$, $T =
\cO(nd\eps^{-2})=\tilde \cO(d^2\eps^{-4})$, we have 
$$
\EE [F(\hat x^T) - F(x^*)]\leq \cO(\eps).
$$
Thus for $\phi(x,\xi) = x\wedge \xi$,  the sample complexity of Algorithm \ref{alg:sgd_on_g} is $n=\tilde \cO(d\eps^{-2})$ and the gradient complexity of Algorithm \ref{alg:sgd_on_g} is $T = \tilde \cO(d^2\eps^{-4})$.
\end{proof1}

{\section{Application: Assemble-to-Order Systems}
\label{appendix:application_ATO}
We show how the Assemble-to-Order system with a random capacity from \citep{chen2018preservation} can be formulated as a speical case of problem \eqref{problem:math_original_extension}.

Consider a dynamic ATO system with $T$ periods. There are $d$ components indexed by $i\in[d]$ and $n$ products indexed by $j\in[m]$. After observing inventory levels $y=(y_1,...,y_d)^\top$, firm decides the up-to-inventory levels $x=(x_1,...,x_d)^\top$. The replenishment lead time is zero, and the delivered quantity is truncated by a random supply capacity $\xi$. The random demand for products is presented as $D=(D_1,...,D_m)^\top$, which is independent of capacities. The bill of materials is defined by $d\times m$ matrix $A$. The unit ordering, holding, and shortage cost are denoted by vectors $c,h,b$. Unsatisfied demand is assumed to be lost, and the objective is to minimize the expected discounted cost with discount factor $\alpha$. Then, the bellman equation can be written as,
$$f_t(y)=\min_{x\geq y}\EE[c^\top(x\wedge(y+\Xi)-y)]+\EE[g_t(x\wedge(y+\Xi)|D)],$$
where $$g_t(z|d)=\min_{u:Au\leq z,0\leq u\leq d}\{h^\top(z-Au)+b^\top(d-u)+\alpha f_{t+1}(z-Au)\}.$$
The boundary condition is $f_{T+1}=0$. Notation $z$ is the inventory level after delivery, and vector $u$ represents the assembled product quantities. For any $t$, \cite{chen2018preservation} show that the cost-to-go function $f_t(y)$ is convex in $y$, and the function $G_t(z):=c^\top z+\EE[g_t(z|D)]$ is convex in $z$, which means the nonconvex minimization problem at each period in ATO dynamic formulation can be solved by our algorithms.
}

\section{Further Discussion for NRM}\label{sec:appendix_NRM_alg}

\subsection{Passenger NRM Modeling}
Passenger NRM is a special class of air-cargo NRM introduced in Section \ref{sec:air_cargo_model}, with one-dimensional capacity (e.g., seats on the plane), deterministic consumption (e.g., one passenger takes one seat of the airplane), and fixed route (e.g., passenger takes the route in the request). We introduce the following notations for the passenger NRM: $A=(a_{ij})_{i\in[m],j\in [d]}$ is the consumption matrix, where each unit of demand class $j$ consumes $a_{ij}$ units of the inventory class $i$. Then the passenger NRM problem under booking limit control policy can be written as follows.
\begin{equation}\label{eqn:basic_model}
	\begin{aligned}
		\max_{x\geq 0} ~~& \mathbb{E}_{\tilde{D}}[f(x\wedge\tilde{D})],
	\end{aligned}
\end{equation}
where $f(x)=r^\top x-\mathbb{E}_{\tilde{Z}(x),\tilde{c}}[\Gamma(\tilde{Z}(x),\tilde{c})]$ and 
\begin{equation}
\label{eq:primal_form_Gamma}
	\begin{aligned}
	\Gamma(z,c)=~\min_w \{&l^\top (z-w) \mid 
		Aw\leq c;
		~0\leq w\leq z\}.
	\end{aligned}
\end{equation}
Similar to the notation in Section \ref{sec:air_cargo_model}, $r$ denotes the revenue per-unit vector, $\tilde Z(x)$ denotes the show-ups given $x$ accepted reservations at the reservation stage, $\tilde c$ is the random capacity, $l$ denotes the penalty for rejecting accepted reservations.

\subsection{Structural Properties of NRM Models}
\label{sec:structural_properties}
In this section, we first reproduce the standard results on the structural properties of our booking limit models (\ref{eqn:basic_model}) and (\ref{eqn:ari_cargo_two_dim}). 

\begin{lemma}\label{lem:convexity}
For our booking limit model, we have the following structure properties,
\begin{enumerate}[(I)]
    \item In model (\ref{eqn:basic_model}), $\Gamma(z,c)$ is convex in $z$ (and $c$).
    \item In model (\ref{eqn:ari_cargo_two_dim}), $\Gamma(z,W, V, c_w,c_v)$ is convex in $z$ (and $c$).
\end{enumerate}
\end{lemma}

\begin{lemma}\label{lem:compo-convexity}
In both model (\ref{eqn:basic_model}) and model (\ref{eqn:ari_cargo_two_dim}), if the random show-up $\tilde{Z}_i(x_i)$ follows Poisson distribution with coefficient $p_ix_i~i=1,\dots,d$,  then $f(x)$ is component-wise concave in $x$. If all reservations show up, i.e., $\tilde{Z}(x)=x$, then $f(x)$ is concave in $x$.
\end{lemma}

Lemma \ref{lem:convexity} follows from standard linear programming theory. The proof of Lemma \ref{lem:compo-convexity} can be found in \cite{karaesmen2004overbooking}. Lemma \ref{lem:compo-convexity} claims that in both passenger and air-cargo NRM, the function $f(x)$ is component-wise concave when the random show-up follows Poisson distribution, and concave in the all-show-up case. 

\subsection{Stochastic Gradient of $f$}
\label{sec:gradient_calculation}
Due to random capacity in NRM problems, computing the exact gradient of $f$ is unpractical. In this subsection, we discuss how to compute the stochastic gradient of $f$ to facilitate the implementation of the proposed stochastic gradient-based algorithms in NRM applications. We reproduce the unbiased gradient estimator from \cite{karaesmen2004overbooking} for completeness. For simplicity, we focus on the stochastic gradient construction of the passenger NRM model (\ref{eqn:basic_model}). The procedure for the air-cargo model (\ref{eqn:ari_cargo_two_dim}) is similar and we directly give its gradient estimator construction.

Recall that $f(x) = r^\top x-\mathbb{E}_{\tilde{Z}(x),\tilde{c}}[\Gamma(\tilde{Z}(x),\tilde{c})]$.  We derive the stochastic gradient of $\mathbb{E}_{\tilde{Z}(x),\tilde{c}}[\Gamma(\tilde{Z}(x),\tilde{c})]$ with respect to $x$, and the remaining is straightforward.
First, we represent $\Gamma(z,c)$ in the dual form of \eqref{eq:primal_form_Gamma}.
\begin{equation}\label{eqn:dual_of_V0}
\begin{aligned}
    \Gamma(z,c)=\max_{v_1, v_2}&~l^Tz-(c^Tv_1+z^Tv_2)\\
    \textit{s.t.}&~A^Tv_1+ v_2\geq l\\
     &~v_1,v_2\geq 0.
\end{aligned}
\end{equation}
Thus $z$ only appears in the objective.

Second, we calculate the partial derivative of $\EE[ \Gamma(\tilde Z(x),\tilde c)]$ when $\tilde{Z}$ satisfies a Poisson distribution.
\begin{align*}\nonumber
& \frac{\partial}{\partial x_i}\mathbb{E}_{\tilde{Z},\tilde{c}}[\Gamma(\tilde{Z}(x),\tilde{c})]
=\lim_{h\rightarrow0}\frac{1}{h}\Big[\mathbb{E}_{\tilde{Z},\tilde{c}}[\Gamma(\tilde{Z}(x+e_ih),\tilde{c})]-\mathbb{E}_{\tilde{Z},\tilde{c}}[\Gamma(\tilde{Z}(x),\tilde{c})] \Big],
\end{align*}
where $e_i$ denotes the $i$-th unit vector in $\mathbb{R}^d$. Let $Y_i(h)$ denote a Poisson random variable with mean $p_ih$ that is  independent of $\tilde{Z}(x)$. We can represent $\mathbb{E}_{\tilde{Z},\tilde{c}}[\Gamma(\tilde{Z}(x+e_ih),\tilde{c})]$ as follows.

\begin{align*}\nonumber
& 
    \mathbb{E}_{\tilde{Z},\tilde{c}}[\Gamma(\tilde{Z}(x+e_ih),\tilde{c})]\\
=&
    \mathbb{E}_{\tilde{Z},\tilde{c},Y_i}[\Gamma(\tilde{Z}(x)+e_iY_i(h),\tilde{c})]\\
= &
    \mathbb{E}_{\tilde{Z},\tilde{c},Y_i}[\Gamma(\tilde{Z}(x),\tilde{c})|Y_i(h)=0]P(Y_i(h)=0)+\mathbb{E}_{\tilde{Z},\tilde{c},Y_i}[\Gamma(\tilde{Z}(x)+e_i,\tilde{c})]P(Y_i(h)=1)+o(h),
\end{align*}
where the first equality uses the property that sum of independent Poisson distribution is still Poisson, and the second equality holds by the law of total expectation and probability mass function of Poisson distribution. Since $Y_i(h)$ is a Poisson random variable with mean $p_ih$, we have $P(Y_i(h)=1)=p_ihe^{p_ih}=p_ih+o(h)$ and $P(Y_i(h)=0)=e^{p_ih}=1-p_ih+o(h)$. As a result, we can represent the partial derivative as follows.
\begin{align*}
    \frac{\partial}{\partial x_i}\mathbb{E}_{\tilde{Z},\tilde{c}}[\Gamma(\tilde{Z}(x),\tilde{c})]
= &
    \lim_{h\rightarrow0}\frac{1}{h}p_ih\Big[\mathbb{E}_{\tilde{Z},\tilde{c}}[\Gamma(\tilde{Z}(x)+e_i,\tilde{c})-\Gamma(\tilde{Z}(x),\tilde{c})]+o(h) \Big]\\
= &
    p_i\Big[\mathbb{E}_{\tilde{Z},\tilde{c}}[\Gamma(\tilde{Z}(x)+e_i,\tilde{c})-\Gamma(\tilde{Z}(x),\tilde{c})]\Big].
\end{align*}
Thus, an unbiased stochastic gradient estimator of $\frac{\partial}{\partial x_i}\mathbb{E}_{\tilde{Z},\tilde{c}}[\Gamma(\tilde{Z}(x),\tilde{c})]$ is $p_i(\Gamma(Z+e_i,c)-\Gamma(Z,c))$ for all $i\in[d]$ given realizations $Z$ and $c$ of $\tilde Z(x)$ and $\tilde c$, respectively.

Algorithm \ref{alg:gradient_air_cargo_NRM} demonstrates how to compute the unbiased stochastic gradient estimator of $f$ for the air-cargo NRM setting. 

\begin{algorithm}[t]
	\caption{Unbiased Stochastic Gradient Estimator for Air-cargo NRM}
	\label{alg:gradient_air_cargo_NRM}
	\begin{algorithmic}[1]
		\REQUIRE Parameters $p_i=1,~i\in[d]$ if all reservations show up. Booking limit $x$.
		\STATE Draw samples $D, W, V, c_w, c_v$.
		\STATE Draw a sample of show-ups $Z$. The $i$-th coordinate of $Z$ satisfies
		\begin{equation*}
		\begin{aligned}
		&\text{ When there is no-shows: } Z_i  \sim \text{Poisson distributed with mean } p_i(x_i\wedge D_i). \\
		&\text{ When all accepted reservations show-up: } Z_i = x_i\wedge D_i.
		\end{aligned}
		\end{equation*}
		\STATE Construct the gradient estimator $v^\mathrm{NRM}(x)=(v_1^\mathrm{NRM}(x_1),\dots,v_d^\mathrm{NRM}(x_d))^T$ where
		\begin{center}
		$
		v_i^\mathrm{NRM}(x)=\mathbf{1}\{x_i\leq D_i\}(r_i(W,V)-p_i(\Gamma(Z+e_i,W,V, c_w,c_v)-\Gamma(Z,W,V, c_w,c_v))) \text{ for } i\in[d].
		$
		\end{center}
		\ENSURE $v^\mathrm{NRM}(x)$.
	\end{algorithmic}
\end{algorithm}

\subsubsection{Practical Computational Issues}
{\color{black}
For both models (\ref{eqn:ari_cargo_two_dim}) and (\ref{eqn:basic_model}), to obtain $\Gamma(Z,W,V,c_w,c_v)$ and $\Gamma(Z,c)$ with given realizations $Z$ and other random variables, one needs to solve one LP. Since an unbiased stochastic gradient estimator of $\mathbb{E}[\Gamma(\tilde{Z}(x),\cdot)]$ requires knowledge of $\Gamma(Z+e_i,\cdot)$ for $i\in[d]$ and $\Gamma(Z,\cdot)$, obtaining such an unbiased gradient estimator requires  solving $d+1$ LPs. When $d$ is large, it could still be costly.

To overcome such computational burden, we use the optimal dual solution $v_{2i}^*,~i\in[d]$, associated with the constraint $w_i\leq z_i$ in the LP to construct an estimator $r_i-l_i+v_{2i}^*$ for given realizations $Z,c$. In the all-show-up case when the dual form admits a unique solution, it holds that $f$ is continuously differentiable and admits an unbiased gradient $r-l+v_{2}^*$. As a result, we only need to solve one LP to obtain the stochastic gradient of $f$ rather than $d+1$ LPs.

When the show-up is Poisson distributed with $p<1$, we still heuristically use $l_i-v_{2i}^*$ to approximate $\Gamma(Z+e_i, c)-\Gamma(Z,c)$ for $i\in[d]$ for reducing the computational cost in our numerical experiments.
}

\subsection{Discussions on Integer Booking Limits and Poisson Show-ups}
\label{sec:discussion_integer_poisson}
{\color{black} We focus on continuous booking limit decisions and Poisson random show-ups in Section \ref{sec:application}. However, the booking limit is generally in the integer space, and the random show-ups follow a binomial distribution in practice. In this subsection, we discuss such inconsistency and how we handle the integer booking limit setting with binomial show-ups, i.e., the setting in our numerical experiments. 

We first discuss the integer booking limits. During implementation, we keep continuous booking limits $\{x^t\}_{t=1}^T$ when running the algorithm and only round the final output of the algorithm to the nearest integer value. This simple rounding procedure works well in our reported numerical experiments with large demands. Although one may identify a better integer solution by enumerating integer solutions near the converging solution through sample average evaluation of the revenue, this procedure is still heuristic, and the exhaustive searching requires $O(2^d)$ times revenue evaluation. On the other hand, for numerical instances with few total demands but a large number of demand classes, when the average demand for each demand class is small (maybe even smaller than $1$), our booking limit control with the simple rounding procedure may not work well. This situation typically happens when there are too many fare classes for a given origin-destination flight, and each fare class has few demands. One heuristic solution is by nesting and collecting multiple fare classes with similar prices and same origin-destination as a new demand class with the replaced mean price. Essentially, our continuous optimization model can be regarded as a fluid relaxation of the integer booking limit model. Thus when the optimal booking limit has large values, the revenue incurred by the fractional part becomes negligible in practice.

Next, we discuss the Binomial random show-ups. Although Poisson can be regarded as a continuous approximation to the binomial show-ups, one drawback of Poisson is that the realized show-ups can be greater than the accepted reservations. Due to this drawback, we need to heuristically adapt our algorithm to binomial show-up case. Note that binomial show-ups require the accepted reservations to be an integer number as a parameter input. In contrast, Poisson show-ups only depend on a mean parameter, which can be non-integer. 
However, in our model formulation, we consider a continuous booking limit $x$, resulting in fractional (non-integer) acceptance $x\wedge\xi$, which is not an ideal parameter input to binomial distribution. To address this issue, we follow the convention in \cite{erdelyi2010dynamic}: with probability $\left \lfloor{x}\right \rfloor+1-x$, the random show-up follows binomial distribution with parameter $(\left \lfloor{x}\right \rfloor,p)$; otherwise, the random show-up follows binomial distribution with parameter $(\left \lfloor{x}\right \rfloor+1,p)$. This procedure guarantees the same expected show-ups. 
}

\section{Appendix for Computation Experiments in Section \ref{sec:numerics}}\label{appendix:numerics}

\subsection{Implementation Details of Benchmark Strategies}\label{appendix:computation_details}

\noindent \textbf{Deterministic Linear Programming (DLP).} We first introduce the standard bid price control policy obtained from the DLP for completeness. The DLP method is a standard method for NRM \citep{talluri1998analysis} and serves as the most famous benchmark. The DLP method solves \eqref{eqn:basic_model} with all random variables replaced by the corresponding expectations, leading to time-independent control policies. The mathematical formulation of DLP is as follows.

\begin{equation}\label{eqn:DLP}
\begin{aligned}
    \max_{w, x\geq 0} ~~~& r^Tx-l^T(p x-w)\\
     \text{s.t.} ~~~& Aw\leq \mathbb{E}[\tilde{c}]\\
     & x\leq \mathbb{E}[\tilde{D}]\\
     & w\leq p x,
\end{aligned}
\end{equation}
where $p$ is the show-up probability rate, decision $x$ is the total number of accepted reservations, and $w$ is the number of accommodated passengers. The revenue collected during the reservation period is $r^Tx$, and the loss induced by rejecting $px-w$ bookings is $l^T(px-w)$. The first constraint specifies the capacity constraint in the expected sense. The second constraint ensures that accepted reservations are no more than the expected demand. Due to the cancellations and no-shows,  only $px$ out of the total $x$ accepted reservations show up. The third constraint means that the number of accommodated passengers is no more than the number of show-up passengers. 

One can use the dual solution of the DLP to construct a policy for accepting and rejecting booking requests. Take the bid price control as an example. Let $\{\pi_j^*:j\in[m]\}$ be the optimal dual solution associated with the capacity constraint $Aw\leq \mathbb{E}[\tilde{c}]$. One can use $\pi^*$ to construct a bid price control policy. If the revenue from a request exceeds the sum of the expected opportunity cost of capacities consumed by this request, i.e., $r_i\geq \sum_{j=1}^m a_{ji}\pi_j^*$,  we accept the request. 

In addition to the bid price control, the primal solution $x$ of DLP can serve as booking limits. Moreover, the booking limit control policy basically accepts all requests until the limits are met. The optimal objective value of DLP is an upper bound of the optimal revenue \citep{erdelyi2009separable}. The formulation of DLP can be easily extended to more complicated settings, including the air-cargo network setup. However, due to its static decision rule and relatively poor performance (comparing to the more sophisticated bid control policies obtained from dynamic programmings as we will discuss later), we only use DLP as one of the benchmarks in the passenger network revenue management and neglect this method in the air-cargo variants.

\noindent \textbf{RSG, MSG and SAA+SG.} First, we specify the common parts shared by RSG, MSG, and SAA+SG, including the stochastic gradient construction of $f(x)$ in the NRM problem, step size, and stopping criteria. As discussed in Appendix \ref{sec:gradient_calculation}, we heuristically use the optimal dual value associated with the constraint $w\leq z$ to approximate the stochastic gradient, which reduces solving $d+1$ LPs to solving $1$ LP at each iteration. The computation indicates that such approximation performs well in the NRM instances since it induces a similar trajectory of $\{x^t\}_{t=1}^T$ to the unbiased gradient estimator. Thus, throughout all of our numerical experiments via RSG, MSG, and SAA+SG, we stick to this dual approximated stochastic gradient. The step size is set as $\gamma_t=a/\sqrt{t}$, where $a$ is tuned for specific instances. As for stopping criteria, we compute the Euclidean distance between two consecutive average solutions of $x^t$ over 100 iterations, and the algorithm stops when this Euclidean distance is less than 0.5 or the number of iterations exceeds the maximum 5,000. Except that for Figure \ref{fig:comparison_of_algorithms}, we stop these algorithms at 3000-th iteration for illustration. In general, three algorithms converge within 3,000 iterations. As for the binomial show-ups, we follow the discussion in Online Appendix \ref{sec:discussion_integer_poisson}. As for algorithm specific parameters, we set the regularization term $\lambda=1/t$ in RSG. In SAA+SG, we randomly sample $1,000$ i.i.d. samples for the sample average construction. In MSG, we set $K=10$. 

We also want to remark that the final convergent continuous booking limit is rounded to the nearest integer value because we do not allow fractional or probabilistic acceptance over all numerical experiments for a fair comparison.

\noindent \textbf{Dynamic Programming Decomposition (DPD).} As mentioned in the literature review, dynamic programming decomposition is widely used to derive the bid-price-based control policies. We compare our methods with the DPD method proposed by \citet{erdelyi2010dynamic} for two reasons: 1) their DPD method considers the random show-up, which is similar to our setting; 2) they provide a public dataset of NRM instances with good quality for a fair comparison. Next, we present their DPD method to illustrate how the decomposition  deals with the curse of dimensionality. The basic idea of DPD is to decompose the NRM problem with $m$ flight legs into $m$ single-leg dynamic models. Formally, the decomposed model becomes
\begin{align*}\nonumber
&V_{t,j}(w)=~\max_{x\in \{0,1\}^d}\sum_{i:a_{j,i}=1} \lambda_{t,i}\{\hat{R}_{i,j}x_i+V_{t-1,j}(w+x_i)\}, ~~\forall 1\leq t\leq T\\
&V_{0,j}(w)=~-\mathbb{E}_{\hat{Z}}[\Gamma_j(\hat{Z}(\alpha_j w))]
\end{align*}
For each $j\in[m]$, the $\Gamma_j$ function is 
\begin{align*}
~\Gamma_j({Z})=~\min & ~\sum_{i:a_{j,i}=1}\hat{L}_{i,j}g_i\\
\textit{s.t.}&~ ~\sum_{i=1}^d a_{j,i}({Z}_i-g_i)\leq c_j\\
&~ g_i\leq {Z}_i,~\forall i.
\end{align*}
One accepts the reservation request only when the revenue of the reservation is more than the implicit cost (revenue loss of the value-to-go function by accepting the request). For a detailed description of the method, please refer to \cite{erdelyi2010dynamic}. This method directly applies to the passenger network variant with deterministic capacity. In our passenger network variant with random capacity, we follow standard methodology to revise the boundary function $V_{0,j}(w)=-\mathbb{E}_{\hat{Z}}[\Gamma_j(\hat{Z}(\alpha_j w))]$ as $V_{0,j}(w)=-\mathbb{E}_{\hat{Z},\hat{c}}[\Gamma_j(\hat{Z}(\alpha_j w),\hat{c})]$ and incorporate the random capacity. 

Since this method does not explicitly consider the random consumption, two-dimensional capacity, and routing flexibility in the air-cargo variant, we only report the numerical results of the DPD method in this passenger NRM case. 

\noindent \textbf{Air Cargo Dynamic Programming Decomposition \citep{barz2016air} (ACDPD).} To compare our booking control policy in the air-cargo NRM setting, we introduce the following state-of-the-art DPD method specifically designed for air-cargo NRM, denoted as ACDPD~\citep{barz2016air}. ACDPD is a variant of DPD, the policy also bases on the bid price control, i.e., if the revenue of the incoming reservation is larger than the total bid price of all inventory classes, the airline accepts the reservation. In the air-cargo network variant, two-dimensional capacity is easy to handle as one may treat the air-cargo NRM as two different inventory classes sharing the same network structure. \cite{barz2016air} deal with the random consumption in a similar way as DLP by taking its expected value in the formulation. To ease the exposition, we write down the decomposed formula when the cargo volume is always $V_i=0$ (one-dimensional capacity).
\begin{align*}\nonumber
& H_{t,j}(w)=\max_{x\in \{0,1\}^d}\sum_{i:a_{j,i}=1} \lambda_{t,i}\{\hat{r}_{i,j}x_i+H_{t-1,j}(w+x_i\mathbb{E}[W_i])\}, ~~\forall 1\leq t\leq T\\
&H_{0,j}(w)=-\hat{l}_{j}\mathbb{E}_{c_j}[(w-c_j)^+].
\end{align*}
In this formulation, one needs to approximate the penalty $\{\hat{l}_j\}_{j=1}^m$ of rejecting one unit weight from the real loss $l=(l_1,l_2,\dots,l_d)^T$, as well as the revenue $\hat{r}_{i,j}$. Theorem 4 from \cite{barz2016air} states that as long as $\sum_{j=1}^m a_{j,i}\hat{r}_{i,j}=r_i$, and $\sum_{j=1}^m a_{j,i}\hat{l}_j\leq l_i$ hold for all reservation class $i\in[d]$, the decomposed model gives an upper bound on the maximum expected revenue. Let $b_j$ be the shadow price of capacity constraint of inventory class $j$. \textcolor{black}{ \citet{barz2016air} suggest using $\hat{r}_{i,j}=r_i\frac{b_j}{\sum_{j=1}^m a_{j,i}b_j}$ as more revenue should be allocated to legs with positive bid prices, i.e., the capacity is tight. Such intuition is similar to what most DPD methods use. Similarly, the loss is set as $\hat{l}_j=\min_{i:a_{j,1}=1}l_i\frac{b_j}{\sum_{j'=1}^m a_{j',i}b_{j'}}$.
}

{\color{black}
It is worth mentioning that ACDPD deals with routing decisions in a heuristic way. During the reservation stage, ACDPD splits requests from each demand class with specified origin-destination pair, but non-designated routes equally into multiple demand subclasses, which have the same origin-destination pair but different designated routes.  
}

\noindent \textbf{Virtual Capacity and Bid Price Policy by \cite{previgliano2021managing} (VCBP).} This benchmark strategy is designed specifically for solving passenger NRM problems. In VCBP control, the airline sets a virtual capacity and a bid-price for each leg and accepts an incoming request if revenue is not less than the sum of bid prices of used inventories and there is sufficient virtual capacities. VCBP consider two different random capacity settings, Resource Allocation (RA) and the Random Capacity (RC) in a unified framework. They formulate the problem as the stochastic optimization model and develop a stochastic gradient-based algorithm, which guarantees the stationary convergence. In the RC setting which is similar to ours, they allow the random capacity to be revealed at any time during the reservation stage rather than at the beginning of the service stage as we assumed. Our method can be easily adapted into this setting via resolving at the capacity revealed time. In addition, our method can also be adapted to their RA setting where the decision-maker has to assign $m$ available resources with realized capacity level to $m$ inventory classes (i.e., make a scheduling decision to allocate $m$ air crafts to serve $m$ different legs) by incorporating the resource allocation decision. 

In our implementation, the stopping criteria of their stochastic gradient-based algorithm for VCBP is set the same as \cite{previgliano2021managing}, which stops at the 2,500-th iteration. Because their stochastic gradient-based algorithm only guarantees the convergence to stationary points, the convergent solution varies with different samples. We implement their algorithm five times in every passenger network instance and report the best one for comparison.

\noindent \textbf{Revenue Evaluation} The expected revenue of control policies is evaluated via 5,000 independent Monte Carlo samples. Since we do not allow fractional acceptance for a fair comparison, booking limits are rounded to the nearest integer value when calculating the expected revenue. Although, in theory, our booking limit model assumes independent demands, to be consistent with VCBP and DPD, we set the random demands among different classes to be slightly negatively correlated due to the multinomial distribution. This negative correlation is extremely small (the average coefficient of correlation among all instances is $-0.0032$) and can be ignored.

\subsection{Discussions on Assumptions in Numerical Studies}
\label{appendix:numerical_assumption_fails}

In the experiments, the following assumptions might be lacking. Even so, the numerical results show that the proposed methods still achieve superior performance against the benchmarks in nearly all test instances. See the comparison in Table \ref{tab:s-h-rand-cap}. 
\begin{itemize}
    \item Assumption on that $x$ should take integer values in passenger NRM by nature. In this work, we consider continuous decision variables, and we perform rounding to integer numbers after obtaining the optimal continuous solution. One may also use continuous $x$ and do a randomized booking limit policy.
    \item Assumption \ref{assumption:general}(a), the compact domain $\mathcal{X}$. It fails as in NRM settings, the domain is just $x\in\RR^d_+$. However, we can always manually add an upper bound $\bar X$ in the NRM setting ~\citep{karaesmen2004overbooking}.
    \item Assumption \ref{assumption:general}(b), the convexity of $f$. When there is random no-shows or cancellation in NRM problems, $f$ is only component-wise convex.
    \item Assumption \ref{assumption:general_2}(b), $\PP_i(\xi\geq\bar X_i)\geq\mu_g$ for all $i\in[d]$. It can fail when the upper bound on $\mathcal{X}$ is not carefully chosen. See discussions in Appendix \ref{appendix:a2.3a_fails}.
    \item Assumption \ref{assumption:general_2}(c) requires that the CDF of $\PP(\xi)$ is $S_g$-Lipschitz continuous. It fails as $\xi$ is a discrete distributed random vector in our numerical studies. We discuss the convergence of MSG under such a setting in Theorem \ref{cor:msg_without_smoothness}.
\end{itemize}

\subsection{Numerical Convergence Comparison of RSG, MSG, and SAA+SG on a Passenger NRM Instance}\label{appendix:convergence_numerical_NRM}
In this appendix, we compare the convergence behavior of RSG, MSG, SAA+SG, and SG (RSG without regularization, i.e., $\lambda = 0$) through a passenger NRM instance $(4,4,4,0,1,1.2,0.1)$, where the show-up probability is $1$ so that $f$ in \eqref{sec:air_cargo_model} is concave as required. 
Note that we assume continuous decision space and allow fractional acceptance. All three algorithms initialize at $x=0$. Furthermore, we evaluate the performance of the solutions from these algorithms using sample average estimation of the objective value over 5,000 independent samples at every 50 iterations and stop all algorithms at the 3,000-th iteration. See more detailed parameter choice of algorithms in Appendix \ref{appendix:computation_details}. We use this test instance to show: 1) RSG, MSG, and SAA+SG converge to the same optimal solution, 2) gradient complexities, 3) the benefits of using a regularization with $\lambda>0$ in the computation.

Figure \ref{fig:comparison_of_algorithms} (Left) demonstrates the convergence in terms of the objective value where each line represents the average revenue (objective) gained by each algorithm, $x$-axis is the index of iteration, and $y$-axis represents expected revenue. In addition, we report the booking limit solution obtained by different algorithms at 3,000-th iteration in Figure \ref{fig:comparison_of_algorithms} (Right), where $x$-axis is the index of 40 demand classes, i.e., the index $i$ in booking limit $x_i$ and $y$-axis is the value of the solution (solution $x$ is rounded to the nearest integer value and truncated at 100 for better illustration). We verify that RSG, MSG, and SAA+SG all converge to the same objective value as indicated by Figure \ref{fig:comparison_of_algorithms} (Left) and the same solution as indicated by Figure \ref{fig:comparison_of_algorithms} (Right). Figure \ref{fig:comparison_of_algorithms} (Left) additionally shows that SAA+SG indeed converges slower than RSG and MSG, as we mentioned in the introduction.

An interesting observation from Figure \ref{fig:comparison_of_algorithms} (Right) is that SG, i.e., RSG without using regularization, has extremely slow convergence. In fact, it fails to converge to an approximate optimal solution after 3000-th iterations since the revenue achieved is much smaller than the revenue achieved by MSG. 
As we have mentioned in Example \ref{example:RSG_on_truncation}, when $\lambda = 0$, SG would update the $i$-th coordinate of the decision variable $x$ only when the event $\{x_i\leq \xi_i\}$ happens. In our test experiment instance, some components of  $x$ could arrive at a very large value as shown in Figure  \ref{fig:comparison_of_algorithms} (Right) (we truncated the booking limit with a value larger than 100 to 100 in this figure for better illustration). As a result, SG encounters a vanishing gradient issue, and could take a long time to update these coordinates. If one knows the upper bound of the support of the random variable well enough, one can choose a small initialization point and a small stepsize to avoid encountering a decision point with large values. However, a small stepsize would also lead to a slow convergence speed. Although SG and RSG have the same gradient complexity $\cO(\eps^{-4})$ under Assumption \ref{assumption:general} and \ref{assumption:general_2} from a theoretical perspective, Figure \ref{fig:comparison_of_algorithms} demonstrates the importance of adding a regularization in RSG from a practical perspective.

\begin{figure}
\centering
\begin{minipage}{0.47\linewidth}
\begin{center}
\includegraphics[width=0.8\linewidth]{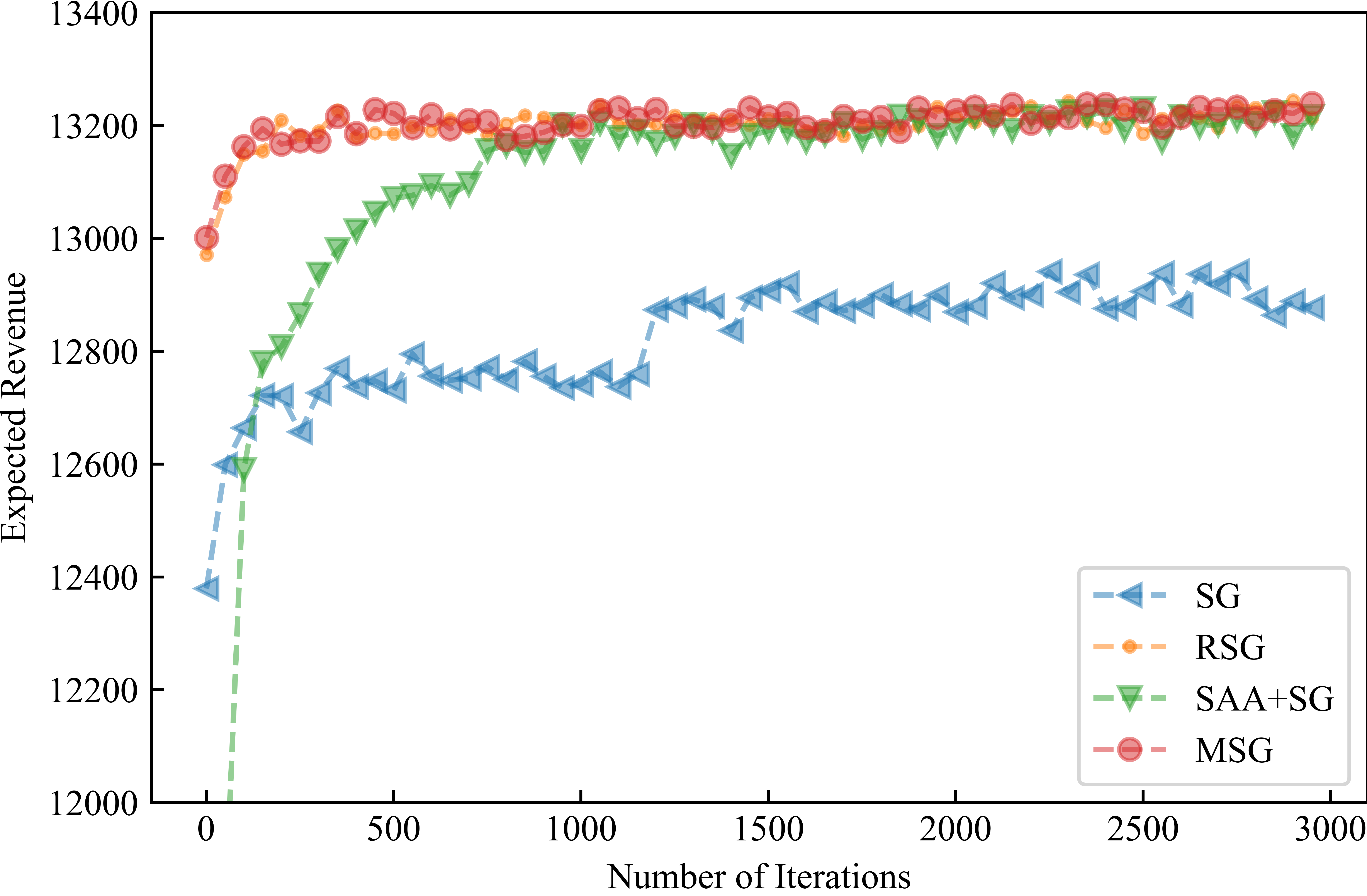}
\label{fig:compare_rev_converge}   
\end{center}
\end{minipage}
~~~~~
\begin{minipage}{0.47\linewidth}
\begin{center}
\includegraphics[width=0.8\linewidth]{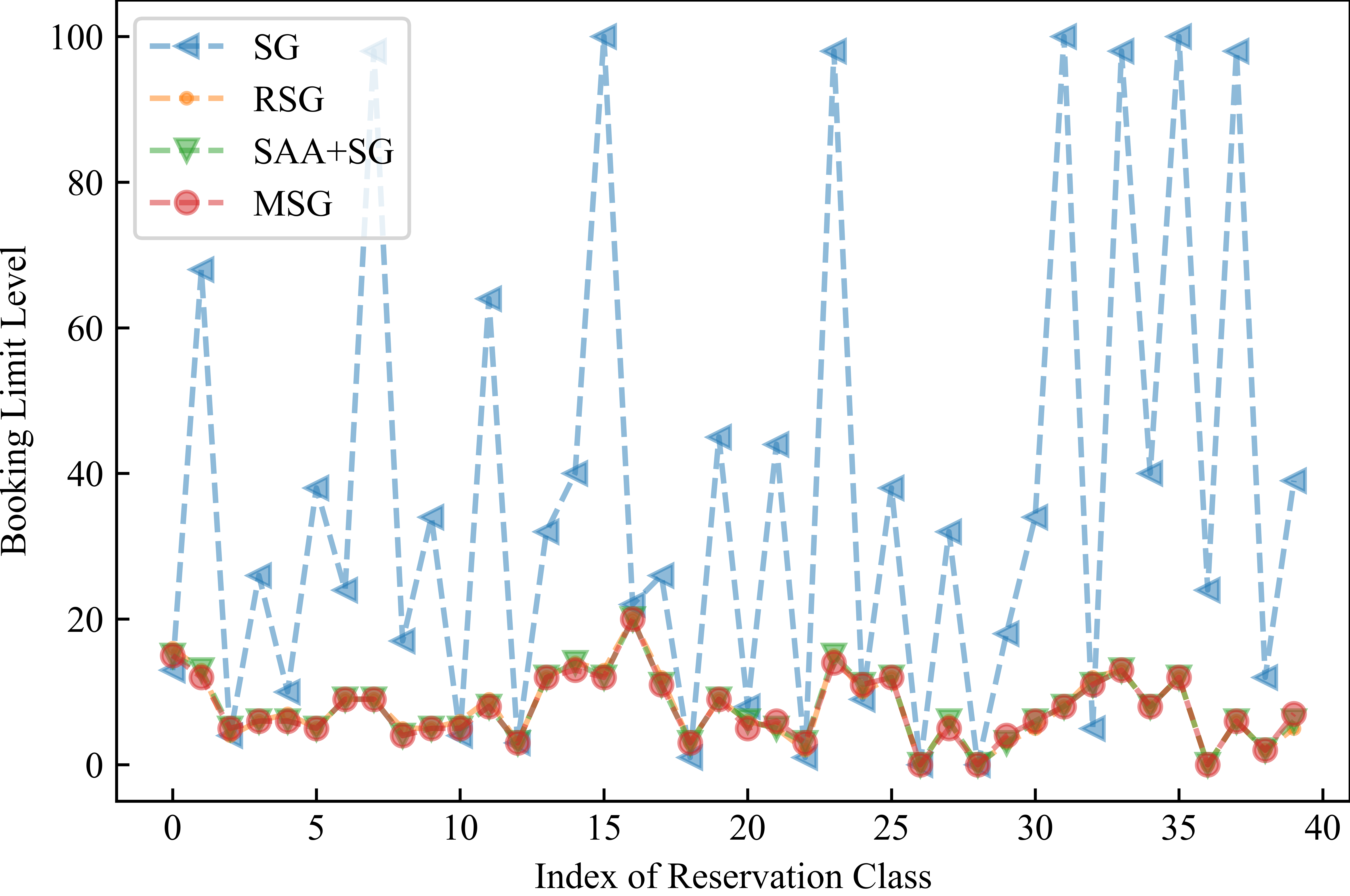}
\label{fig:compare_sol_converge}    
\end{center}
\end{minipage}
\caption{Revenue Convergence (Left) and Booking limit Solution (Right) Comparison by Different Algorithms under the instance (4,4,4,0,1,1.2,0.1).}
\label{fig:comparison_of_algorithms}
\vskip -0.2in
\end{figure}

\subsection{Complete Numerical Results in Passenger NRM}\label{appendix:complete_numerical_passenger_NRM}
Table \ref{tab:s-h-rand-cap} summarizes the complete numerical results in passenger NRM.
The first column in Table \ref{tab:s-h-rand-cap} is the parameter setting of the test instance. The second to the seventh columns give the expected revenue obtained by DLP, DPD, VCBP RSG, MSG, and SAA+SG. The remaining columns are the percentage of improvements in the expected revenue achieved by MSG over other methods. Note that $\odot$ means there is no significant difference at 95\% confidence level.

\begin{table}[ht]
\centering 
\scriptsize
\caption{Computation Results of Expected Revenue for Passenger NRM} \label{tab:s-h-rand-cap} 
\begin{tabular}{lrrrrrrrrrrrr}
%\begin{tabular}{lllllllllllll}
\\[-1.8ex]\hline 
\hline\\[-1.8ex]
\multirow{2}{*}{parameters} & \multirow{2}{*}{DLP} & \multirow{2}{*}{DPD} & \multirow{2}{*}{VCBP} & \multirow{2}{*}{RSG} & \multirow{2}{*}{MSG} & \multirow{2}{*}{SAA+SG} & \multirow{2}{*}{} & \multicolumn{5}{l}{Percentage of Improvements of MSG over}      \\ \cline{9-13} 
                            &                      &                      &                       &                      &                      &                         &                   & DLP    & DPD   & VCBP   & RSG & SAA+SG \\ 

\\[-1.8ex] \hline\\[-1.8ex]
(4,4,4,0,0.90,1.2,0.1) & 11,599 & 13,911 & 12,688 & 13,857 & 13,859 & 13,857 &  & 19.5\% & $\odot$ & 9.2\%   & $\odot$ & $\odot$ \\
(4,4,4,0,0.90,1.6,0.1) & 17,006 & 19,012 & 17,990 & 18,965 & 18,983 & 18,949 &  & 11.5\% & $\odot$ & 5.5\%   & $\odot$ & $\odot$ \\
(4,4,4,0,0.95,1.2,0.1) & 19,507 & 21,318 & 20,705 & 21,458 & 21,477 & 21,466 &  & 10.0\% & $\odot$ & 3.7\%   & $\odot$ & $\odot$ \\
(4,4,4,0,0.95,1.6,0.1) & 16,367 & 19,757 & 17,179 & 19,642 & 19,682 & 19,641 &  & 20.0\% & $\odot$ & 14.6\%  & $\odot$ & $\odot$ \\
\\[-1.8ex] \hline\\[-1.8ex]
(4,4,8,0,0.90,1.2,0.1) & 19,174 & 20,167 & 20,628 & 20,666 & 20,629 & 20,616 &  & 7.8\%  & 2.3\%   & $\odot$ & $\odot$ & $\odot$ \\
(4,4,8,0,0.90,1.6,0.1) & 8,792  & 10,372 & 9,418  & 10,598 & 10,600 & 10,581 &  & 20.5\% & 2.2\%   & 12.5\%  & $\odot$ & $\odot$ \\
(4,4,8,0,0.95,1.2,0.1) & 14,066 & 16,491 & 16,253 & 16,939 & 16,955 & 16,979 &  & 20.4\% & 2.8\%   & 4.3\%   & $\odot$ & $\odot$ \\
(4,4,8,0,0.95,1.6,0.1) & 12,396 & 13,612 & 13,422 & 13,893 & 13,896 & 13,890 &  & 12.1\% & 2.1\%   & 3.5\%   & $\odot$ & $\odot$ \\
\\[-1.8ex] \hline\\[-1.8ex]
(4,4,1,1,0.90,1.2,0.1) & 12,759 & 15,946 & 15,575 & 16,288 & 16,282 & 16,264 &  & 27.7\% & 2.1\%   & 4.5\%   & $\odot$ & $\odot$ \\
(4,4,1,1,0.90,1.6,0.1) & 8,738  & 9,401  & 9,261  & 9,473  & 9,502  & 9,491  &  & 8.4\%  & $\odot$ & 2.6\%   & $\odot$ & $\odot$ \\
(4,4,1,1,0.95,1.2,0.1) & 16,320 & 17,228 & 17,613 & 17,541 & 17,568 & 17,530 &  & 7.5\%  & 2.0\%   & $\odot$ & $\odot$ & $\odot$ \\
(4,4,1,1,0.95,1.6,0.1) & 13,754 & 15,182 & 14,306 & 15,274 & 15,311 & 15,273 &  & 11.1\% & $\odot$ & 7.0\%   & $\odot$ & $\odot$ \\
\\[-1.8ex] \hline\\[-1.8ex]
(4,8,4,0,0.90,1.2,0.1) & 25,631 & 29,101 & 28,352 & 29,034 & 29,010 & 29,069 &  & 13.3\% & $\odot$ & 2.3\%   & $\odot$ & $\odot$ \\
(4,8,4,0,0.90,1.6,0.1) & 27,229 & 30,110 & 29,456 & 30,020 & 30,114 & 30,013 &  & 10.3\% & $\odot$ & 2.2\%   & $\odot$ & $\odot$ \\
(4,8,4,0,0.95,1.2,0.1) & 25,845 & 27,403 & 26,988 & 27,348 & 27,373 & 27,312 &  & 5.8\%  & $\odot$ & 1.4\%   & $\odot$ & $\odot$ \\
(4,8,4,0,0.95,1.6,0.1) & 21,985 & 22,741 & 22,541 & 22,713 & 22,741 & 22,711 &  & 3.3\%  & $\odot$ & 0.9\%   & $\odot$ & $\odot$ \\
\\[-1.8ex] \hline\\[-1.8ex]
(4,8,8,0,0.90,1.2,0.1) & 23,354 & 24,722 & 23,887 & 24,973 & 24,920 & 24,947 &  & 6.9\%  & $\odot$ & 4.3\%   & $\odot$ & $\odot$ \\
(4,8,8,0,0.90,1.6,0.1) & 26,401 & 28,050 & 27,239 & 28,287 & 28,322 & 28,257 &  & 7.1\%  & $\odot$ & 4.0\%   & $\odot$ & $\odot$ \\
(4,8,8,0,0.95,1.2,0.1) & 28,555 & 30,299 & 30,146 & 30,807 & 30,802 & 30,806 &  & 7.9\%  & 1.7\%   & 2.2\%   & $\odot$ & $\odot$ \\
(4,8,8,0,0.95,1.6,0.1) & 23,915 & 25,685 & 25,241 & 25,838 & 25,813 & 25,774 &  & 8.0\%  & $\odot$ & 2.3\%   & $\odot$ & $\odot$ \\
\\[-1.8ex] \hline\\[-1.8ex]
(4,8,1,1,0.90,1.2,0.1) & 23,128 & 24,052 & 23,548 & 24,501 & 24,536 & 24,468 &  & 5.9\%  & 2.0\%   & 4.2\%   & $\odot$ & $\odot$ \\
(4,8,1,1,0.90,1.6,0.1) & 26,435 & 28,559 & 27,521 & 28,825 & 28,823 & 28,809 &  & 9.0\%  & $\odot$ & 4.7\%   & $\odot$ & $\odot$ \\
(4,8,1,1,0.95,1.2,0.1) & 18,365 & 19,543 & 20,218 & 20,168 & 20,163 & 20,137 &  & 9.8\%  & 3.2\%   & $\odot$ & $\odot$ & $\odot$ \\
(4,8,1,1,0.95,1.6,0.1) & 21,704 & 23,490 & 22,525 & 23,840 & 23,816 & 23,720 &  & 9.8\%  & 1.4\%   & 5.7\%   & $\odot$ & $\odot$ \\
\\[-1.8ex] \hline\\[-1.8ex]
(8,4,4,0,0.90,1.2,0.1) & 18,845 & 19,896 & 18,725 & 19,848 & 19,822 & 19,849 &  & 5.3\%  & $\odot$ & 5.9\%   & $\odot$ & $\odot$ \\
(8,4,4,0,0.90,1.6,0.1) & 13,441 & 15,335 & 14,545 & 15,223 & 15,260 & 15,242 &  & 13.3\% & $\odot$ & 4.9\%   & $\odot$ & $\odot$ \\
(8,4,4,0,0.95,1.2,0.1) & 17,092 & 19,328 & 18,287 & 19,169 & 19,200 & 19,100 &  & 12.1\% & $\odot$ & 5.0\%   & $\odot$ & $\odot$ \\
(8,4,4,0,0.95,1.6,0.1) & 12,537 & 13,559 & 12,987 & 13,465 & 13,456 & 13,387 &  & 7.4\%  & $\odot$ & 3.6\%   & $\odot$ & $\odot$ \\
\\[-1.8ex] \hline\\[-1.8ex]
(8,4,8,0,0.90,1.2,0.1) & 15,197 & 16,525 & 15,487 & 16,847 & 16,784 & 16,713 &  & 10.9\% & 1.6\%   & 8.4\%   & $\odot$ & $\odot$ \\
(8,4,8,0,0.90,1.6,0.1) & 10,556 & 11,402 & 10,866 & 11,408 & 11,386 & 11,358 &  & 8.1\%  & $\odot$ & 4.8\%   & $\odot$ & $\odot$ \\
(8,4,8,0,0.95,1.2,0.1) & 15,912 & 17,698 & 16,615 & 17,941 & 17,997 & 17,949 &  & 12.7\% & 1.7\%   & 8.3\%   & $\odot$ & $\odot$ \\
(8,4,8,0,0.95,1.6,0.1) & 12,839 & 14,118 & 13,375 & 14,054 & 14,074 & 14,047 &  & 9.5\%  & $\odot$ & 5.2\%   & $\odot$ & $\odot$ \\
\\[-1.8ex] \hline\\[-1.8ex]
(8,4,1,1,0.90,1.2,0.1) & 15,388 & 17,141 & 16,410 & 17,144 & 17,186 & 17,131 &  & 11.4\% & $\odot$ & 4.7\%   & $\odot$ & $\odot$ \\
(8,4,1,1,0.90,1.6,0.1) & 12,082 & 12,916 & 12,226 & 12,931 & 12,950 & 12,959 &  & 7.0\%  & $\odot$ & 5.9\%   & $\odot$ & $\odot$ \\
(8,4,1,1,0.95,1.2,0.1) & 11,265 & 13,516 & 12,967 & 13,643 & 13,646 & 13,603 &  & 21.1\% & $\odot$ & 5.2\%   & $\odot$ & $\odot$ \\
(8,4,1,1,0.95,1.6,0.1) & 14,912 & 16,504 & 16,030 & 16,327 & 16,399 & 16,340 &  & 9.5\%  & $\odot$ & 2.3\%   & $\odot$ & $\odot$ \\
\\[-1.8ex] \hline\\[-1.8ex]
(8,8,4,0,0.90,1.2,0.1) & 31,329 & 32,411 & 32,283 & 32,392 & 32,395 & 32,322 &  & 3.4\%  & $\odot$ & $\odot$ & $\odot$ & $\odot$ \\
(8,8,4,0,0.90,1.6,0.1) & 20,851 & 22,197 & 21,731 & 22,158 & 22,181 & 22,126 &  & 6.3\%  & $\odot$ & 2.1\%   & $\odot$ & $\odot$ \\
(8,8,4,0,0.95,1.2,0.1) & 30,770 & 31,370 & 30,714 & 31,349 & 31,270 & 31,286 &  & 1.9\%  & $\odot$ & 1.8\%   & $\odot$ & $\odot$ \\
(8,8,4,0,0.95,1.6,0.1) & 22,556 & 23,605 & 23,286 & 23,653 & 23,586 & 23,509 &  & 4.9\%  & $\odot$ & 1.3\%   & $\odot$ & $\odot$ \\
\\[-1.8ex] \hline\\[-1.8ex]
(8,8,8,0,0.90,1.2,0.1) & 26,975 & 28,331 & 27,610 & 28,821 & 28,799 & 28,797 &  & 6.8\%  & 1.7\%   & 4.3\%   & $\odot$ & $\odot$ \\
(8,8,8,0,0.90,1.6,0.1) & 31,464 & 33,425 & 32,749 & 33,579 & 33,457 & 33,320 &  & 6.7\%  & $\odot$ & 2.2\%   & $\odot$ & $\odot$ \\
(8,8,8,0,0.95,1.2,0.1) & 24,107 & 25,812 & 25,421 & 25,847 & 25,837 & 25,793 &  & 7.2\%  & $\odot$ & 1.6\%   & $\odot$ & $\odot$ \\
(8,8,8,0,0.95,1.6,0.1) & 27,140 & 29,291 & 28,923 & 29,184 & 29,263 & 29,225 &  & 7.5\%  & $\odot$ & 1.2\%   & $\odot$ & $\odot$ \\
\\[-1.8ex] \hline\\[-1.8ex]
(8,8,1,1,0.90,1.2,0.1) & 22,200 & 25,310 & 25,004 & 25,963 & 25,939 & 25,997 &  & 16.9\% & 2.5\%   & 3.7\%   & $\odot$ & $\odot$ \\
(8,8,1,1,0.90,1.6,0.1) & 23,612 & 24,547 & 24,365 & 24,706 & 24,698 & 24,653 &  & 4.6\%  & $\odot$ & 1.4\%   & $\odot$ & $\odot$ \\
(8,8,1,1,0.95,1.2,0.1) & 19,993 & 23,486 & 23,326 & 23,726 & 23,733 & 23,700 &  & 18.7\% & 1.1\%   & 1.7\%   & $\odot$ & $\odot$ \\
(8,8,1,1,0.95,1.6,0.1) & 21,875 & 24,691 & 24,022 & 24,704 & 24,581 & 24,498 &  & 12.9\% & $\odot$ & 2.3\%   & $\odot$ & $\odot$\\
\hline 
\hline \\[-1.8ex] 
\end{tabular}
\begin{tablenotes}
\item $\odot$ denotes there is no statistically significant difference between MSG and the alternative, all at 95\% confidence level. All other comparisons are significant at 95\% confidence level.
\end{tablenotes}
\end{table}

\begin{table}[ht]
\centering \ContinuedFloat
\caption{Computation Results for Passenger NRM (Continued)} %\label{tab:s-h-rand-cap-2}
\scriptsize
\begin{tabular}{lrrrrrrrrrrrr}
%\begin{tabular}{lllllllllllll}
\\[-1.8ex]\hline 
\hline\\[-1.8ex]
\multirow{2}{*}{parameters} & \multirow{2}{*}{DLP} & \multirow{2}{*}{DPD} & \multirow{2}{*}{VCBP} & \multirow{2}{*}{RSG} & \multirow{2}{*}{MSG} & \multirow{2}{*}{SAA+SG} & \multirow{2}{*}{} & \multicolumn{5}{l}{Percentage of Improvements of MSG over}      \\ \cline{9-13} 
                            &                      &                      &                       &                      &                      &                         &                   & DLP    & DPD   & VCBP   & RSG & SAA+SG \\ 

\\[-1.8ex] \hline\\[-1.8ex]
(4,4,4,0,0.90,1.2,0.5) & 8,518  & 9,133  & 9,253  & 9,642  & 9,669  & 9,634  &  & 13.2\%  & 5.9\%   & 4.5\%   & $\odot$ & $\odot$ \\
(4,4,4,0,0.90,1.6,0.5) & 11,884 & 11,971 & 11,996 & 12,330 & 12,332 & 12,307 &  & 3.8\%   & 3.0\%   & 2.8\%   & $\odot$ & $\odot$ \\
(4,4,4,0,0.95,1.2,0.5) & 13,236 & 13,883 & 14,824 & 14,821 & 14,882 & 14,897 &  & 12.0\%  & 7.2\%   & $\odot$ & $\odot$ & $\odot$ \\
(4,4,4,0,0.95,1.6,0.5) & 12,466 & 13,227 & 12,525 & 13,272 & 13,335 & 13,392 &  & 6.5\%   & $\odot$ & 6.5\%   & $\odot$ & $\odot$ \\
\\[-1.8ex] \hline\\[-1.8ex]
(4,4,8,0,0.90,1.2,0.5) & 3,591  & 9,061  & 10,630 & 11,689 & 11,670 & 11,688 &  & 225.5\% & 28.8\%  & 9.8\%   & $\odot$ & $\odot$ \\
(4,4,8,0,0.90,1.6,0.5) & 4,960  & 5,008  & 5,333  & 5,434  & 5,451  & 5,440  &  & 9.6\%   & 8.8\%   & 2.2\%   & $\odot$ & $\odot$ \\
(4,4,8,0,0.95,1.2,0.5) & 4,024  & 7,588  & 8,186  & 8,499  & 8,541  & 8,522  &  & 111.2\% & 12.6\%  & 4.3\%   & $\odot$ & $\odot$ \\
(4,4,8,0,0.95,1.6,0.5) & 3,249  & 6,419  & 6,402  & 6,752  & 6,800  & 6,783  &  & 107.8\% & 5.9\%   & 6.2\%   & $\odot$ & $\odot$ \\
\\[-1.8ex] \hline\\[-1.8ex]
(4,4,1,1,0.90,1.2,0.5) & 9,750  & 10,485 & 10,811 & 10,966 & 10,993 & 10,967 &  & 12.5\%  & 4.9\%   & 1.7\%   & $\odot$ & $\odot$ \\
(4,4,1,1,0.90,1.6,0.5) & 6,217  & 6,667  & 6,462  & 6,742  & 6,719  & 6,753  &  & 8.4\%   & $\odot$ & 4.0\%   & $\odot$ & $\odot$ \\
(4,4,1,1,0.95,1.2,0.5) & 7,506  & 10,527 & 11,224 & 11,320 & 11,355 & 11,318 &  & 50.8\%  & 7.9\%   & 1.2\%   & $\odot$ & $\odot$ \\
(4,4,1,1,0.95,1.6,0.5) & 10,164 & 10,590 & 10,380 & 10,854 & 10,842 & 10,893 &  & 6.8\%   & 2.4\%   & 4.4\%   & $\odot$ & $\odot$ \\
\\[-1.8ex] \hline\\[-1.8ex]
(4,8,4,0,0.90,1.2,0.5) & 20,455 & 21,474 & 20,726 & 21,888 & 21,785 & 21,747 &  & 7.0\%   & 1.4\%   & 5.1\%   & $\odot$ & $\odot$ \\
(4,8,4,0,0.90,1.6,0.5) & 20,880 & 21,665 & 19,806 & 21,520 & 21,558 & 21,664 &  & 3.1\%   & $\odot$ & 8.8\%   & $\odot$ & $\odot$ \\
(4,8,4,0,0.95,1.2,0.5) & 20,073 & 20,271 & 20,449 & 20,698 & 20,807 & 20,793 &  & 3.1\%   & 2.6\%   & 1.7\%   & $\odot$ & $\odot$ \\
(4,8,4,0,0.95,1.6,0.5) & 16,189 & 16,206 & 15,822 & 16,508 & 16,579 & 16,551 &  & 2.0\%   & 2.3\%   & 4.8\%   & $\odot$ & $\odot$ \\
\\[-1.8ex] \hline\\[-1.8ex]
(4,8,8,0,0.90,1.2,0.5) & 14,504 & 15,316 & 14,519 & 15,655 & 15,646 & 15,611 &  & 7.9\%   & 2.2\%   & 7.8\%   & $\odot$ & $\odot$ \\
(4,8,8,0,0.90,1.6,0.5) & 11,908 & 14,397 & 12,576 & 15,207 & 15,260 & 15,267 &  & 27.7\%  & 6.0\%   & 21.3\%  & $\odot$ & $\odot$ \\
(4,8,8,0,0.95,1.2,0.5) & 11,287 & 17,260 & 17,371 & 17,854 & 17,929 & 17,917 &  & 58.2\%  & 3.9\%   & 3.2\%   & $\odot$ & $\odot$ \\
(4,8,8,0,0.95,1.6,0.5) & 11,938 & 12,617 & 12,705 & 12,900 & 12,950 & 12,943 &  & 8.1\%   & 2.6\%   & 1.9\%   & $\odot$ & $\odot$ \\
\\[-1.8ex] \hline\\[-1.8ex]
(4,8,1,1,0.90,1.2,0.5) & 12,413 & 17,799 & 18,686 & 18,558 & 18,667 & 18,583 &  & 49.5\%  & 4.9\%   & $\odot$ & $\odot$ & $\odot$ \\
(4,8,1,1,0.90,1.6,0.5) & 15,850 & 22,496 & 22,039 & 22,651 & 22,587 & 22,596 &  & 42.9\%  & $\odot$ & 2.5\%   & $\odot$ & $\odot$ \\
(4,8,1,1,0.95,1.2,0.5) & 9,690  & 14,897 & 15,053 & 15,142 & 15,206 & 15,195 &  & 56.3\%  & 2.1\%   & 1.0\%   & $\odot$ & $\odot$ \\
(4,8,1,1,0.95,1.6,0.5) & 14,056 & 18,118 & 18,242 & 18,433 & 18,539 & 18,577 &  & 31.1\%  & 2.3\%   & 1.6\%   & $\odot$ & $\odot$ \\
\\[-1.8ex] \hline\\[-1.8ex]
(8,4,4,0,0.90,1.2,0.5) & 11,029 & 11,350 & 12,136 & 13,171 & 13,148 & 13,162 &  & 19.4\%  & 15.8\%  & 8.3\%   & $\odot$ & $\odot$ \\
(8,4,4,0,0.90,1.6,0.5) & 9,033  & 8,942  & 8,718  & 9,684  & 9,764  & 9,760  &  & 7.2\%   & 9.2\%   & 12.0\%  & $\odot$ & $\odot$ \\
(8,4,4,0,0.95,1.2,0.5) & 11,194 & 11,275 & 12,538 & 12,720 & 12,746 & 12,703 &  & 13.6\%  & 13.0\%  & 1.7\%   & $\odot$ & $\odot$ \\
(8,4,4,0,0.95,1.6,0.5) & 7,120  & 7,797  & 8,227  & 8,543  & 8,561  & 8,565  &  & 20.0\%  & 9.8\%   & 4.0\%   & $\odot$ & $\odot$ \\
\\[-1.8ex] \hline\\[-1.8ex]
(8,4,8,0,0.90,1.2,0.5) & 2,020  & 4,629  & 7,266  & 8,439  & 8,473  & 8,460  &  & 317.8\% & 83.1\%  & 16.6\%  & $\odot$ & $\odot$ \\
(8,4,8,0,0.90,1.6,0.5) & 3,357  & 3,870  & 4,600  & 5,436  & 5,419  & 5,449  &  & 61.9\%  & 40.0\%  & 17.8\%  & $\odot$ & $\odot$ \\
(8,4,8,0,0.95,1.2,0.5) & 2,398  & 4,528  & 8,380  & 9,002  & 8,931  & 8,964  &  & 275.4\% & 97.2\%  & 6.6\%   & $\odot$ & $\odot$ \\
(8,4,8,0,0.95,1.6,0.5) & 2,914  & 4,119  & 5,808  & 6,378  & 6,368  & 6,333  &  & 118.9\% & 54.6\%  & 9.6\%   & $\odot$ & $\odot$ \\
\\[-1.8ex] \hline\\[-1.8ex]
(8,4,1,1,0.90,1.2,0.5) & 6,633  & 9,694  & 11,109 & 11,644 & 11,693 & 11,617 &  & 75.5\%  & 20.6\%  & 5.3\%   & $\odot$ & $\odot$ \\
(8,4,1,1,0.90,1.6,0.5) & 6,687  & 7,375  & 8,519  & 8,735  & 8,782  & 8,780  &  & 30.6\%  & 19.1\%  & 3.1\%   & $\odot$ & $\odot$ \\
(8,4,1,1,0.95,1.2,0.5) & 2,705  & 6,319  & 8,799  & 8,794  & 8,879  & 8,899  &  & 225.1\% & 40.5\%  & 0.9\%   & $\odot$ & $\odot$ \\
(8,4,1,1,0.95,1.6,0.5) & 7,130  & 9,097  & 10,134 & 11,014 & 11,065 & 11,067 &  & 54.5\%  & 21.6\%  & 9.2\%   & $\odot$ & $\odot$ \\
\\[-1.8ex] \hline\\[-1.8ex]
(8,8,4,0,0.90,1.2,0.5) & 21,514 & 22,051 & 22,295 & 23,359 & 23,390 & 23,395 &  & 8.6\%   & 6.1\%   & 4.9\%   & $\odot$ & $\odot$ \\
(8,8,4,0,0.90,1.6,0.5) & 15,074 & 14,772 & 14,315 & 15,575 & 15,516 & 15,606 &  & 3.3\%   & 5.0\%   & 8.4\%   & $\odot$ & $\odot$ \\
(8,8,4,0,0.95,1.2,0.5) & 20,492 & 21,035 & 22,462 & 22,619 & 22,591 & 22,712 &  & 10.4\%  & 7.4\%   & $\odot$ & $\odot$ & $\odot$ \\
(8,8,4,0,0.95,1.6,0.5) & 15,167 & 14,962 & 16,337 & 16,507 & 16,504 & 16,490 &  & 8.8\%   & 10.3\%  & 1.0\%   & $\odot$ & $\odot$ \\
\\[-1.8ex] \hline\\[-1.8ex]
(8,8,8,0,0.90,1.2,0.5) & 8,655  & 12,340 & 14,967 & 16,647 & 16,710 & 16,696 &  & 92.4\%  & 35.4\%  & 11.6\%  & $\odot$ & $\odot$ \\
(8,8,8,0,0.90,1.6,0.5) & 9,735  & 11,891 & 12,804 & 15,324 & 15,454 & 15,357 &  & 57.4\%  & 30.0\%  & 20.7\%  & $\odot$ & $\odot$ \\
(8,8,8,0,0.95,1.2,0.5) & 7,088  & 11,449 & 14,049 & 14,837 & 14,875 & 14,816 &  & 109.3\% & 29.9\%  & 5.9\%   & $\odot$ & $\odot$ \\
(8,8,8,0,0.95,1.6,0.5) & 7,903  & 12,001 & 13,535 & 14,287 & 14,265 & 14,293 &  & 80.8\%  & 18.9\%  & 5.4\%   & $\odot$ & $\odot$ \\
\\[-1.8ex] \hline\\[-1.8ex]
(8,8,1,1,0.90,1.2,0.5) & 2,525  & 15,561 & 19,049 & 19,049 & 19,047 & 19,073 &  & 654.5\% & 22.4\%  & $\odot$ & $\odot$ & $\odot$ \\
(8,8,1,1,0.90,1.6,0.5) & 13,177 & 15,605 & 17,254 & 18,189 & 18,153 & 18,189 &  & 38.0\%  & 16.3\%  & 5.2\%   & $\odot$ & $\odot$ \\
(8,8,1,1,0.95,1.2,0.5) & 3,319  & 14,344 & 17,265 & 17,797 & 17,700 & 17,730 &  & 436.2\% & 23.4\%  & 2.5\%   & $\odot$ & $\odot$ \\
(8,8,1,1,0.95,1.6,0.5) & 9,876  & 15,342 & 17,495 & 18,276 & 18,332 & 18,288 &  & 85.1\%  & 19.5\%  & 4.8\%   & $\odot$ & $\odot$\\
\hline 
\hline \\[-1.8ex] 
\end{tabular}
\begin{tablenotes}
\item $\odot$ denotes there is no statistically significant difference between MSG and the alternative, all at 95\% confidence level. All other comparisons are significant at 95\% confidence level.
\end{tablenotes}
\end{table}

\subsection{Parameters of Air Cargo NRM Instances}\label{appendix:air_cargo_param}

Since the full information regarding the reservation classes in Appendix of \cite{barz2016air} is truncated, we construct similar instances based on the following parameters listed in Table \ref{tab:example_air_cargo}.

\begin{table}[ht]
\centering 
\scriptsize
\caption{Parameters of Air Cargo NRM Instances} \label{tab:example_air_cargo}
\begin{tabular}{rrrrrr}
\\[-1.8ex]\hline 
\hline\\[-1.8ex]
Class & \begin{tabular}[c]{@{}l@{}}Mean\\ weight\end{tabular} & \begin{tabular}[c]{@{}l@{}}Mean \\ volume\end{tabular} & Origin & Destination & \begin{tabular}[c]{@{}l@{}}Per-unit\\ revenue\end{tabular} \\ 
\hline\\[-1.8ex]
1     & 5                                                     & 3                                                      & 1      & 5           & 1.4                                                        \\
2     & 5                                                     & 4                                                      & 2      & 1           & 1.4                                                        \\
3     & 5                                                     & 5                                                      & 3      & 1           & 1.4                                                        \\
4     & 5                                                     & 2                                                      & 1      & 2           & 1.4                                                        \\
5     & 10                                                    & 6                                                      & 3      & 5           & 1.4                                                        \\
6     & 10                                                    & 5                                                      & 2      & 3           & 1.4                                                        \\
7     & 10                                                    & 7                                                      & 1      & 3           & 1.4                                                        \\
8     & 10                                                    & 8                                                      & 3      & 4           & 1.4                                                        \\
9     & 20                                                    & 13                                                     & 3      & 1           & 1.4                                                        \\
10    & 20                                                    & 12                                                     & 2      & 3           & 1.4                                                        \\
11    & 20                                                    & 10                                                     & 2      & 1           & 1.4                                                        \\
12    & 25                                                    & 15                                                     & 5      & 4           & 1.4                                                        \\
13    & 25                                                    & 14                                                     & 1      & 5           & 1.4                                                        \\
14    & 30                                                    & 18                                                     & 4      & 3           & 1.4                                                        \\
15    & 40                                                    & 24                                                     & 2      & 4           & 1.4                                                        \\
16    & 50                                                    & 28                                                     & 1      & 3           & 1.4                                                        \\
17    & 100                                                   & 60                                                     & 2      & 1           & 1.4                                                        \\
18    & 150                                                   & 90                                                     & 5      & 4           & 1.4                                                        \\
19    & 250                                                   & 149                                                    & 1      & 2           & 1.4                                                        \\
20    & 350                                                   & 208                                                    & 4      & 1           & 1.4                                                        \\
21    & 7                                                     & 23                                                     & 2      & 3           & 1.4                                                        \\
22    & 7                                                     & 2                                                      & 3      & 1           & 1.4                                                        \\
23    & 21                                                    & 70                                                     & 1      & 2           & 1.4                                                        \\
24    & 21                                                    & 5                                                      & 4      & 3           & 1.4                                                        \\
25    & 5                                                     & 3                                                      & 3      & 4           & 0.7                                                        \\
26    & 5                                                     & 4                                                      & 1      & 3           & 0.7                                                        \\
27    & 5                                                     & 5                                                      & 3      & 5           & 0.7                                                        \\
28    & 5                                                     & 2                                                      & 3      & 2           & 0.7                                                        \\
29    & 10                                                    & 6                                                      & 1      & 5           & 0.7                                                        \\
30    & 10                                                    & 5                                                      & 5      & 4           & 0.7                                                        \\
31    & 10                                                    & 7                                                      & 1      & 2           & 0.7                                                        \\
32    & 10                                                    & 8                                                      & 2      & 3           & 0.7                                                        \\
33    & 20                                                    & 13                                                     & 3      & 2           & 0.7                                                        \\
34    & 20                                                    & 12                                                     & 1      & 4           & 0.7                                                        \\
35    & 20                                                    & 10                                                     & 3      & 2           & 0.7                                                        \\
36    & 25                                                    & 15                                                     & 2      & 1           & 0.7                                                        \\
37    & 25                                                    & 14                                                     & 1      & 2           & 0.7                                                        \\
38    & 30                                                    & 18                                                     & 3      & 5           & 0.7                                                        \\
39    & 40                                                    & 24                                                     & 4      & 1           & 0.7                                                        \\
40    & 50                                                    & 28                                                     & 1      & 5           & 0.7                                                        \\
41    & 100                                                   & 60                                                     & 1      & 2           & 0.7                                                        \\
42    & 150                                                   & 90                                                     & 3      & 4           & 0.7                                                        \\
43    & 250                                                   & 149                                                    & 4      & 3           & 0.7                                                        \\
44    & 350                                                   & 208                                                    & 2      & 3           & 0.7                                                        \\
45    & 7                                                     & 23                                                     & 2      & 3           & 0.7                                                        \\
46    & 7                                                     & 2                                                      & 4      & 1           & 0.7                                                        \\
47    & 21                                                    & 70                                                     & 5      & 4           & 0.7                                                        \\
48    & 21                                                    & 5                                                      & 1      & 3           & 0.7                                                        \\
49    & 5                                                     & 3                                                      & 3      & 5           & 1                                                          \\
50    & 5                                                     & 4                                                      & 1      & 4           & 1                                                          \\
51    & 5                                                     & 5                                                      & 3      & 2           & 1                                                          \\
52    & 5                                                     & 2                                                      & 1      & 2           & 1                                                          \\
53    & 10                                                    & 6                                                      & 5      & 4           & 1                                                          \\
54    & 10                                                    & 5                                                      & 2      & 1           & 1                                                          \\
55    & 10                                                    & 7                                                      & 2      & 3           & 1                                                          \\
56    & 10                                                    & 8                                                      & 3      & 5           & 1                                                          \\
57    & 20                                                    & 13                                                     & 1      & 5           & 1                                                          \\
58    & 20                                                    & 12                                                     & 4      & 1           & 1                                                          \\
59    & 20                                                    & 10                                                     & 4      & 3           & 1                                                          \\
60    & 25                                                    & 15                                                     & 2      & 3           & 1                                                          \\ 
\hline 
\hline \\[-1.8ex] 
\end{tabular}
\end{table}

\end{APPENDICES}

\end{document}